\documentclass[oneside,reqno]{amsart}
\usepackage{amsmath}
\usepackage{amssymb}
\usepackage{amsfonts}
\usepackage{pspicture} 
\usepackage[normalem]{ulem} 
\usepackage{verbatim}
\usepackage{eucal} 

\newtheorem{theorem}{Theorem}[section]

\newtheorem{corollary}[theorem]{Corollary}

\newtheorem{lemma}[theorem]{Lemma}

\newtheorem{proposition}[theorem]{Proposition}

\theoremstyle{definition}

\newtheorem{definition}[theorem]{Definition}
\newtheorem{example}[theorem]{Example}
\newtheorem{notation}[theorem]{Notation}

\newtheorem{remark}[theorem]{Remark}

\newcommand\ms{\hspace{.5 em}}
\newcommand\ibid{[ibid]}
\newenvironment{mylist}
{\begin{list}{}
{
\setlength{\labelwidth}{.3truein}
\setlength{\labelsep}{.1truein}
\setlength{\leftmargin}{.6truein}
\setlength{\rightmargin}{.3truein}
}}
{\end{list}}

\newcommand\CD{\operatorname{CD}}
\newcommand\ch{\operatorname{ch}}

\newcommand\GL{\operatorname{GL}}
\newcommand\Elem{\operatorname{E}}
\newcommand\id{\operatorname{id}}
\newcommand\SL{\operatorname{SL}}
\newcommand\supp{\operatorname{supp}}

\newcommand\bq{{\mathbf{q}}}

\newcommand\bbZ{\mathbb Z}

\renewcommand\th{{\tilde h}}
\newcommand\tN{{\tilde N}}

\newcommand\al{\alpha}
\newcommand\ep{\varepsilon}
\newcommand\gm{\gamma}
\newcommand\Gm{\Gamma}
\newcommand\lm{\lambda}
\newcommand\Lm{\Lambda}
\newcommand\sg{\sigma}
\newcommand\ph{\varphi}

\newcommand\bal{\bar{\alpha}}
\newcommand\bbe{\bar{\beta}}
\newcommand\bde{\bar{\delta}}
\newcommand\bgm{\bar{\gamma}}
\newcommand\blm{\bar{\lambda}}
\newcommand\bmu{\bar{\mu}}
\newcommand\bLm{\bar{\Lambda}}

\newcommand\cA{\mathcal A}
\newcommand\cB{\mathcal B}
\newcommand\cC{\mathcal C}
\newcommand\cD{\mathcal D}
\newcommand\cE{\mathcal E}

\newcommand\cH{\mathcal H}
\newcommand\cI{\mathcal I}
\newcommand\cJ{\mathcal J}

\newcommand\cM{\mathcal M}
\newcommand\cO{\mathcal O}
\newcommand\cP{\mathcal P}
\newcommand\cQ{\mathcal Q}
\newcommand\cT{\mathcal T}
\newcommand\cU{\mathcal U}
\newcommand\cV{\mathcal V}
\newcommand\cW{\mathcal W}
\newcommand\cX{\mathcal X}
\newcommand\cZ{\mathcal Z}

\newcommand\andd{\quad \text{and} \quad}
\newcommand\barinv{\bar{\phantom{a}}}
\newcommand\Cay{\cC}
\newcommand\cyc{\circlearrowleft}
\newcommand\hermk{k}
\newcommand\Herm{\cH(\cC_{3})}
\newcommand\ig{\simeq_{ig}}
\newcommand\ot{\otimes}
\newcommand\orr{\quad \text{or} \quad}
\newcommand\ord[1]{\left\vert#1 \right\vert}
\newcommand\set[1]{\{#1\}}
\newcommand\suchthat{\: : \: }
\newcommand\tcJ{{\tilde{\mathcal J}}}
\newcommand\tcZ{{\tilde{\mathcal Z}}}
\newcommand\ttrace{{\tilde{t}}}
\newcommand\Zn{\mathbb Z^n}

\begin{document}
\title{Structurable Tori}

\author{Bruce Allison}
\address[Bruce Allison]{Department of Mathematical and Statistical Sciences\\
University of Alberta\\Edmonton, AB, Canada T6G 2G1}
\email{ballison@ualberta.ca}
\author{John Faulkner}
\address[John Faulkner]
{Department of Mathematics\\
University of Virginia\\Kerchof Hall, P.O. Box 400137\\Charlottesville VA 22904-4137 USA}
\email{jrf@virginia.edu}
\author{Yoji Yoshii}
\address[Yoji Yoshii]
{Akita National College of Technology\\
1-1 Iijimabunkyocho Akita-shi, Akita Japan 011-8511}
\email{yoji\_yoshii@yahoo.co.jp}

\dedicatory{Dedicated to the memory
of Professor Issai Kantor}
\thanks{Bruce Allison gratefully acknowledge the support
of the Natural Sciences and Engineering Research Council of Canada.}
\date{}

\begin{abstract}
The classification of structurable tori with nontrivial involution,
which was begun by Allison and Yoshii, is completed.  New
examples of structurable tori are obtained using a construction of structurable
algebras from a semilinear version of cubic forms satisfying the
adjoint identity. The classification uses techniques borrowed from
quadratic forms over $\mathbb{Z}_{2}$ and from the geometry of
generalized quadrangles.  Since structurable tori are the coordinate algebras for the
centreless cores of
extended affine Lie algebras of type BC$_1$, the results of this paper provide
a classification and new examples for this class of Lie algebras.
\end{abstract}

\maketitle
\section{Introduction}
\label{sec:intro}

The purpose of is to complete the classification of
structurable tori with nontrivial involution begun in \cite{AY}.
Structurable algebras are (in general nonassociative) algebras with
involution that are defined by an identity that is needed for their
use in the construction of $5$-graded Lie algebras \cite{K2,A2}.  If
$\Lm$ is a finitely generated abelian group, a \textit{structurable
$\Lm$-torus} is a $\Lm$-graded structurable algebra with the
property that any nontrivial homogeneous component of $\cA$ is
spanned by an invertible element (as well as the convenient
assumption that the support of $\cA$ generates the group $\Lm$).

Structurable tori with nontrivial involution
arise as coordinate algebras of  centreless cores of
extended affine Lie algebras of type BC$_1$.
The class of extended affine Lie
algebras (EALAs) is an axiomatically defined class of  Lie algebras
over a field of characteristic 0  that contains finite dimensional
split simple Lie algebras (the nullity 0 EALAs) and affine Kac-Moody
Lie algebras (the nullity 1 EALAs) as motivating
examples.
Any EALA $\cE$ has two compatible gradings, one
by a finitely generated abelian group $\Lm$ (whose rank is called
the \textit{nullity} of $\cE$) and the other by a finite root system
(whose type is called the \textit{type} of $\cE$).  Also, the structure of an EALA
$\cE$ is determined by the structure of its centreless core
$\cE_\text{cc}$, in the sense that there is a general construction
that produces all EALAs with a given centreless core \cite{N}. For a
EALA $\cE$ of type BC$_1$, it is shown in \cite{AY} that
$\cE_\text{cc}$ can be constructed from a structurable $\Lm$-torus
with nontrivial involution using a Lie algebra construction due to
Issai Kantor \cite{K2}. (The result in \cite{AY} is stated over a
field $\mathbb C$ of complex numbers, but the same proof works for any field of
characteristic 0.)  Hence, our classification of structurable tori
with nontrivial involution gives a corresponding classification of
the centreless cores of EALAs of type BC$_1$. Previously,
the
classification of the centreless cores over $\mathbb C$ had been completed for all
types except BC$_r$ \cite{BGK,BGKN,Y1,AG}.  So this paper, combined
with the work in \cite{F2} on type BC$_2$ and the work in \cite{ABG}
on type BC$_r$ ($r\ge 3$), completes the classification of the
centreless cores of EALAs over $\mathbb C$.
(We  note that the terminology for EALAs has not
fully stabilized.  We are using the term EALA in the sense of \cite{N}
where the base field is an arbitrary field of characteristic 0.
Over $\mathbb C$ these form a slightly more general class of algebras
than the class of tame EALAs studied earlier in \cite{AABGP}, but the centreless cores
obtained from algebras in the two classes coincide.
Also, EALAs in the sense of \cite{N} coincide with tame
EALAs of finite null rank studied recently in \cite{MY}.)

We now outline the contents of this paper in more detail.
Since composition algebras will appear frequently in our constructions,
we begin in \S2 by recalling some basic definitions and facts
about composition algebras.
Then, in
\S\ref{sec:tori}, we recall the definition and give some examples of
structurable algebras and structurable tori. If the involution is
trivial, a structurable torus is a Jordan torus and these were
classified by Yoshii \cite{Y1}. The remaining tori divide naturally
into three classes: I, II, and III, and we recall that trichotomy
from \cite{AY} in \S\ref{sec:classes}.

In \S\ref{sec:classI}, we classify structurable tori of
class I.  These are the tori generated by the skew elements.  We introduce
maps $\varepsilon$ measuring skewness, $\beta$ measuring commutativity, and
$\alpha$ measuring associativity.  Indeed, $\beta$ and $\alpha$ are
multiplicative versions of the commutator and the associator.  We show that
$\varepsilon$ is roughly (a multiplicative version of) a quadratic form on a
$\mathbb{Z}_{2}$-vector space with polarization $\beta$.  In the associative
case, this is precisely correct and allows a classification of tori directly from
the classification of quadratic forms over $\mathbb{Z}_{2}$.  In the
nonassociative case, although $\ep$ is not a quadratic form, we are able to
use similar methods to obtain the classification.  The main result of the section,
which was announced without proof in  \cite{AFY},
states roughly that a torus of class I is the tensor product of composition algebras
over the algebra of Laurent polynomials.

In \S\ref{sec:classII}, we briefly recall from \cite{AY} the classification
of structurable tori of class~II.  These are all constructed from a diagonal
graded hermitian form over a class I associative torus $\cB$ with involution.  We
include this result since we now better understand the possibilities for $\cB$
(in view of our work in \S\ref{sec:classI}) and since we need the
construction using hermitian forms in any case to obtain some of the tori of class III
(the class III(a) tori).

To prepare for the classification of tori of class III,
we introduce in \S\ref{sec:Cubicforms} a
semilinear version of cubic forms satisfying the adjoint identity.  We show
that they can be used to construct certain structurable algebras.  We connect
this construction with previous constructions of structurable algebras
known as the Cayley-Dickson process and the $2\times2$ matrix algebra construction.
This section is written so that it can be read independently of the rest of the paper,
since we believe that the material is of interest by itself in the theory of structurable algebras.

In \S\ref{sec:geometry}, we associative a geometry $\cI$ to a class
III torus $\cA$.  The incidence structure in $\cI$ is defined to
encode the multiplication of homogeneous components in $\cA$.
Properties of the geometry $\cI$ give a natural subdivision of class
III into III(a), III(b), and III(c). In class III(b), $\cI$ turns
out to be a generalized quadrangle of order $(2,t)$. Although in the end
we made no use of that fact, it did guide our thinking in developing
geometric properties.

In \S\ref{sec:IIIa}, we use a recognition theorem from \cite{AY} to show that
all tori of class III(a) are constructed from a diagonal hermitian form over
the quaternion algebra with nonstandard involution over the algebra of Laurent polynomials.

In \S\ref{sec:IIIb}, we classify tori of class III(b).  We see that any such torus
can be constructed using the Cayley-Dickson doubling process from a Jordan torus of degree 4 over
the algebra of Laurent polynomials.

In the final three sections we turn our attention to
class III(c) which we show contains precisely two examples.
These examples are obtained using the construction in \S\ref{sec:Cubicforms} of a structurable algebra
from a cubic form. The cubic form is defined on the space $\Herm$ of
$3\times 3$-hermitian matrices over $\cC$, where $\cC$ is either
a quaternion torus or an octonion torus.
The first of the three sections contains two
preparatory lemmas on $\Herm$, and the last two
contain the construction and classification of class III(c) tori.

For ease of future reference, the classification theorems in this
paper, Theorem \ref{thm:classI} for class I, Theorem
\ref{thm:classII} for class II, Theorem \ref{thm:classIIIa} for
class III(a), Theorem \ref{thm:IIIb} for class III(b), and Theorem
\ref{thm:IIIc} for class III(c), are presented to stand alone as
much as possible.  We note also that the tori of class III(a) that
are not associative, as well as the tori of class III(b) and
III(c), are new in this work. As discussed above, each of these new
structurable tori in turn determines a new family of EALAs of type
BC$_1$.

Throughout the paper, \textit{we assume that $F$ is a field of characteristic $\neq2$ or
$3$}.  (Including the positive characteristic case requires no additional work.)

\section{Composition algebras}
\label{sec:composition}

At  several points in this work, we will make use of composition algebras to construct examples.
In this section, we recall some of the basic definitions and facts that we will  need about composition
algebras.  We first fix some terminology and notation  for algebras.

By an \emph{algebra}, we mean a (not necessarily associative) unital algebra $\cA$ over~$F$.
If $\cA$ is an algebra and $x,y,z\in\cA$, we use the notation $[x,y] = xy-yx$ and
$(x,y,z) = (xy)z - x(yz)$ for the commutator and associator.
If $x\in \cA$, we define endomorphisms  $L_x$ and $R_x$  of $\cA$ by
$L_x y = xy$ and $R_x y = yx$.

An \emph{algebra
with involution} is a pair $(\cA,*)$ consisting of an
algebra $\cA$ with an anti-automorphism $*$ of period 2.
When no confusion will result,
we usually denote an algebra with involution $(\cA,*)$ simply by $\cA$.
For $\sg=\pm$ we let
\[\cA_{\sg}=\{x\in\cA:x^{*}=\sg x\},\]
and we call elements of $\cA_-$ (resp.~$\cA_+$) \emph{skew} (resp.~\emph{hermitian}).
We define the \emph{centre}  of the algebra with involution $\cA$ to be
\begin{equation}
\label{eq:centredef}
\cZ(\cA) :=
\{z\in\cA_{+}:[z,\cA]
=(z,\cA,\cA)=(\cA,z,\cA)=(\cA,\cA,z)=0 \}.
\end{equation}

If $\cA$ is a commutative algebra, then the identity map is an involution on $\cA$
called the \textit{trivial involution}.

For the rest of this section, assume that \emph{$K$ is a commutative associative
algebra over $F$.}

If $\cU$ is a $K$-module, then a map $q: \cU \to K$ is called
a \textit{quadratic form} over $K$ if
$q(ra) = r^2 q(a)$ for $r\in K$, $a\in \cU$, and
if $f(a,b) := q(a+b) - q(a) -q(b)$ defines a $K$-bilinear map
$f:\cU\times\cU \to K$ called the \emph{linearization} of $q$.
In that case, $q$ is said to be \textit{nondegenerate} if $f$
is nondegenerate ($f(a,\cU) = 0 \implies a = 0$ for $a\in \cU$).
It is usual to abuse notation and denote the linearization of a quadratic form
$q$ by $q$.

Recall \cite[p.156]{Mc} that a \textit{composition
algebra} over $K$ is an  algebra $\cC$ over $K$
with a nondegenerate quadratic form $n:\cC%
\rightarrow K$, called the \textit{norm}, which permits \textit{composition}:%
\[
n(1)=1 \andd n(ab)=n(a)n(b),\text{ for }a,b\in\cC.
\]
We then define the \textit{trace} $t:\cC\rightarrow K$ by
$t(a)=n(a,1)$.
It is known \cite[p.156]{Mc} that the map $\barinv : \cC \to \cC$
defined by
\begin{equation*}
\bar{a}=t(a)1-a
\end{equation*}
is an involution over $K$ called the \emph{canonical involution} on $\cC$.
Unless mentioned to the contrary, \emph{we will always regard a composition algebra
as an algebra with involution using the canonical involution}.
It also known \ibid\ that $\cC$ is an
alternative algebra and that
$n(a)1=a\bar{a}=\bar{a}a$. Therefore,
$a\in\cC$ is invertible if and only if $n(a)$ is invertible in $K$, in
which case $a^{-1}=n(a)^{-1}\bar{a}$. Also, we have $n(\bar a) = n(a)$, $t(\bar a) = t(a)$ and
\[n(a,b)=t(a\bar{b}).\]
Further,
\[t((ab)c)=t(a(bc)),\]
so we can write $t(abc)=t((ab)c)$. Also,
alternativity shows that $(ab)c=1$ if and only if $a(bc)=1$, which we can thus
write as $abc=1$. Finally, if $abc=1$, then $bca=cab=1$.

To construct examples of composition algebras  we recall the classical \textit{Cayley-Dickson process}.
(See for example \cite[\S 1]{BGKN} or \cite[pp.~160-163]{Mc}).
If $\cD$ is an algebra with involution $\barinv$, and $\mu$ is a unit in the centre
of $\cD$, we let
\[\CD(\cD,\mu) = \cD \oplus u \cD\]
be the algebra with involution whose product
and involution are defined by
\begin{equation}
\label{eq:CDdef}
(a+ub)(c+ud) = (ac + \mu d \bar b) + u(\bar a d + cb)  \andd \overline{a+ub} = \bar{b} - ub.
\end{equation}

If  $\mu_1,\mu_2,\mu_3$
are units in $K$, we can iterate the CD-process starting at $K$ with the trivial involution $\barinv$.
We successively construct algebras with  involution:
\begin{equation}
\label{eq:CDprocess}
\begin{aligned}
&\hspace{1.5em}K,\\
\CD(K,\mu_1) &= K\oplus x_1 K,\\
\CD(K,\mu_1,\mu_2) &= \CD(\CD(K,\mu_1),\mu_2) = \CD(K,\mu_1) \oplus x _2 \CD(K,\mu_1),\\
\CD(K,\mu_1,\mu_2,\mu_3) &= \CD(\CD(K,\mu_1,\mu_2),\mu_3)\\
&= \CD(K,\mu_1,\mu_2) \oplus x_3 \CD(K,\mu_1,\mu_2).
\end{aligned}
\end{equation}
It is well known (see for example
the argument in \cite[p.~58]{S}) that these algebras
are composition algebras over
$K$ with norm defined by $n(a) = a\bar a$ and that their involution
(defined in \eqref{eq:CDdef}) is the canonical involution.
The  first two of these algebras are commutative and associative,
the second is associative and the last is alternative.
They have $K$-bases $\set{1}$, $\set{1,x_1}$, $\set{1,x_1,x_2,x_2x_1}$ and
$\set{1,x_1,x_2,x_2x_1,x_3,\allowbreak x_3x_1,x_3x_2,x_3(x_2x_1)}$ respectively;
and they have generating sets $\set{1}$, $\set{x_1}$, $\set{x_1,x_2}$ and
$\set{x_1,x_2,x_3}$ respectively as algebras over $K$.  The elements in these generating sets
are called \textit{canonical generators} and they satisfy $x_i^2 = \mu_i$,
$x_ix_j = -x_jx_i$ for $i\ne j$ and $(x_ix_j)x_k = -x_i(x_jx_i)$ for  $i,j,k\ne$.
The algebras with involution $\CD(K,\mu_1,\mu_2)$ and
$\CD(K,\mu_1,\mu_2,\mu_3)$ are called \textit{quaternion} and \textit{octonion} algebras
over $K$ respectively.

The theorem of Hurwitz \cite[Theorem 3.25]{S} tells us that if $K$ is a field,
every composition algebra over $K$ is obtained as in the previous paragraph.
In particular, the
\emph{2-dimensional composition algebras} over $F$ are the algebras $E$
over $F$ that have a basis $1,e$ satisfying  $e^2 \in F^\times$.
In that case the canonical involution of $E$ is
the automorphism $\sg_E$ satisfying $\sg_E(e) = -e$.  There are only two
possibilities for a 2-dimensional composition algebra $E$: either
$E/F$ is a quadratic field extension and $\sg_E$ generates the
Galois group of $E/F$, or $E$ is isomorphic to $E = F
\oplus F$ with the exchange involution. In the second case we say
that $E$ is \textit{split}.

\section{Structurable algebras and tori}
\label{sec:tori}

In  this section, we recall the definition and give some examples of structurable tori.
We first recall some terminology and notation for graded algebras.

Assume that \emph{$\Lm$ is an arbitrary abelian group}.
A vector
space $\cV$ over $F$ is \textit{graded }by $\Lm$
if $\cV =\bigoplus_{\lm\in\Lm}\cV^{\lm}$, where $\cV^\lm$ is a subspace of $\cV$ for $\lm\in\Lm$.
If $M$ is a subgroup of $\Lm$, let
\[
\cV^{M}:=\bigoplus_{\lm\in M}\cV^{\lm},
\]
an $M$-graded space. The \textit{support of } $\cV$ is
\[
\supp(\cV):=\{\lm\in\Lm:\cV^{\lm}\neq0\}.
\]
$\cV$ is said to be \emph{finely graded} if
$\dim(\cV^\lm) = 1$ for all $\lm\in \supp(\cV)$.

Two gradings of $\cV$ by abelian groups $\Lm$ and $\Lm'$ are
\textit{isomorphic gradings} if there is a group isomorphism $\theta
:\Lm\rightarrow\Lm'$ with $\cV^{\lm}%
=\cV^{\theta(\lm)}$. More generally, an \textit{isograded
isomorphism} is a linear isomorphism $\ph:\cV\rightarrow\cW$
with $\ph(\cV^{\lm})=\cW^{\theta(\lm)}$, where
$\cV$ is $\Lm$-graded, $\cW$ is $\Lm'$-graded,
and $\theta:\Lm\rightarrow\Lm'$ is a group isomorphism. In
this case, we write $\cV\ig \cW$. If $\Lm
=\Lm'$ and $\theta=\id$, this reduces to the usual notion of a
\textit{graded isomorphism}, and we write $\cV\simeq_{\Lm
}\cW$.

A \textit{graded algebra} is an algebra $\cA$ which is a graded space
with $\cA^{\lm}\cA^{\mu}\subset\cA^{\lm+\mu}$.
For a \textit{graded algebra with involution}, we also require the involution
to be a graded vector space isomorphism; i.e. $(\cA^{\lm})^{*}=\cA^{\lm}$.
In that case, the subspaces $\cA_+$ and $\cA_-$
are graded subspaces of $\cA$.  Moreover, if $\cA$ is a finely graded algebra with involution
then $x^* = \pm x$ for each homogeneous element $x$ of $ \cA$.

We extend the notions of isograded isomorphism and graded isomorphism  to
isomorphisms of graded algebras and isomorphisms of
graded algebras with involution in the obvious fashion.  If $\cA$ and
$\cA'$ are graded algebras with involution, we write $\cA \ig \cA'$
to mean that $\cA$ and $\cA'$ are isograded isomorphic as algebras with involution.
(If we mean that $\cA$ and $\cA'$ are isograded isomorphic as algebras without involution,
we will say this specifically.)

Next we recall
\cite[(3) and Corollary 5]{A1}\footnote{In \cite{A1} there is an
overriding assumption that algebras are finite dimensional.
However that assumption is not used in the sections (1--5 and 8(iii)) of \cite{A1}
that we use in the infinite dimensional setting.}
that a \textit{structurable algebra}
is a (unital) algebra
$\cA$ with an involution $*$ satisfying
\[
\{xy\{zwq\}\}-\{zw\{xyq\}\}=\{\{xyz\}wq\}-\{z\{yxw\}q\}
\]
where%
\[
\{xyz\}:=(xy^{*})z+(zy^{*})x-(zx^{*})y.
\]
Examples  of structurable algebras include associative algebras with involution,
alternative algebras with involution and Jordan algebras
with the trivial involution \cite[Theorem 13]{A1}.
In fact, unless stated to the contrary, we will regard Jordan algebras as
algebras with the trivial involution, and hence as structurable algebras.

As the term suggests, a \emph{$\Lm$-graded structurable algebra} is a $\Lm$-graded algebra
with involution $\cA$ such that $\cA$ is structurable.

Suppose that $\cA$ is a finely $\Lm$-graded structurable algebra.  A homogeneous
element $x$ of $\cA$ is said to be \textit{invertible} if there exists $x^{-1}\in \cA$ so that
\[xx^{-1}=x^{-1}x=1 \andd [L_{x},L_{x^{-1}}]=[R_{x},R_{x^{-1}}]=0.\]
In that case, $x^{-1}$ is unique and if $x\in \cA_\sg^\lm$, where
$\lm\in \Lm$, $\sg = \pm$, we have $x^{-1}\in \cA_\sg^{-\lm}$.
\cite[Prop.~3.1]{AY} .  We call $x^{-1}$ the \emph{inverse} of $x$.

\begin{remark}
\label{rem:invertible}
Suppose that $\cA$ is a finely $\Lm$-graded structurable algebra,
and let  $x$ be a homogeneous element of $\cA$.

(a) To check that $x$ is invertible with inverse
$y$, it is sufficient to check that $xy=1$ and that $[L_x,L_y]= 0$.  Indeed, in that case, we have
$yx = L_yL_x 1 = L_xL_y1 = xy =~1$.   Also, conjugating the equality $[L_x,L_y]= 0$ by
the involution, we get $[R_x,R_y]= 0$.

(b) $x$ is invertible if and only if there exists an element $\hat x$ in $\cA$ such that
$\{x,\hat x,z\} = z$ for all $z\in \cA$.  (Elements $x$ with the latter property
are said to be \emph{conjugate invertible}.) In that case, $\hat x$ is unique and $\hat x = \ep x^{-1}$,
where $x^* = \ep x$ with $\ep = \pm 1$ \cite[Prop.~3.1]{AY}.

(c) If $x = s$ is skew, then $s$ is invertible if and only if $L_s$
(or equivalently $R_s$) is invertible, in which case
$L_s^{-1} = L_{s^{-1}}$ and $R_s^{-1} = R_{s^{-1}}$
\cite[Lemma~2.9]{AY}.  Moreover, in that case $sw$ is invertible
for any invertible homogeneous element $w$ of
$\cA$ (by (b) and \cite[Prop.~8.2 and Theorem 11.4]{AH}).

(d) If $\cA$ is associative, alternative or Jordan, then the notion of invertibility used
here coincides with the usual notion.
\end{remark}

\begin{definition}
A \textit{structurable $\Lm$-torus} is a $\Lm$-graded structurable algebra $\cA$ satisfying
\begin{mylist}
\item[(ST1)] $\cA$ is a finely graded.
\item[(ST2)] Each nonzero homogeneous element of $\cA$ is invertible.
\item[(ST3)] $\Lm$ is generated as a group by $\supp(\cA)$.
\end{mylist}
\end{definition}

\begin{remark}
\label{rem:subtori}
If $\cA$ is a structurable $\Lm$-torus and $M$ is a subgroup of
$\Lm$, then $\cA^{M}$ is clearly a subalgebra of $\cA$
satisfying (ST1) and (ST2). Thus, if $M$ is generated by $\supp(\cA^{M})$,
then $\cA^{M}$ is a structurable $M$-torus.  In particular,
if $\supp(\cA) = \Lm$, then $\cA^M$ is a structurable torus for any subgroup $M$ of $\Lm$.
\end{remark}

\begin{remark}
\label{rem:powers}
If $\cA$ is a structurable $\Lm$-torus and $x$ is a nonzero homogeneous element
of $\cA$, then the (unital) subalgebra of $\cA$ that is generated by $x$ and $x^{-1}$
is commutative and associative \cite[Corollary~7.8]{AY}.
Hence, the powers $x^k$ for $k\in \bbZ$ are all nonzero.
\end{remark}

Examples of structurable $\Lm$-tori include $\Lm$-graded associative algebras with involution
satisfying (ST1)--(ST3), $\Lm$-graded alternative algebras with involution satisfying (ST1)--(ST3),
and $\Lm$-graded Jordan algebras (with the trivial involution) satisfying (ST1)--(ST3)
\cite[Examples 4.2 and 4.3]{AY}. These are called respectively
\emph{associative $\Lm$-tori with involution, alternative $\Lm$-tori with involution
and Jordan $\Lm$-tori}.

More specifically we have the following basic examples:

\begin{example}
\label{ex:laurent}
Let
\[\cP(n) := F[t_1^{\pm1},\dots,t_n^{\pm1}]\]
be the algebra of Laurent polynomials over $F$, where $n\ge 0$. We
call the elements $t_1,\dots,t_n$  \textit{canonical
(Laurent) generators} for $\cP(n)$. The algebra $\cP(n)$ has a
unique $\Zn$-grading so that $\deg(t_i) = \ep_i$ for each $i$, where
$\set{\ep_1,\dots,\ep_n}$ is the standard basis for $\Zn$. Then
$\cP(n)$ is a commutative associative $\Zn$-torus with trivial
involution, and hence a structurable $\Zn$-torus with trivial
involution. (If $n=0$, then $\bbZ^0 = 0$ and $\cP(0) = F$.)
\end{example}

\begin{example}
\label{ex:C(n)} (\cite[\S1]{BGKN})
Suppose that $n=1,2,3$, and
let
\[
\Cay(n):=\CD(\cP(n),t_{1},\ldots t_{n}).
\]
We give the algebra $\Cay(n)$ the \textit{canonical} $\Zn$-grading by
assigning the degrees in the standard basis $\ep_1,\dots,\ep_n$ to the
canonical generators $x_1,\dots,x_n$ of $\Cay(n)$ (in order).
We use the notation $\Cay(n)$ for
$(\Cay(n),\natural)$, where $\natural$ is the canonical involution  (satisfying $x_i^\natural = -x_i)$.
We also set
\[\Cay(0)=F\]
(where, as in Example \ref{ex:laurent},
we regard $F$ as an algebra with trivial involution
graded by the trivial group $\bbZ^0 = 0$).  Then,
$\Cay(n)$ is an alternative $\Zn$-torus with involution and
hence it is a structurable $\Zn$-torus for $n=0,1,2,3$.
\end{example}

\begin{example}
\label{ex:A*(2)} The graded algebra $\Cay(2)$ described in
Example \ref{ex:C(n)} also has the involution $*_{m}$,
called the \textit{main involution}, that fixes the canonical
generators. We use the notation $\Cay_*(2)$ for $(\Cay(2),*_{m})$.
Then $\Cay_*(2)$ is a structurable $\bbZ^2$-torus.
\end{example}

\begin{remark}
\label{rem:altgrading}  Unless indicated to the contrary, we will
regard $\cP(n)$, $\Cay(n)$ and (if $n=2$) $\Cay_*(2)$
as $\bbZ^n$-graded algebras with involution as above.
However it is occasionally convenient to grade these algebras with involution
by a free abelian group $\Lm$ with basis $\lm_1,\dots,\lm_n$
by assigning the degree $\lm_i$ to the $i^\text{th}$ canonical generator.
In that case, we say that $\cP(n)$, $\Cay(n)$ or $\Cay_*(2)$
has the $\Lm$-grading \textit{determined by the basis}
$\lm_1,\dots,\lm_n$.
\end{remark}

\begin{remark}
\label{rem:assoctoriwithout} We sometimes refer to \emph{associative
$\Lm$-tori} (without involution). By definition these are
$\Lm$-graded associative algebras (without involution) satisfying
(ST1)--(ST3).
\end{remark}

\section{Classes I, II and III}
\label{sec:classes}

For the remainder of this paper, \textit{we assume that $\Lm$ is a free abelian
group of finite rank.}

In \cite{AY}, structurable tori with nontrivial involution
were placed into 3 classes I, II and III.  To recall that trichotomy,
we introduce some notation.

Suppose that $\cA$ is a structurable $\Lm$-torus.  Let
\[
S = S(\cA) :=\supp(\cA)
\]
and
\[
S_{\sg}=S_{\sg}(\cA):= \supp(\cA_\sg) = \{\lm\in S:\cA^{\lm
}\subset\cA_{\sg}\}.
\]
We also set
\[\Lm
_{-}=\Lm_{-}(\cA):=\left\langle S_{-}\right\rangle,\]
the subgroup generated by $S_{-}$.
Next let
\[
\cZ = \cZ(\cA)\]
be the centre of the algebra with involution $\cA$ (as defined
in \eqref{eq:centredef}).
Then $\cZ$ is a graded subalgebra of $\cA$ and we let
\[\Gm=\Gm(\cA):= \supp(\cZ) = \{\lm\in S:\cA^{\lm
}\subset\cZ\}.\]
By \cite[Prop.~6.7]{AY}, $\Gm$ is a
subgroup of $\Lm$ called the \textit{central grading group}.

If $S_{-}=\emptyset$ (that is the involution is trivial), then $\cA$ is
a Jordan torus \cite[Example 4.3]{AY}.  The classification of Jordan tori was done by Yoshii
in \cite{Y1}, and consequently we are interested in this work in the
case when $S_{-}\ne \emptyset$.

If $S_{-}\neq\emptyset$, let
\[
\cE:=\cA^{\Lm_{-}} \andd \cW:=\cA%
^{\Lm\smallsetminus\Lm_{-}},
\]
so
\[\cA=\cE\oplus\cW\]
with $\cE\cE \subset\cE$ and
$\cE\cW+\cW\cE\subset\cW$. Then, $\cE$ is the subalgebra of $\cA$
that is generated by $\cA_-$ \cite[Prop.~8.1]{AY}.

\begin{definition}
If $\cA$ is a structurable $\Lm$-torus, we say that $\cA$ has
\textit{class} I, II, or III, if $S_{-}\neq\emptyset$ and the
corresponding condition below holds:
\begin{mylist}
\item[I.] $\cE=\cA$,
\item[II.] $\cE\neq \cA$ and $\cW\cW\subset\cE$,
\item[III.] $\cE\neq \cA$ and $\cW\cW\not \subset \cE$.
\end{mylist}
\end{definition}

Allison and Yoshii \cite[Theorem 9.5]{AY} classified structurable tori of
class II and obtained basic results which we shall use to obtain the
classification in the remaining two classes.

\section{Structurable tori of class I}
\label{sec:classI}

In this section we obtain
the construction and classification of structurable tori of class I.
We will construct these tori as tensor products of the basic examples described
in \S\ref{sec:tori}.  Thus we begin with a few simple observations about tensor
products of graded algebras.

\begin{definition}
\label{tpgrading}
If $\cA'$ and $\cA''$ are algebras with involutions (denoted by $*'$ and $*''$ respectively)
then the tensor product algebra $\cA'\ot\cA''$ is an algebra with involution $*'\ot *''$.
(Here and subsequently an unadorned symbol $\ot$ means $\ot_F$.)
Also, if $\cA'$ is a $\Lm'$-graded algebra  with involution
and $\cA''$ is a $\Lm''$-graded algebra with involution, then the
$\cA'\ot\cA''$ is a $(\Lm'\oplus\Lm'')$-graded
algebra with involution, where
\[
(\cA'\ot\cA'')^{(\lm^{\prime
},\lm'')}=\cA^{\prime\lm'}\ot
\cA^{\prime\prime\lm''}.
\]
We call this grading on $\cA'\ot\cA''$ the \textit{tensor product}
grading. These definitions extend in the obvious fashion to tensor
products of several graded algebras with involution. Of course all
of these definitions are made in the same way for graded algebras
without involution.
\end{definition}

\begin{lemma}
\label{lem:tensor} Suppose
that $\cA$ is a finely $\Lm$-graded algebra
with involution such that $\supp(\cA) = \Lm$ and the product of any two nonzero homogeneous
elements of $\cA$ is not zero.
Suppose $\Lm'$ and $\Lm''$
are subgroups of $\Lm$ such that $\Lm=\Lm'\oplus
\Lm''$,
\[
\lbrack\cA^{\Lm'},\cA^{\Lm''%
}]=0\quad\text{and}\quad(\cA^{\Lm'},\cA%
^{\Lm''},\cA)=(\cA,\cA^{\Lm
'},\cA^{\Lm''})=0.
\]
Then $\cA\simeq_{\Lm}\cA^{\Lm'}\ot
\cA^{\Lm''}$.
\end{lemma}

\begin{proof}
Clearly, $\eta:\cA^{\Lm'}\ot\cA^{\Lm
''}\rightarrow\cA$ with $\eta:x\ot y\rightarrow xy$
is a homomorphism of $\Lm$-graded vector spaces. Now, if $x\in\cA%
^{\Lm'}$ and $y\in\cA^{\Lm''}$, we have
\[
(\eta(x\ot y))^{*}=(xy)^{*}=y^{*}x^{*}=x^{*}y^{*}%
=\eta(x^{*}\ot y^{*})=\eta((x\ot y)^{*}).\]
Also, if $x_{i}\in\cA^{\Lm'}$ and $y_{i}\in
\cA^{\Lm''}$, we have
\begin{align*}
(x_{1}y_{1})(x_{2}y_{2}) &=x_{1}(y_{1}(x_{2}y_{2}))
=x_{1}((y_{1}x_{2})y_{2})
=x_{1}((x_{2}y_{1})y_{2})
=x_{1}(x_{2}(y_{1}y_{2}))\\
&=(x_{1}x_{2})(y_{1}y_{2})
\end{align*}
It remains then to show that $\eta$ is a bijection.  Since
$\eta$ is a homomorphism of finely graded vector spaces, it is enough
to show that the restriction of $\eta$ to $\cA^{\lm'}\ot \cA^{\lm''}$ is nonzero
for $\lm'\in\Lm'$ and $\lm''\in\Lm''$.  This follows from our assumptions that
$\supp(\cA) = \Lm$ and the product of any two nonzero homogeneous
elements of $\cA$ is not zero.
\end{proof}

We also need the following lemma about base ring extension for tensor products.

\begin{lemma}
\label{lem:baseringtensor} Suppose we are given an
$F$-algebra homomorphism $K_1\ot K_2\to K$, where
$K_1$, $K_2$ and $K$ are associative and commutative $F$-algebras.
If $\cA_i$ is a $K_i$-algebra for $i=1,2$,
then
\begin{equation}
\label{eq:baseringtensor}
(K\ot_{K_1}\cA_1)\ot_{K}(K \ot_{K_2}\cA_2)
\simeq
K \ot_{K_1\ot K_2} (\cA_1\ot \cA_2).
\end{equation}
as $K$-algebras.  Moreover, if $\cA_i$ is a $K_i$-algebra with involution for $i=1,2$,
then \eqref{eq:baseringtensor} is an isomorphism of $K$-algebras with involution.
\end{lemma}

\begin{proof}  Note that, for $i=1,2$, we have the $F$-algebra homomorphism $K_i\to K_1\ot K_2 \to K$,
so the base ring extensions $K\ot_{K_i}\cA_i$ on the left hand side of
\eqref{eq:baseringtensor} make sense.
Also, $A_1\ot A_2$ is naturally a $K_1\ot K_2$-algebra, so the base ring extension
on the right hand side also makes sense.

It follows from the universal properties of the tensor product that there exists a unique
additive map $\ph$ from the left hand side to the right hand side under which
\[(\al_1\ot_{K_1} a_1)\ot_{K}(\al_2\ot a_2)
 \mapsto
\al_1\al_2\ot_{K_1\ot K_2}(a_1\ot a_2) \]
for $\al_1,\al_2\in K$, $a_1\in \cA_1$, $a_2\in \cA_2$.
Similarly there exists a unique additive map
$\psi$ from the right hand side to the left hand side under which
\[\al\ot_{K_1\ot K_2}(a_1\ot a_2)
\mapsto
(\al \ot_{K_1} a_1)\ot_{K}(1\ot a_2)
 \]
for
$\al\in K$, $a_1\in \cA_1$, $a_2\in \cA_2$.
One immediately sees that $\ph$ and $\psi$ are  homomorphisms of $K$-algebras  with
$\ph\circ \psi = \id$ and $\psi\circ \ph = \id$.  Furthermore, in the involutory case,
these maps preserve the involutions.\end{proof}

\begin{remark} Lemma \ref{lem:baseringtensor} has an obvious extension (with the same proof)
allowing $n$ tensor factors $\cA_1,\dots,\cA_n$, where $n\ge 1$.
\end{remark}

\begin{remark}
\label{rem:tensor} (a)\ms If $\cA_i$
is a structurable $\Lm_i$-torus for $i=1,\dots,k$ and if
$\cA_1\ot\cA_2\ot \cdots\ot\cA_k$ is a structurable
algebra, then it is clear that $\cA_1\ot\cA_2\ot \cdots\ot\cA_k$ satisfies (ST1)-(ST3)
and hence $\cA_1\ot\cA_2\ot \cdots\ot\cA_k$ is a structurable $\Lm_1\oplus\dots\oplus\Lm_k$-torus.

(b)\ms If $\cA$ is a structurable $\Lm$-torus, then it follows from (a) that
$\cA\ot \cP(r)$
is a structurable $\Lm\oplus \bbZ^r$-torus for $r\ge 0$.
\end{remark}

Finite dimensional simple structurable algebras can be constructed by taking the tensor
product of two composition algebras over the base field
\cite{K1}, \cite[\S 8]{A1}.  (In \cite{K1}, Kantor worked with
\textit{conservative algebras of second degree with left unit} which were
later seen to be structurable algebras with a modified product \cite{AH}.) We now present an
infinite dimensional adaptation of this for class I structurable tori.

\begin{proposition}
\label{prop:constructionI}
If $0\le k\le 3$ and $r \ge 0$, then
\[\Cay(3)\ot \Cay(k) \ot \cP(r)\]
is a structurable $\bbZ^{3+k+r}$-torus of class I .
\end{proposition}

\begin{proof} Let $\cA = \Cay(3)\ot \Cay(k) \ot \cP(r)$.
Then $\cA_-$ is spanned by elements of the form
$s_1\ot h_2\ot f$ or $h_1\ot s_2\ot f$, where $s_i$ is skew and $h_i$ is hermitian.
It is clear from this that $\cA_-$ generates $\cA$ as an algebra. Thus,
it suffices to show that $\cA$ is a structurable torus.  By Remarks
\ref{rem:subtori} and \ref{rem:tensor}(b), we can assume that
$\cA = \Cay(3)\ot \Cay(3)$.  Furthermore, by Remark \ref{rem:tensor}(a), it suffices to show
that $\cA$ is a structurable algebra.  For this, we let $K_1 = K_2 = \cP(3)$
and we let $K$ be the quotient field of the integral domain $K_1\ot K_2$.
Then, by Lemma \ref{lem:baseringtensor}, we have
\[
(K\ot_{K_1}\Cay(3))\ot_{K}(K \ot_{K_2}\Cay(3))
\simeq
K \ot_{K_1\ot K_2} \cA
\]
as $K$-algebras. Now the left hand side is the tensor product of two octonion algebras over the field $K$
and hence it is structurable \cite[\S 8(iv)]{A1}. Thus, since the $F$-algebra $\cA$
embeds in the right hand side, it is structurable.
\end{proof}

We will show in the main theorem of this section (Theorem
\ref{thm:classI}) that any class~I structurable
torus that is not associative is (up to isograded isomorphism)
one of the tori described in the preceding proposition.
In the same theorem we will show that any class~I structurable
torus that is associative
arises as a tensor product of the tori
$\Cay(1)$, $\Cay(2)$, $\Cay_*(2)$ and $\cP(r)$.

With this goal in mind, we next discuss
some properties of skew-elements in structurable
algebras.
If $\cA$ is a structurable algebra then $\cA$ is \textit{skew-alternative}, which means that
\[(s,x,y) = -(x,s,y) = (x,y,s)\]
for $x,y\in \cA$ and $s\in \cA_-$ \cite[Prop.~1]{A1}.  It follows that if $s\in \cA_-$ and
$x\in \cA$, we have
\[s(sx) = s^2 x,\quad (sx)s = s(xs) \andd (xs)s = xs^2.\]

Using these facts we can prove the following two lemmas:

\begin{lemma}
\label{lem:strid} The following hold in a structurable algebra $\cA$
with $r,s,t\in\cA_{-}$, $x,y\in~\cA$:

(a)\ $r(s,r,x)+(r,sr,x)=0,$

(b)\ $(r^{2},x,y)=(r,rx,y)+r(r,x,y),$

(c)\ $(r^{2},s,x)=2r(r,s,x)+(r,[r,s],x).$
\end{lemma}

\begin{proof}
By expanding associators, (a) becomes $-r(s(rx))+(r(sr))x=0$, which is the
left Moufang identity \cite[(43)]{A1}. Also, after expanding associators, (b)
follows from $r(rx)=r^{2}x$. Letting $x=s$ and $y=x$ in (b), we have by (a)
\begin{align*}
(r^{2},s,x) &  =(r,rs,x)+r(r,s,x)\\
&  =(r,rs,x)+r(r,s,x)-r(s,r,x)-(r,sr,x)\\
&  =2r(r,s,x)+(r,[r,s],x)
\end{align*}
by skew alternativity.
\end{proof}

In the following lemma, we use the notation $x\circ y=xy+yx$.

\begin{lemma}
\label{lem:skewgen}
Let $\cA$ be a structurable algebra that is
generated as an algebra by $\cA_{-}$.  Then:

(a)\ms $\cA=\cA_{-}\oplus\cA_{-}\circ\cA_{-}.$

(b)\ms If $(\cA_{-},\cA_{-},\cA_{-})=0$ then
$\cA$ is associative.

(c)\ms If $x\in\cA$ and $[x,\cA]=0$, then $(x,\cA%
,\cA)=(\cA,x,\cA)=(\cA,\cA,x)=0$.
\end{lemma}

\begin{proof}
(a) has been shown in \cite[Lemma 14]{A1}. In particular, $\cA_{+}$
is spanned by $r^{2}$ with $r\in\cA_{-}$. \ For (b), assume that
$(\cA_{-},\cA_{-},\cA_{-})=0$. \ Since $[r,s]\in
\cA_{-}$, letting $x=t$ in Lemma \ref{lem:strid}(c) gives
$(r^{2},s,t)=0$. Hence, $(\cA,\cA_{-},\cA_{-})=0$,
so also $(\cA_{-},\cA,\cA_{-})=(\cA_{-},\cA_{-},\cA)=0$,
by skew alternativity. Now Lemma
\ref{lem:strid}(c) gives $(r^{2},s,x)=0$. Thus $(\cA,\cA%
_{-},\cA)=0$, so also $(\cA_{-},\cA,\cA%
)=(\cA,\cA,\cA_{-})=0$. \ Lemma \ref{lem:strid}(b) now
gives $(r^{2},x,y)=0$, so $(\cA,\cA,\cA)=0$.
\ Finally, (c) is proved in \cite[Lemma 24]{A1}.
\end{proof}

For the rest of this section, \emph{we assume that
$\cA$ is a class I structurable $\Lm$-torus}.
\emph{We use the notation $S$, $S_\sg$, $\Lm_-$, $\cZ$, $\Gm$,
\dots\ from \S\ref{sec:classes}}.

We will  make frequent use of the following facts
from \cite[Prop.~8.1 and 8.11]{AY} for class I tori:
\begin{itemize}
\item $\cA$ is generated as an algebra by $\cA_{-}$.
\item $S=\Lm$.
\item If $0\neq s\in\cA_{-}^{\lm}$, then $s^{2}\in\cZ$.
\end{itemize}

\begin{lemma}
\label{lem:associate} Let $\cA$ be a structurable $\Lm$-torus of
class I and let $r$, $h$, $x$, $y$ be homogeneous in $\cA$ with
$r\in\cA_{-}$ and $h\in\cA_{+}$.  Then:

(a)\ms If $[r,x]=0$, then $r$, $x$ and $y$ associate in all orders.

(b)\ms $h = st$ for some commuting homogenous elements $s,t$ of $\cA%
_{-} $.

(c)\ms If $x\neq0\neq y$, then $xy\neq0$.
\end{lemma}

\begin{proof}
For (a), we can assume that $r\neq0$. Suppose first that $x=s\in
\cA_{-}$ and $[r,s]=0$. Since $r^2\in \cZ$,
Lemma \ref{lem:strid}(c) shows
$(r,s,y)=0$. \ Suppose next that $x=h\in\cA_{+}$ and $[r,h]=0$. Then
$r(r^{-1}h-hr^{-1})r=r(r^{-1}h)r-r(hr^{-1})r=hr-rh=0$, so $r^{-1}$ commutes
with $h$.  Set $s=hr^{-1}=r^{-1}h\in\cA_{-}$ in which case $[r,s]=0$.
Now the first case and Lemma \ref{lem:strid}(a) give
$(r,h,y)=(r,sr,y)=-r(s,r,y)=0$. \ Thus, $(r,x,y)=0$ in all cases. \ Applying
the involution gives $(y,x,r)=0$ and (a) holds by skew alternativity.

For (b), we know from Lemma \ref{lem:skewgen}(a) that we can write $h=st$ with
$s,t$ homogeneous in $\cA_{-}$. \ Also, $h=h^{*}=ts$, so $[s,t]=0$.

Clearly (c) holds if $x\in\cA_{-}$. \ If $x=h\in\cA_{+}$, we
can write $h=st$ as in (b). \ Now $hy=(st)y=s(ty)\neq0$ by (a).
\end{proof}

We next
define functions
\[
\ep  :\Lm\rightarrow\{\pm1\},\ \
\beta  :\Lm\times\Lm\rightarrow F,\ \
\al  :\Lm\times\Lm\times\Lm\rightarrow F,\ \
\mu  :\Lm\times\Lm\times\Lm\rightarrow\{\pm1\}.
\]
We will see in the next lemma that $\beta$ and $\al$ also take values in
$\{\pm1\}$.

Since the support of $\cA$ is $\Lm$, if $\lm\in\Lm$ we can
choose a unique $\ep(\lm)=\pm1$ so that
\[
x^{*}=\ep(\lm)x
\]
for $x\in\cA^{\lm}$. \ Also, if $\lm_{i}\in\Lm$, we can
define $\beta(\lm_{1},\lm_{2})\in F^{\times}$ and $\al(\lm
_{1},\lm_{2},\lm_{3})\in F^{\times}$ by
\begin{align*}
x_{2}x_{1} &  =\beta(\lm_{1},\lm_{2})x_{1}x_{2},\\
x_{1}(x_{2}x_{3}) &  =\al(\lm_{1},\lm_{2},\lm_{3})(x_{1}%
x_{2})x_{3},
\end{align*}
for $x_{i}\in\cA^{\lm_{i}}$. \ The scalars $\beta(\lm
_{1},\lm_{2})$ and $\al(\lm_{1},\lm_{2},\lm_{3})$ are
uniquely determined by Lemma \ref{lem:associate}(c). \ We can view $\beta$ as
a \textit{multiplicative commutator} and $\al$ as a \textit{multiplicative
associator}.
We also define
\[
\mu(\lm_{1},\lm_{2},\lm_{3})=\ep(\lm_{1}%
)\ep(\lm_{2})\ep(\lm_{3})\ep(\lm
_{1}+\lm_{2})\ep(\lm_{1}+\lm_{3})\ep(\lm
_{2}+\lm_{3})\ep(\lm_{1}+\lm_{2}+\lm_{3}).
\]

We set
\[
\bar{\Lm}=\Lm/2\Lm
\]
and denote the canonical projection from $\Lm$ to $\bar{\Lm}$ by
$\lm\mapsto\blm$. \ We regard $\bar{\Lm}$ as a vector space
over the field $\mathbb{Z}_{2}=\{0,1\}$. \ Since $s^{2}\in\cZ$ for
homogeneous $s\in\cA^{-}$ and since $\cA^{-}$ generates
$\cA$, we have $2\Lm\subset\Gm$, the support of $\cZ$.
\ Thus, $\ep$, $\al$, $\beta$, and $\mu$ are unchanged by adding an
element $\gm\in2\Lm$ to any one of their arguments. \ So $\ep$,
$\al$ $\beta$ and $\mu$ induce maps $\ep:\bar{\Lm}%
\rightarrow\{\pm1\}$, $\beta:\bar{\Lm}\times\bar{\Lm}\rightarrow F$,
$\al:\bar{\Lm}\times\bar{\Lm}\times\bar{\Lm}\rightarrow F$ and
$\mu:\bar{\Lm}\times\bar{\Lm}\times\bar{\Lm}\rightarrow\{\pm1\}$
with arguments from the $\mathbb{Z}_{2}$-vector space $\bar{\Lm}$.

\begin{lemma}
\label{lem:multid} Let $\cA$ be a structurable $\Lm$-torus of
class I and let $\lm_{i}\in\Lm$ for $1\leq i\leq4$. Then:

(a)\ms The commutator $\beta$ takes values in $\{\pm1\}$,
and
\[\beta(\lm_2,\lm_1) = \beta(\lm_{1},\lm_{2})=\ep(\lm_{1})\ep
(\lm_{2})\ep(\lm_{1}+\lm_{2}).\]

(b)\ms $\mu(\lm_{1},\lm_{2},\lm_{3})=\beta(\lm_{1}+\lm
_{2},\lm_{3})\beta(\lm_{1},\lm_{3})\beta(\lm_{2},\lm
_{3})$.

(c)\ms $\mu(\lm_{1},\lm_{2},\lm_{3})=\al(\lm_{1},\lm
_{2},\lm_{3})\al(\lm_{3},\lm_{1},\lm_{2})\al(\lm
_{2},\lm_{3},\lm_{1})$.

(d)\ms $\al(\lm_{2},\lm_{3},\lm_{4})\al(\lm_{1}%
+\lm_{2},\lm_{3},\lm_{4})^{-1}\al(\lm_{1},\lm
_{2}+\lm_{3},\lm_{4})\newline
{}\qquad\qquad\qquad\qquad\qquad\qquad\al(\lm_{1},\lm_{2}%
,\lm_{3}+\lm_{4})^{-1}\al(\lm_{1},\lm_{2},\lm_{3})=1$.

(e)\ms If $\blm_{1},\blm_{2},\blm_{3}$ are dependent
then $\mu(\lm_{1},\lm_{2},\lm_{3})=1$.

(f)\ms If $\ep(\lm_{i})=-1$ and $\beta(\lm_{i},\lm_{j})=1$
for some $i\neq j$ then $\al(\lm_{1},\lm_{2},\lm_{3})=1$

(g)\ms If some $\ep(\lm_{i})=-1$ then $\al(\lm_{1}%
,\lm_{2},\lm_{3})=\mu(\lm_{1},\lm_{2},\lm_{3})$ .

(h)\ms The associator $\al$ take values in
$\{\pm1\}$ and
$\al(\lm_{3},\lm_{2},\lm_{1}) = \al
(\lm_{1},\lm_{2},\lm_{3})$.

\end{lemma}

\begin{proof}
To simplify notation, write $\ep_{1}=\ep(\lm_{1})$,
$\ep_{12}=\ep(\lm_{1}+\lm_{2})$, $\beta_{1,2}%
=\beta(\lm_{1},\lm_{2})$, $\beta_{12,3}=\beta(\lm_{1}+\lm
_{2},\lm_{3})$, etc.

For (a), if $x_{i}\in\cA^{\lm_{i}}$, then
$x_{2}x_{1}=(x_{1}^{*}x_{2}^{*})^{*}=\ep_{1}\ep
_{2}\ep_{12}x_{1}x_{2}$.
Thus, $\beta_{1,2} = \ep_{1}\ep
_{2}\ep_{12} \in\{\pm1\}$, so $\beta_{2,1} = \beta_{1,2}$.

(b) follows by writing $\beta$ in terms of $\ep$ using (a).

For (c), we have $(x_{2}x_{3})x_{1}=\beta_{1,23}x_{1}(x_{2}x_{3})=\beta_{1,23}\al
_{1,2,3}(x_{1}x_{2})x_{3}$,
so
\[
\gm_{1,2,3}\gm_{3,1,2}\gm_{2,3,1}=1
\]
for $\gm_{1,2,3}=\beta_{1,23}\al_{1,2,3}$. But using (a), we have $\beta_{1,23}%
\beta_{3,12}\beta_{2,31}=\mu_{1,2,3}=\mu_{1,2,3}^{-1}$, so (c) follows.

For (d), we have
\begin{align*}
x_{1}(x_{2}(x_{3}x_{4}))  & = \al_{2,3,4}x_{1}((x_{2}x_{3})x_{4})\\
& = \al_{2,3,4}\al_{1,23,4}(x_{1}(x_{2}x_{3}))x_{4}\\
& = \al_{2,3,4}\al_{1,23,4}\al_{1,2,3}((x_{1}x_{2})x_{3})x_{4}\\
& = \al_{2,3,4}\al_{1,23,4}\al_{1,2,3}\al_{12,3,4}^{-1}(x_{1}%
x_{2})(x_{3}x_{4})\\
& = \al_{2,3,4}\al_{1,23,4}\al_{1,2,3}\al_{12,3,4}^{-1}%
\al_{1,2,34}^{-1}x_{1}(x_{2}(x_{3}x_{4})).
\end{align*}

For (e), let $U=\mathbb{Z}_{2}^{3}$ have basis $u_{i}$ and let $\ph$ be the
linear map with $\ph(u_{i})=\blm_{i}$. \ Since $\ep(0)=1$,
we see that
\[
\mu_{1,2,3}=\prod_{u\in U}\ep(\ph(u))\text{.}%
\]
If $\blm_{1},\blm_{2},\blm_{3}$ are dependent, then
each $\ph(u)$ occurs $\ord{\ker(\ph)} =2,4,$ or $8$ times,
so $\mu_{1,2,3}=1$.

(f) is a restatement of Lemma \ref{lem:associate}(a).

Before continuing, we observe that applying the involution to the defining equation
for $\al$ yields
\begin{equation}
\label{eq:alphaout}
\al_{1,2,3}=\al_{3,2,1}^{-1}.
\end{equation}

Now for (g) suppose that $\ep_i = -1$.
If $\beta_{i,j}=1$ for some $j\neq i$, then each factor on the right
side of (c) is $1$, so $\mu_{1,2,3}=1=\al_{1,2,3}$. \ Thus, we may assume
that $\beta_{i,j}=\beta_{i,k}=-1$ for $\{i,j,k\}=\{1,2,3\}$. \ Note that
\[
(x_{i},x_{j},x_{k})=(1-\al_{i,j,k})(x_{i}x_{j})x_{k},
\]
whereas
\[
(x_{i},x_{j},x_{k})=-(x_{j},x_{i},x_{k})=-(1-\al_{j,i,k})(x_{j}x_{i}%
)x_{k}=(1-\al_{j,i,k})(x_{i}x_{j})x_{k},
\]
so $\al_{i,j,k}=\al_{j,i,k}=\al_{k,i,j}^{-1}$. \ Replacing, $1,2,3$
in (c) by $i,j,k$ gives $\al_{j,k,i}=\mu_{i,j,k}=\mu_{1,2,3}$.
Interchanging the roles of $j$ and $k$, we also have $\al_{i,k,j}%
=\al_{k,i,j}$ and $\al_{k,j,i}=\mu_{1,2,3}$. \ Since $\al
_{p,q,r}=\al_{r,q,p}^{-1}$, we see $\al_{p,q,r}=\mu_{1,2,3}$ for all
$\{p,q,r\}=\{i,j,k\}$. \ In particular, $\al_{1,2,3}=\mu_{1,2,3} $.

For (h) it follows from (g) that $\al_{1,2,3}=\pm1$ if $\ep_{1}=-1$. \ If
$\ep(\lm)=1$, by Lemma \ref{lem:associate}(b) we can write
$\lm=\lm_{1}+\lm_{2}$ with $\ep_{1}=\ep_{2}=-1$.
\ Now (d) shows $\al(\lm,\lm_{3},\lm_{4})=\pm1$.  Thus
by \eqref{eq:alphaout} we have (h).
\end{proof}

\begin{notation} We now introduce some notation and terminology that we
will apply both to the pair $(\Lm,\ep)$ and the pair
$(\bar{\Lm},\ep)$.\footnote{The motivation for the terminology
comes from the case when $(M,\ep)=(\bar{\Lm},\ep)$ and
$\cA$ is associative. In that case, $\ep$ is a multiplicative
version of a quadratic form (see Lemma \ref{lem:quadratic}) and the
terminology is standard.} To do this efficiently, suppose that
$(M,\ep)$ is a pair consisting of an abelian group $M$ and a map $\ep
:M\rightarrow\{\pm1\}$ so that
\begin{equation}
\label{eq:pairdef}
\ep(0)=1\quad\text{and}\quad\ep\text{ is constant on cosets of
}2M.
\end{equation}
This last condition means that $\ep(u)=\ep(v)$ whenever
$u-v\in2M$. \ Let $\beta:M\times M\rightarrow\{\pm1\}$ be defined by
\[
\beta(u,v)=\ep(u+v)\ep(u)\ep(v).
\]
If $(M',\ep)$ is another pair satisfying \eqref{eq:pairdef},
we say $(M,\ep)$ is \textit{isomorphic to} $(M',\ep
')$, written $(M,\ep)\simeq(M',\ep')$,
if there is a group isomorphism $\ph:M\rightarrow M'$ with
$\ep'(\ph(u))=\ep(u)$ for $u\in M$. \ If $u\in M$, we
say that $u$ is \textit{anisotropic} if $\ep(u)=-1$. \ If $N$ is a
subgroup of $M$, we say that $N$ is anisotropic if all vectors $u\in N$ that
are not in $2N$ are anisotropic. \ Furthermore we say that the pair
$(M,\ep)$ is anisotropic if $M$ an anisotropic subgroup of itself.
\ If $u,v\in M$, we say $u$ is \textit{orthogonal} to $v$, written $u\perp v$,
if $\beta(u,v)=1$.

Finally suppose that $M_{1},\dots,M_{k}$ are subgroups of $M$. \ If
$M=M_{1}\oplus\dots\oplus M_{k}$ and $\beta(M_{i},\sum_{j\neq i}M_{j})=1$ for
all $i$, we write $M=M_{1}\perp\dots\perp M_{k}$. \ If $M=M_{1}\oplus
\dots\oplus M_{k}$ it is easy to check that
\begin{align*}
M=M_{1}\perp\dots\perp M_{k} &  \iff\ep(\sum_{i}u_{i})=\prod
_{i}\ep(u_{i})\text{ for }u_{i}\in M_{i}\\
&  \iff\beta(M_{i},\sum_{j>i}M_{j})=1\text{ for all }i<k.
\end{align*}
\end{notation}

\begin{lemma}
\label{lem:lift} If $V_{1},\dots,V_{k}$ are subgroups of $\bar{\Lm}$ with
$\bar{\Lm}=V_{1}\perp\dots\perp V_{k}$, then there exist subgroups
$\Lm_{1}$,\dots,$\Lm_{k}$ of $\Lm$ so that $\Lm=\Lm
_{1}\perp\dots\perp\Lm_{k}$ and $\bar{\Lm}_{i}=V_{i}$ for
$i=1,\dots,k$.
\end{lemma}

\begin{proof}
We can assume without loss of generality that $\Lm=\mathbb{Z}^{n}$ and
identify $\bar{\Lm}=\mathbb{Z}_{2}^{n}$. \ We first note that
$\GL(n,\mathbb{Z})\rightarrow \GL(n,\mathbb{Z}_{2})$ is surjective. \ Indeed,
$\GL(n,\mathbb{Z}_{2})=\SL(n,\mathbb{Z}_{2})=\Elem(n,\mathbb{Z}_{2})$, the
elementary group, and $\Elem(n,\mathbb{Z})\rightarrow \Elem(n,\mathbb{Z}_{2})$ is surjective.

For $1\leq i\leq k$, let $d_{i}=\dim(V_{i})$ and choose a basis $\bar{B}_{i} $
for the subspace $V_{i}$ of $\bar{\Lm}$. \ Let $\bar{A}$ be the element of
$\GL(n,\mathbb{Z}_{2})$ whose rows are the elements of $\bar{B}_{1} $, \dots,
$\bar{B}_{k}$ in order. \ Choose $A\in \GL(n,\mathbb{Z})$ with $A\rightarrow
\bar{A}$. \ Let $B_{1}$ be the set of the first $d_{1}$ rows of $A$, let
$B_{2}$ be the set of the next $d_{2}$ rows of $A$, and so on.  Let
$\Lm_{i}$ be the subgroup of $\Lm$ generated by $B_{i}$ for
$i=1,\dots,k$. \ Clearly, $\Lm_{1},\dots,\Lm_{k}$ have the properties claimed.
\end{proof}

\begin{lemma}
\label{lem:anisotropic} Let $\cA$ be a structurable $\Lm$-torus of
class I. Then the dimension of an anisotropic subspace of $\bar{\Lm}$ is at
most $3$. Moreover, for $\lm_1,\lm_2,\lm_3\in\Lm$,
we have $\al(\lm_{1},\lm_{2},\lm_{3}
)=\ep(\lm_{1})=\ep(\lm_2) = \ep(\lm_3) = -1$
if and only if $\blm_{1},\bar{\lm
}_{2},\blm_{3}$ is a basis for an anisotropic subspace of
$\bar{\Lm}$.
\end{lemma}

\begin{proof}
If $\blm_{1},\blm_{2},\blm_{3}$ is a basis for an
anisotropic subspace, then
$\al(\lm_{1},\lm_{2},\lm_{3})=\mu(\lm_{1},\lm
_{2},\lm_{3})=-1$
by Lemma \ref{lem:multid}(g) and the definition of $\mu$. \ Thus,
$\bar{\Lm}$ does not contain an anisotropic subspace of dimension $4$,
since otherwise, using Lemma \ref{lem:multid}(d), we would have $(-1)^{5}=~1$.

It remains to show that if $\al(\lm_{1},\lm_{2},\lm
_{3})=\ep(\lm_{i})=-1$, then $\blm_{1},\blm%
_{2},\blm_{3}$ is a basis for an anisotropic subspace of
$\bar{\Lm}$. \ Indeed Lemma \ref{lem:multid}(f) shows $\beta(\lm
_{i},\lm_{j})=-1$ for $i\neq j$, so $\ep(\lm_{i}+\lm
_{j})=-1$. \ Also, by Lemma \ref{lem:multid}(g), $\mu(\lm_{1},\lm
_{2},\lm_{3})=-1$, so $\ep(\lm_{1}+\lm_{2}+\lm
_{3})=-1$ by the definition of $\mu$. \ Thus, $\blm_{1},\bar{\lm
}_{2},\blm_{3}$ span an anisotropic subspace $V$. \ Moreover,
$\dim(V)=3$ by Lemma \ref{lem:multid}(e).
\end{proof}

\begin{lemma}
\label{lem:quadratic} Let $\cA$ be a structurable $\Lm$-torus of
class I and suppose that $\cA$ is associative. \ Let ${\hat
{\ep}}:\bar{\Lm}\rightarrow\mathbb{Z}_{2}$ and ${\hat{\beta}}%
:\bar{\Lm}\times\bar{\Lm}\rightarrow\mathbb{Z}_{2}$ be the additive
versions of the maps $\ep$ and $\beta$ on $\bar{\Lm}$; that is,
define ${\hat{\ep}}:\bar{\Lm}\rightarrow\mathbb{Z}_{2}$ and
${\hat{\beta}}:\bar{\Lm}\times\bar{\Lm}\rightarrow\mathbb{Z}_{2}$ by
\[
\ep(v)=(-1)^{{\hat{\ep}}(v)}\quad\text{and}\quad\beta
(v_{1},v_{2})=(-1)^{{\hat{\beta}}(v_{1},v_{2})}%
\]
for $v,v_{1},v_{2}\in\bar{\Lm}$. Then ${\hat{\ep}}$ is a
quadratic form over $\mathbb{Z}_{2}$ with linearization ${\hat{\beta}}$.
(A quadratic form over $\mathbb{Z}_{2}$ is defined exactly as in
\S \ref{sec:composition} even though $\mathbb{Z}_{2}$
has characteristic~2.)
\end{lemma}

\begin{proof}
By Lemma \ref{lem:multid}(a), we have
${\hat{\beta}}(v_{1},v_{2})={\hat{\ep}}(v_{1}+v_{2})+{\hat{\ep}}(v_{1})+{\hat{\ep}}(v_{2}).$
In other words ${\hat{\beta}}$ is the linearization of ${\hat{\ep}}$.  Also
by Lemma \ref{lem:multid}(b) and (c), we have $\beta(v_{1}+v_{2},v_{3})
=\beta(v_{1},v_{3})\beta(v_{2},v_{3})$, so
${\hat{\beta}}(v_{1}+v_{2},v_{3})={\hat{\beta}}(v_{1},v_{3})+{\hat{\beta}%
}(v_{2},v_{3})$.
Hence ${\hat{\ep}}$ is a quadratic form.
\end{proof}

\begin{notation}
In the next results, we will make use of two particular
pairs $(M,\ep)$ that satisfy the conditions \eqref{eq:pairdef} above.
\ If $r\geq0$ we let $(\mathbb{Z}^{r},\ep_{-})$ be the anisotropic
pair with underlying group $\mathbb{Z}^{r}$. \ In other words
\[
\ep_{-}(\lm)=
\begin{cases}
\hphantom{-}1 & \text{if }\lm\in2\mathbb{Z}^{r}\\
-1 & \text{otherwise}.
\end{cases}
\]
(If $r=0$, $\mathbb{Z}^{r}=0$ and $\ep_{-}=1$.) We let $(\mathbb{Z}%
^{2},\ep_{0})$ be the pair with
\[
\ep_0(\lm_{1})=\ep_0(\lm_{2})=1\quad\text{and}%
\quad\ep_0(\lm_{1}+\lm_{2})=-1,
\]
where $\lm_{1},\lm_{2}$ is the standard basis for $\mathbb{Z}^{2}$.
\end{notation}

\begin{proposition}
\label{prop:assoc} Let $\cA$ be a class I structurable $\Lm$-torus.
If $\cA$ is associative, then there exist subgroups $\Lm
_{1}',\dots,\Lm_{k}',\Lm'',\Lm
'''$ of $\Lm$ so that
\[
\Lm=\Lm_{1}'\perp\ldots\perp\Lm_{k}'\perp
\Lm''\perp\Lm''',
\]
where $k\geq0$,
\begin{align*}
&(\Lm_{i}',\ep)  \simeq(\mathbb{Z}^{2},\ep_{-})
\text{ for }1\leq i\leq k,\\
&(\Lm'',\ep)  \simeq 0,
(\mathbb{Z},\ep_{-})
\text{ or }(\mathbb{Z}^{2},\ep_{0})
\text{ with }(\Lm'',\ep)\simeq(\mathbb{Z},\ep_{-})
\text{ if } k=0,\\
&\Lm'''  \subset\Gm.
\end{align*}
\end{proposition}

\begin{proof}
Let ${\hat{\ep}}:\bar{\Lm}\rightarrow\mathbb{Z}_{2}$ and
${\hat{\beta}}:\bar{\Lm}\times\bar{\Lm}\rightarrow\mathbb{Z}_{2}$ be
defined as in Lemma \ref{lem:quadratic}. The classification of quadratic
forms over $\mathbb{Z}_{2}$ (see \cite[Chapitre I, \S 16]{D}
or Remark \ref{rem:quadraticforms}
below) shows that
\begin{equation*}
\label{eq:Lmtdecomp}
\bar{\Lm}=V_{1}'\perp\ldots\perp V_{k}'\perp
V''\perp V'''
\end{equation*}
where $k\geq0$, each $V_{i}'$ is anisotropic of dimension $2$,
$V''$ is either $0$, anisotropic of dimension $1$, or has a basis
$v_{1},v_{2}$ with ${\hat{\ep}}(v_{1})={\hat{\ep}}(v_{2})=0$
and ${\hat{\beta}}(v_{1},v_{2})=1$, and ${\hat{\ep}}(V^{\prime
\prime\prime})=0$. \ Moreover, since $\cA$ has class I, $\bar{\Lm
}$ is spanned by anisotropic vectors. \ Thus, if $k=0$, $V''$
must be anisotropic. \ The result now follows from Lemma~\ref{lem:lift}.
\end{proof}

\begin{remark}
\label{rem:quadraticforms} The classification of quadratic forms over
$\mathbb{Z}_{2}$ can be done as follows.  Suppose that ${\hat{\ep}}$
is an arbitrary quadratic form on a finite dimensional space $V$ over
$\mathbb{Z}_{2}$, and let ${\hat{\beta}}$ be the linearization of ${\hat{\ep}}
$. \ Since ${\hat{\beta}}$ is an alternating bilinear form, $V$ is the
orthogonal direct sum, relative to ${\hat{\beta}}$, of copies of $F$ with
basis $x$ and copies of the hyperbolic plane $H$ with basis $x,y$ and
${\hat{\beta}}(x,y)=1$. \ The quadratic form ${\hat{\ep}}$ is
determined by specifying $(F,{\hat{\ep}}(x))$ or $(H,{\hat
{\ep}}(x),{\hat{\ep}}(y))$ for each copy.  We have
\begin{align*}
\text{(a)} &  \text{\ }(F,1)\perp(F,1)\simeq(F,1)\perp(F,0),\\
\text{(b)} &  \text{\ }(H,1,0)\simeq(H,0,0),\\
\text{(c)} &  \text{ }(H,0,0)\perp(H,0,0)\simeq(H,1,1)\perp(H,1,1),\\
\text{(d)} &  \text{ }(H,0,0)\perp(F,1)\simeq(H,1,1)\perp(F,1).
\end{align*}
Indeed, for (a) we can start with a basis $x_{1}$, $x_{2}$ for the left hand
side, where $x_{1}$ is a basis for the first summand and $x_{2}$ is a basis
for the second.  The new basis $x_{1},x_{1}+x_{2}$ then gives a decomposition
as on the right hand side. Similarly (with obvious notation): for (b) use
$x+y,y$, for (c) use $x_{1}+y_{1},x_{1}+x_{2}+y_{2},x_{2}+y_{2},x_{1}%
+y_{1}+x_{2}$, and for (d) use $x_{1}+x_{2},y_{1}+x_{2},x_{2}$. \ It follows
from (a), (b), (c) and (d) that we can arrange to have no copies of $(H,1,0)$
and at most one copy of $(F,1)$ or $(H,0,0)$, but not both.
\end{remark}

\begin{proposition}
\label{prop:nassoc} Let $\cA$ be a class I structurable $\Lm
$-torus.  If $\cA$ is not associative, then there exist subgroups
$\Lm'$, $\Lm''$ and $\Lm^{\prime\prime\prime
}$ of $\Lm$ so that
\[
\Lm=\Lm'\perp\Lm''\perp\Lm^{\prime
\prime\prime}%
\]
where $(\Lm',\ep)\simeq(\mathbb{Z}^{3},\ep_{-})$,
$(\Lm'',\ep)\simeq(\mathbb{Z}^{k},\ep_{-})$
with $0\leq k\leq3$, and $\Lm'''\subset\Gm$.
\end{proposition}

\begin{proof}
In view of Lemmas \ref{lem:lift} and \ref{lem:anisotropic}, it suffices to
show that there are subspaces $V'$, $V''$ and $V'''$ of $\bar \Lm$ such that
$V$ and $V'$ are anisotropic of dimensions $3$ and $k$ respectively,
$V'''\subset\Gm/2\Lm$ and
\[
\bar{\Lm}=V'\perp V''\perp V'''.
\]

Since $\cA$ is not associative, there are $s_{i}\in\cA^{-}$
with $(s_{1},s_{2},s_{3})\neq0$ by Lemma \ref{lem:skewgen}(b). \ We can assume
$s_{i}$ is homogeneous with $s_{i}\in\cA^{\sg_{i}}$. \ Thus,
$\ep(\sg_{i})=-1$ and $\al(\sg_{1},\sg_{2},\sg_{3}%
)=-1$. \ Let $V'$ be the anisotropic subspace with basis
$\bar{\sg}_{1},\bar{\sg}_{2},\bar{\sg}_{3}$ as in Lemma \ref{lem:anisotropic}.

We next observe that if $\lm_{1},\lm_{2},\lm_{3}\in\Lm$ and some
$\ep(\lm_i) = -1$ then
\begin{equation}
\label{eq:orthog}
\lm_{1}\perp\lm_{3},\ \lm_{2} \perp\lm_{3}\implies(\lm_{1}+\lm_{2})\perp\lm_{3}
\end{equation}
Indeed we have $\mu(\lm_{1},\lm_{2},\lm_{3})=\al(\lm
_{1},\lm_{2},\lm_{3})=1$ by Lemma \ref{lem:multid}(g) and (f). Thus,
$\beta(\lm_{1}+\lm_{2},\lm_{3})=
\beta(\lm_{1},\lm_{3})\beta(\lm_{2},\lm_{3}) = 1$ by Lemma
\ref{lem:multid}(b), so $(\lm_{1}+\lm_{2})\perp\lm_{3}$.

Now let
\[
W=V^{\prime\perp}=\{\blm\in\bar{\Lm}:\beta(\bar{\lm
},V')=1\}.
\]
Then, it follows from \eqref{eq:orthog} that $W$ is a subspace of
$\bar{\Lm}$. We will show that
\begin{equation}
\label{eq:decomp1}
\bar{\Lm}=V'\perp W.
\end{equation}

Indeed, $V'\cap W=\{0\}$ since $\beta(\lm_1,\lm_2)=-1$
for any two distinct nonzero elements $\blm_1$,
$\blm_2$ in $V'$. Moreover, $\beta(V',W)=1$ by
definition of $W$.  So to prove \eqref{eq:decomp1} it remains to show that
$\bar{\Lm}=V'+W$. Since $S_{-}$ spans $\Lm$, we can do this
by showing that $\bar{\sg}\in V'+W$ for all $\bar{\sg}\notin
V'\cup W$ with $\ep(\sg)=-1$. Since $\bar{\sg}\notin
W$, we have $\beta(\sg,\lm_{3})=-1$ and hence $\ep(\lm_{3}+\sg)=-1$
for some $\blm_{3}\in V'$.
Let
$\bmu_1$, $\bmu_2$, $\bmu_3$
be a basis for $V'$.
Since $\al(\mu_1,\mu_2,\mu_3) = -1$,
Lemma \ref{lem:multid}(d) gives
\[\al(\mu_{2},\mu_{3},\sg)\al(\mu_{1}
+\mu_{2},\mu_{3},\sg)\al(\mu_{1},\mu_{2}+\mu_{3},\sg)
\al(\mu_{1},\mu_{2},\mu_{3}+\sg)= -1.\]
At least one of the factors on the left hand side of this equality must equal $-1$,
so we can replace $\sg$ by $\mu_3+\sg$ or rechoose the basis $\bmu_1$, $\bmu_2$, $\bmu_3$
of $V'$ to get
$\al(\mu_1,\mu_2,\sg) = -1$.
Thus,
by Lemma \ref{lem:anisotropic}, $X+\mathbb{Z}_{2}\bar{\sg}$ is anisotropic
where $X = \mathbb{Z}_{2}\bmu_1 + \mathbb{Z}_{2}\bmu_2$.
Since $V'+\mathbb{Z}_{2}
\bar{\sg}$ is not anisotropic by Lemma \ref{lem:anisotropic},
there is some $\bar \tau \in V'\setminus X$ with $\ep(\tau + \sg) = 1$.
Then Lemma \ref{lem:multid}(a) shows
\[\beta(\tau+\sg,\tau) = 1.\]
Next, if $i=1,2$, the space $\mathbb{Z}_{2}\bar\tau + \mathbb{Z}_{2}\bar\sg + \mathbb{Z}_{2}\bmu_i$
is not anisotropic since  $\ep(\tau + \sg) = 1$.
Therefore, by Lemma \ref{lem:anisotropic}, we have $\al(\tau,\sg,\mu_i) = 1$. So,
by Lemma \ref{lem:multid}(g),
$\mu(\tau,\sg,\mu_i) = 1$, and hence, by Lemma \ref{lem:multid}(b),
$\beta(\tau+\sg,\mu_i) = \beta(\tau,\mu_i)\beta(\sg,\mu_i)$.
But, since $V'$ and $X+\mathbb{Z}_{2}\bar{\sg}$ are anisotropic, we have $\beta(\tau,\mu_i)=-1$
and $\beta(\sg,\mu_i)=-1$.  So
\[\beta(\tau+\sg,\mu_i) = 1\]
for $i=1,2$. Since $\bar \tau$, $\bmu_1$, $\bmu_2$ span $V'$, we see that $\bar \tau + \bar \sg \in W$
by \eqref{eq:orthog}.  Thus,
$\bar \sg = \bar \tau + (\bar \tau + \bar \sg)\in V' + W$,  proving \eqref{eq:decomp1}.

Next since $\beta(V',W)=1$, we have $\ep(\lm+\tau)=\ep(\lm)\ep(\tau)$ for $\blm\in W$ and
$\bar{\tau}\in V'$. Thus, for $\blm_{i}\in W$ and $\bar\tau_i\in V'$,
we have
\begin{equation}
\label{eq:lmsg1}
\mu(\lm_{1}+\tau_{1},\lm_{2}+\tau_{2},\lm_{3}+\tau_{3})
=\mu(\lm_{1},\lm_{2},\lm_{3})\mu(\tau_{1},\tau_{2},\tau_{3})
\end{equation}
and
\begin{equation}
\label{eq:lmsg2}
\beta(\lm_{1}+\tau_{1},\lm_{2}+\tau_{2})
=\beta(\lm_{1},\lm_{2})\beta(\tau_{1},\tau_{2}).
\end{equation}

Now let
\[
V'''=\{\blm\in W:\ep(\lm)=1\}.
\]

We first argue that
\begin{equation}
\label{eq:VW}
\beta(V''',V)=1.
\end{equation}
Indeed if $\blm_{1}\in V'''$ and $\blm_{2}\in W$,
we have
\[
\mu(\lm_{1}+\sg_{1},\lm_{2}+\sg_{2},\sg_{3})
=\mu(\lm_{1},\lm_{2},0)\mu(\sg_{1},\sg_{2},\sg_3)=(1)(-1) = -1
\]
by \eqref{eq:lmsg1}, and
\[
\ep(\lm_{1}+\sg_{1})=\ep(\lm_{1})\ep
(\sg_{1})\beta(\lm_{1},\sg_{1})=(1)(-1)(1)=-1.
\]
(Here $\bar\sg_1,\bar\sg_2,\bar\sg_3$ is the basis for $V'$ chosen
at the beginning of this proof.) So by Lemma \ref{lem:multid}(f) and
(g), we have $\beta(\lm_{1}+\sg _{1},\lm_{2}+\sg_{2})=-1$. \ Thus,
by \eqref{eq:lmsg2},
$\beta(\lm_{1},\lm_{2})\beta(\sg_{1},\sg_{2})=-1$, so
\[\beta(\lm_{1},\lm_{2})=1.\]
Then, if $\blm_{1}\in V'''$, $\blm_{2}\in W$ and $\bar{\tau}\in V'$, we have using
\eqref{eq:lmsg2} that $\beta(\lm_{1},\lm_{2}+\tau)=\beta(\lm_{1},\lm_{2})\beta(0,\tau)= (1)(1) = 1$
which shows \eqref{eq:VW}.

Next if $\blm_{1},\blm_{2}\in V'''$ then
$\beta(\lm_{1},\lm_{2})=1$ by (\ref{eq:VW}), so $\ep
(\lm_{1}+\lm_{2})=\ep(\lm_{1})\ep(\lm_{2})=1$. \ This
shows that $V'''$ is a subspace of $W$.

We next show that
\begin{equation}
\label{eq:subGamma}
V'''\subset\Gm/2\Lm
\end{equation}
For this suppose that $\blm_{1}\in V'''$. Then,  by (\ref{eq:VW}),
$\beta(\lm_{1},\Lm)=1$.  Hence, if $x\in\cA^{\lm_{1}}$,
we have $[x,\cA]=0$. So, by Lemma \ref{lem:skewgen}(c), $x$
associates with any two elements of $\cA$. \ Therefore, since
$\ep(\blm_{1})=1$, we have $\lm_{1}\in\Gm$, proving
\eqref{eq:subGamma}.

Now let $V''$ be a maximal anisotropic subspace of $W$. It
remains to show that $W=V''\perp V'''$.
Certainly $V''\cap V'''=0$ and, by
(\ref{eq:VW}), $\beta(V'',V''')=1$. \ So we
only have to show that $W=V''+V'''$. \ To do
this let $\blm\in W$ with $\blm\notin V''\cup
V'''$. \ Then $\ep(\lm)=-1$. \ Also,
$\mathbb{Z}_{2}\blm+V''$ is not anisotropic, so
$\ep(\lm+\sg)=1$ for some $\bar{\sg}\in V''$.
\ Thus, $\blm=-\bar{\sg}+(\blm+\bar{\sg})\in
V''+V'''$.
\end{proof}

The following is the main result of this section:

\begin{theorem}
\label{thm:classI}
Let $\cA$ be a structurable $\Lm$-torus of class I.  Then:

(a)\ms If $\cA$ is associative, then $\cA$ is isograded
isomorphic to
\begin{equation*}
\label{eq:isom1}
\cA_{1}\ot\ldots\ot\cA_{k}\ot\cA_{k+1}\ot \cP(r)
\end{equation*}
where $k\geq0,r\geq0$, and
\begin{align}
\label{eq:mainas}
\cA_{i}  & =\Cay(2)\,\text{for }1\leq i\leq k,
\\
\label{eq:mainas'}
\cA_{k+1}  & =F\text{, }\Cay(1)\text{, or }\Cay_*(2)\text{ with }\cA_{k+1}=\Cay(1)\text{ if }k=0.
\end{align}

(b)\ms If $\cA$ is not associative, then $\cA$ is isograded
isomorphic to%
\begin{equation*}
\label{eq:isom2}
\Cay(3)\ot\Cay(k)\ot \cP(r)
\end{equation*}
where $0\leq k\leq3$ and $r\geq0$.

Conversely, each of the graded algebras with involution in (a)
and (b) are class I structurable tori.
\end{theorem}

\begin{proof}
We prove (a) and (b) together. First, we have orthogonal decompositions of
$\Lm$ as in Theorems \ref{prop:assoc} and \ref{prop:nassoc}. By Remark
\ref{rem:subtori}, we obtain the
$M$-tori $\cA^{M}$ with $M=\Lm_{1}',\dots,\Lm
_{k}',\Lm'',\Lm'''$, if
$\cA$ is associative, and $M=\Lm',\Lm^{\prime\prime
},\Lm'''$, if $\cA$ is not associative.
\ Repeated use of Lemma \ref{lem:tensor} now yields the isomorphisms%
\begin{align}
\label{eq:associsom}
& \cA\simeq_{\Lm}\cA^{\Lm_{1}'}\ot
\ldots\ot\cA^{\Lm_{k}'}\ot\cA^{\Lm''}\ot\cA^{\Lm'''}\\
\label{eq:nonassocisom}
& \cA\simeq_{\Lm}\cA^{\Lm'}\ot
\cA^{\Lm''}\ot\cA^{\Lm^{\prime
\prime\prime}},
\end{align}
in the two cases.  Indeed, if $\cA$ is associative, the argument for
this is clear.  If $\cA$ is not associative, we first get
$\cA\simeq_{\Lm}\cA^{\Lm'\oplus\Lm
''}\ot\cA^{\Lm'''}$ by Lemma
\ref{lem:tensor}, since $\Lm'''\subset \Gm$.  Next $[\cA^{\Lm'},\cA%
^{\Lm''}]=0$ since $\beta(\Lm',\Lm
'')=1$. \ Also since $\Lm'$ is spanned by anisotropic
vectors and $\beta(\Lm',\Lm'')=1$, we have
$\al(\Lm',\Lm'',\Lm)=\al(\Lm
,\Lm',\Lm'')=1$ by Lemma \ref{lem:multid}(f).
\ So we have $(\cA^{\Lm'},\cA^{\Lm
''},\cA)=(\cA,\cA^{\Lm'},\cA^{\Lm''})=0$
and we can apply Lemma
\ref{lem:tensor} to get $\cA^{\Lm'\oplus\Lm''}
\simeq_{\Lm'\oplus\Lm''}\cA^{\Lm'}\ot\cA^{\Lm''}$
which gives the isomorphism (\ref{eq:nonassocisom}).

To obtain our conclusions in (a) and (b),
it now suffices to show that if $M$ is a subgroup of $\Lm$, then
$\cA^{M}$ is isograded isomorphic to

\begin{itemize}
\item[(i)] $\Cay(k),$ if $(M,\ep)\simeq(\mathbb{Z}%
^{k},\ep_{-})$ with $1\leq k\leq3$,

\item[(ii)] $\Cay_*(2),$ if $(M,\ep)\simeq(\mathbb{Z}%
^{2},\ep_{0})$ and $\cA^{M}$ is associative,

\item[(iii)] $\cP(r)$, if $\ep(M)=1$ and
$M\subset \Gm$.
\end{itemize}

\noindent In cases (ii) and (iii) $\cA^{M}$ is associative and the
conclusions are well known and easy to verify (using an argument similar to
the one below for (i)).

So we consider only the case (i) and suppose that
$(M,\ep)\simeq(\mathbb{Z}^{k},\ep_{-})$ with $1\leq k\leq3$.
Choose a basis
$B=\{\lm_{1},\dots,\lm_{k}\}$ for $M$ and $0\neq x_{i}\in
\cA^{\lm_{i}}$. \ We get a basis $\{x^{\lm}\}_{\lm\in M}$
for $\cA^{M}$ by choosing for each $\lm\in M$ an element
$x^{\lm}\in\cA^{\lm}$ as follows: Write $\lm$ in a fixed
way as
\begin{equation}
\label{eq:lmexp}
\lm=\lm_{i_{1}}+\ldots+\lm_{i_{p}}-\lm_{j_{1}}-\ldots
-\lm_{j_{q}},
\end{equation}
and choose $x^{\lm}$ to be a product of $x_{i_{1}},\ldots,x_{i_{p}%
},x_{j_{1}}^{-1},\ldots,x_{j_{q}}^{-1}$ in some fixed order and association.
A change in the expression \eqref{eq:lmexp} for $\lm$ or in the order and
association of the product changes $x^{\lm}$ by a sign which is determined
by $\al$ and $\beta$ on $M$. Thus, the multiplication table for
$\cA^{M}$ is determined by $\al$ and $\beta$ while the involution
table is determined by $\ep$, $\al$ and $\beta$. \ But $\beta$ and
$\mu$ are determined by $\ep$, and, by Lemma \ref{lem:multid}(e) and
(g), $\al=\mu$ on $\cA^{M} $. \ Thus, the multiplication and
involution tables for $\cA^{M} $ are determined by $\ep$.
This argument shows that if $\cB$ is any other class I structurable
$N$-torus with $(M,\ep)\simeq(N,\ep)$, then $\cA%
^{M}\ig \cB$.  In particular, $\cA^{M}\ig \Cay(k)$.

Conversely, suppose that $\cA$ is a $\Lm$-graded algebra with involution $\cA$
satisfying the conclusions in (a) or (b).  In case (b),
$\cA$ is a structurable torus of class I by Proposition \ref{prop:constructionI}.
In case (a), $\cA$ is associative, hence structurable, and hence a structurable torus by
Remark \ref{rem:tensor}(a).  Moreover,
one easily checks using \eqref{eq:mainas} and \eqref{eq:mainas'}
that $\cA$ is generated by its skew-elements.
\end{proof}

\begin{remark}
\label{rem:assoctori}
If $\cA$ is an associative
$\Lm$-torus with involution (not necessarily of class I), then
the arguments in this section with almost no change show that $\cA$ is
isograded isomorphic to
\[\cA_{1}\ot\ldots\ot\cA%
_{k}\ot\cA_{k+1}\ot \cP(r),\]
where $k\ge 0$, $r\ge 0$, $\cA_{i}   =\Cay(2)$ for $1\le i \le k$,
and $\cA_{k+1}  =\Cay(0)$, $\Cay(1)$ or $\Cay_*(2)$.
(In other words we drop the restriction in Theorem \ref{thm:classI} on $\cA_{k+1}$
when $k=0$.)
This result has previously been obtained by Yoshii in \cite{Y2} using different methods.
It can also be deduced from K.-H.~Neeb's recent
classification of rational associative tori (without involution) \cite{Neeb} along with
(a)-(d) in Remark \ref{rem:quadraticforms}.
\end{remark}

\section{Structurable tori of class II}
\label{sec:classII}

In  this section we briefly recall from \cite{AY} the
construction and classification of structurable tori of class II.

We first recall a construction of structurable algebras from hermitian forms.
Suppose that $\cB$ is an associative algebra with involution $*$, and let
$\cX$ be a left $\cB$-module.
Suppose that $\hermk:\cX\times \cX \to \cB$ is a hermitian form over $\cB$;
that is $$\hermk(au,v) = a\hermk(u,v),\quad \hermk(u,a v) =
\hermk(u,v)a^*$$
and
\[\hermk(v,u) = \hermk(u,v)^*\]
for $u,v\in \cX$ and $a\in \cB$.  Let
\[\cA(\hermk) := \cB\oplus \cX\]
with product and involution
defined respectively by
\[(a,u)(b,v) = (a\cdot_{\text{op}}b + \hermk(u,v), av + b^* u)
\andd
(a,u)^* = (a^*,u)
\]
for $a,b\in\cB$ and $u,v\in \cX$, where $a\cdot_{\text{op}}b := ba$
is the product in the opposite algebra $\cB^\text{op}$ of $\cB$.%
\footnote{The product on $\cA(\hermk)$ defined in \cite{A1} and
\cite{AY} is the opposite of the product defined here.
However, the involution on $\cA(\hermk)$ is an isomorphism
of $\cA(\hermk)$ with $\cA(\hermk)^\text{op}$  as algebras with involution.}
The algebra with involution $\cA(\hermk)$ is a structurable algebra called the
\emph{structurable algebra associated with $\hermk$}
\cite[\S 8]{A1}.

\begin{example}
\label{ex:herm}
Let $M$ be a subgroup of
$\Lm$ such that $2\Lm \subset M$,
and let $\cB$ be an associative $M$-torus with involution $*$.
We regard $\cB$ as a $\Lm$-graded algebra with involution
by setting $\cB^\lm = 0$
for $\lm\in \Lm\setminus M$.
Let $m\ge 1$; let $\rho_1,\dots,\rho_m$ be elements of $\Lm\setminus M$
such that
\[
\Lm = \langle M,\rho_1,\dots, \rho_m\rangle,\quad \rho_i-\rho_j\notin M
\text{ for } i\ne j, \andd 2\rho_i\in S_+(\cB) \text{ for all } i;
\]
and let $b_1,\dots,b_m$ be elements of $\cB$ such that
\[0\ne b_i\in \cB^{2\rho_i} \text{ for all $i$}.\]
We construct a graded hermitian form from this data as follows.
Let $\cX$ be the left $\cB$-module
that is free of rank $m$ with basis $v_1,\dots,v_m$,
and assign $\cX$ the $\Lm$-grading such that
$bx_i$ has degree $\sg + \rho_i$ if $b\in \cB^\sg$, $\sg\in M$
and $1\le i \le m$.  Then \textit{the module $\cX$ is $\Lm$-graded};
that is $\cB^\lm \cX^{\lm'} \subset \cX^{\lm+\lm'}$
for $\lm,\lm'\in\Lm$.
Let  $\hermk : \cX\times \cX \to \cB$ be the
hermitian form such that
\[\hermk(v_i,v_j) = \delta_{i,j} b_i\]
for $1\le i,j \le m$.
Then \textit{the hermitian form $\hermk$ is $\Lm$-graded};
that is,
$\hermk(\cX^\lm,\cX^{\lm'}) \subset \cB^{\lm+\lm'}$ for $\lm,\lm'\in\Lm$.
We call a graded hermitian form $\hermk$ constructed in this way (using some choice
of $\rho_1,\dots,\rho_m$ and $b_1,\dots,b_m$) a \textit{diagonal $\Lm$-graded
hermitian form over~$\cB$}.
Now $\cA(\hermk) = \cB \oplus \cX$ is a structurable $\Lm$-torus with $\Lm$-grading
extending the $\Lm$-grading on $\cB$ and
$\cX$ \cite[Example 4.6]{AY}.\footnote{In \cite{AY}, the $M$-torus $\cB$ is realized as a quantum torus with
involution (see \cite[Prop.~4.5]{AY}).}  We call
$\cA(\hermk)$ the
\emph{structurable $\Lm$-torus
associated with $\hermk$}.
\end{example}

In the next theorem, proved in \cite{AY}, we restrict our attention to
coordinate algebras $\cB$ that are class I associative tori with involution.
These are precisely the associative  tori with involution that
were classified in Theorem \ref{thm:classI}(a).

\begin{theorem}
\label{thm:classII}
If  $M$ is a subgroup of
$\Lm$ such that $2\Lm \subset M$,
$\cB$ is an class I associative $M$-torus with involution,
and $\hermk$ is a diagonal $\Lm$-graded hermitian form
over~$\cB$, then the structurable $\Lm$-torus $\cA(\hermk)$ associated
with $\hermk$ (as in Example \ref{ex:herm}) is of class II.
Moreover, any structurable $\Lm$-torus of class II is graded-isomorphic
to a structurable $\Lm$-torus $\cA(\hermk)$  obtained in this way.
\end{theorem}

\begin{proof}  For the first statement, since $\cB$ is of class I, we have
$\Lm_- = M$, $\cE = \cB^\text{op}$  and $\cW = \cX$
(using the Notation from \S\ref{sec:classes} for $\cA = \cA(\hermk)$).
So, since $\cX \ne 0$, it is clear that $\cA(\hermk)$ is of class II.
The second statement is Theorem~9.5
of~\cite{AY}.
\end{proof}

\section{Cubic forms and structurable algebras}
\label{sec:Cubicforms}

Before beginning our study of class III structurable tori,
we pause to introduce a new construction of structurable algebras.
As background for this construction, we mention that there are three
constructions of finite dimensional simple structurable algebras $\cA$
with $\cA_-$ one dimensional: (i) a construction of a $2\times 2$-matrix algebra with
coefficients from a cubic Jordan algebra \cite{K1}, \cite[\S 8]{A1}; (ii)
a Cayley Dickson process for doubling a degree 4 Jordan algebra \cite{AF};
(iii) a construction using a self-adjoint norm semisimilarity of
a cubic Jordan algebra over a quadratic extension
of the base field \cite[p.1861]{A2},
\cite[Chapter 8]{Se}.  In this section we adapt the third
of these constructions to the infinite dimensional setting.
However, rather than initially using a norm similarity,
we base our construction
on a pair $(h,N)$ consisting of a hermitian form and a
cubic form satisfying the adjoint identity.
We then see that the first two constructions (adapted to the infinite dimensional setting)
can be viewed as special cases of this construction.

Since $F$ is a field of characteristic $\neq 2$ or $3$, we
may define homogeneous polynomial maps of degree $n$ for $1\le n \le
4$ in the following simple fashion. Suppose  $K$  is a commutative,
associative algebra over $F$, and suppose $\cU$ and $\cV$ are
left $K$-modules.  If $1\le n \le 4$, a map $g : \cU \to \cV$  is
called a \textit{homogeneous polynomial map of degree $n$} over $K$
if there is a symmetric $K$-multilinear map $f : \cU^n \to \cV$ with
\[g(x) = \frac 1{n} f(x,\ldots,x)\]
for $x\in\cU$.
In that case, given $u,v\in \cU$, there exist unique elements
$c_{i,u,v}\in \cV$ for $1\le i\le n$ such that
\[g(v+au) = g(v) + \sum_{i=1}^n a^i c_{i,u,v},\]
for $a\in F$.  We define $\partial_{u}g|_{v} = c_{1,u,v}$,
in which case
\[\partial_{u}g|_{v} = \frac 1{(n-1)} f(u,v,\dots,v).\]
Thus, for $u\in \cU$, we have the map $\partial_{u}g : \cU \to \cV$.
If $n\ge 2$, $\partial_{u}g$ is a homogeneous polynomial map of
degree $n-1$; whereas, if $n=1$, $\partial_{u}g$ is the constant map
$x\mapsto g(u)$.
Finally, we have
\[\partial_{u_1} \dots \partial_{u_{n-1}} g|_{u_n}  = f(u_1,\dots,u_n)\]
for $u_1,\dots,u_n \in \cU$, which shows that we can recover $f$ from $g$.
We call $f$ the \textit{full linearization} of $g$.

If $K$ is a  commutative, associative algebra over $F$ and if
$\cU$ is a left $K$-module, a \textit{cubic} (resp.~\textit{quadratic})
form on $\cU$ over $K$
is a homogeneous polynomial map $g: \cU \to K$ of  degree~$3$ (resp.~2) over
$K$.  (This notion of quadratic form is clearly equivalent to the one
described in \S \ref{sec:composition}.)

Suppose now that $\cE$ is a commutative, associative algebra with
involution $*$ over~$F$ and $\cW$ is a left  $\cE$-module.
Let $N:\cW\rightarrow \cE$ be a cubic form over $\cE$ and let $h:\cW\times\cW\rightarrow\cE$
be a hermitian form (see \S\ref{sec:classII}) that is nondegenerate
(i.e. $h(\cW,u)  =0\Longrightarrow u=0$ for $u\in \cW$).
We say that $(h,N)$ has an
\textit{adjoint} if for each $v\in\cW$ there is an element $v^{\natural}%
\in\cW$ with%
\[
\partial_{u}N|_{v}=h(u,v^{\natural})\text{ for all }u\in\cW%
\text{.}
\]
Clearly, the nondegeneracy of $h$ shows that the adjoint map $v\mapsto v^{\natural}$ is
unique if it exists.\footnote{Nondegeneracy of $h$ is equivalent to the map
$\al:u\rightarrow h(\ ,u)$ being an injection of $\cW$ into its
dual space for $\cE $. If $\al$ is a bijection, then $(h,N)$
automatically has an adjoint.}

If $(h,N)$ has an adjoint, we let
\[
u\diamond v=(u+v)^{\natural}-u^{\natural}-v^{\natural},%
\]
in which case
\[
h(u,v\diamond w)=\partial_{v}\partial_{u}N|_{w}%
\]
is a symmetric trilinear form (the full linearization of $N$).
Moreover,
\[v^\natural = \frac 12 v\diamond v \andd N(u) = \frac 16 h(u,u\diamond u).\]
Since $h(u,v\diamond w)$ is symmetric in its arguments, $\diamond$ is a
\emph{bi-semilinear product}
on $\cW$ (that is $\diamond$ is biadditive and
$(au)\diamond v = u\diamond (av) = a^* (u\diamond v)$
for $u,v\in\cW$, $a\in \cE$). Hence the map $\natural:\cW\to \cW$ is \emph{semiquadratic}
(that is $v^\natural = \frac 12 v\diamond v$ where $\diamond$ is symmetric
and bi-semilinear).

Conversely, if $h$ is a nondegenerate hermitian
form on $\cW$ and $\diamond$ is a bi-semilinear product on $\cW$ with
$h(u,v\diamond w)$ symmetric, then $N(u)=\frac{1}{6}h(u,u\diamond u)$ defines
a cubic form and $(h,N)$ has adjoint $u^{\natural}=\frac{1}{2}(u\diamond u)$.

\begin{remark}
\label{rem:classicaladjoint}
In the special case when the involution $*$ on $\cE$ is trivial ($=\id$), pairs $(h,N)$ with an adjoint
play an important role in the study of cubic Jordan algebras \cite[Chapter II.4]{Mc}.  In that
case, we usually write $T$ for $h$, $\#$ for $\natural$, and $\times$ for $\diamond$ (thereby following the
usual notational conventions).
\end{remark}

Given $(h,N)$ with an adjoint, we define
\[\cA(h,N)=\cE\oplus\cW\]
with product and involution given respectively by
\begin{equation}
\label{eq:defAhN}
(a,v)(b,w)   =(ab+h(v,w),aw+b^{*}v+v\diamond w) \andd (a,v)^{*}   =(a^{*},v).
\end{equation}
We will write $a+v$ for $(a,v)$ in $\cA(h,N)$.%

We note that if $h$ is a nondegenerate hermitian form and $N=0$ then $(h,0)$
has an adjoint (the zero map) and, since $\cE$ is commutative, $\cA(h,N) = \cA(h)$ (see \S\ref{sec:classII}).

We now wish to determine when the algebra with involution $\cA(h,N)$ is structurable.

We say that
$(h,N)$ satisfies the \textit{adjoint identity} if $(h,N)$ has an
adjoint and
\begin{mylist}
\item[(ADJ)]  $(v^{\natural})^{\natural}=N(v)v$ for all $v\in \cW$.
\end{mylist}
If
$(h,N)$ satisfies the \textit{adjoint identity}
then, since $\ord{F} \geq5$,
$(h,N)$ also satisfies the polarizations \label{page:ADJ}
\begin{mylist}
\item[(ADJ1)]
$(w\diamond v)\diamond v^{\natural}=N(v)w+h(w,v^{\natural})v$,
\item[(ADJ2)]
$(w\diamond u)\diamond v^{\natural}+(w\diamond v)\diamond(u\diamond v)
=h(u,v^{\natural})w+h(w,u\diamond v)v+h(w,v^{\natural})u,$
\item[(ADJ3)]
$(w\diamond u)\diamond(x\diamond v)+(w\diamond
x)\diamond(u\diamond v)+(w\diamond v)\diamond(u\diamond x)
\\{} \qquad \qquad \qquad =h(u,x\diamond v)w+h(w,u\diamond x)v+h(w,u\diamond v)x+h(w,x\diamond v)u,$
\item[(ADJ4)] $(v\diamond w)^{\natural}+v^{\natural}\diamond
w^{\natural}=h(w,v^{\natural})w+h(v,w^{\natural})v.$
\end{mylist}
Also, $h(v\diamond w,u) = h(u,v\diamond w)^* = h(v,u\diamond w)^*$.  Hence,
applying $h(\ ,u)$ to (ADJ1) gives
$h(w,(u\diamond v^{\natural})\diamond v)=h(w,N(v)^{*}u)+h(w,h(u,v)v^{\natural})$,
so
\begin{mylist}
\item[(ADJ5)]
$(u\diamond v^{\natural})\diamond v=N(v)^{*}%
u+h(u,v)v^{\natural}$
\end{mylist}

\begin{theorem}
\label{thm:cubic} Suppose that $\cE$ is a commutative associative algebra
with involution, $N:\cW\rightarrow \cE$ is a cubic form, and $h:\cW\times\cW\rightarrow\cE$
is a nondegenerate hermitian form such that $(h,N)$ has an adjoint. If $(h,N)$ satisfies
the adjoint identity, then $\cA(h,N)$ is a structurable
algebra. If $\cE$ contains a skew element that is invertible in $\cE$,
then the converse holds.
\end{theorem}

\begin{proof}
We recall from \cite[Theorem 13]{A1} that $\cA$ is structurable
if and only if%
\begin{itemize}
\item[(i)]  $\cA$  is skew-alternative,
\item[(ii)]  $(a,b,c)-(c,a,b) = (b,a,c) - (c,b,a)$ for $a,b,c\in \cA_+$,
\item[(iii)] $D_{a^{2},a}(b)=0$ for $a,b\in \cA_+$,
\end{itemize}
where
\[
D_{a,b}(c)=\frac{2}{3}[[a,b],c]+(c,b,a)-(c,a,b)
\]
for $a,b,c\in\cA_{+}$.

We shall first show that (i) and (ii) hold
automatically for $\cA(h,N)$ without assuming the adjoint
identity. We note that $\cW$ is a bimodule for $\cE$ using
the product in $\cA(h,N)$, so%
\[
(\cE,\cE,\cW)= (\cE,\cW,\cE)=(\cW,\cE,\cE)= 0.
\]
Skew-alternativity reduces  to
\begin{align*}
(s,u,v)+(u,s,v)  &= h(su,v) + (su)\diamond v - sh(u,v)-s(u\diamond v)\\
&\quad\quad + h(us,v) + (us)\diamond v-h(u,sv)-u\diamond(sv) =0
\end{align*}
for $s\in\cA(h,N)_{-}=\cE_{-}$ and $u,v\in\cW$, which
holds since $us=-su$, $h$ is sesquilinear and $\diamond$ is  bi-semilinear.
Since $(a,b,c)^{*}=-(c,b,a)$ for $a,b,c\in\cA_{+}$, (ii) is equivalent to
$(a,b,c)-(c,a,b)\in\cA_{+}$ for $a,b,c\in\cA_{+}$. \ We note
that $\cE_{+}$ is in the centre of the algebra with involution
$\cA(h,N)$ and that $\cA(h,N)_{+}=\cE%
_{+}\oplus\cW$, so it suffices to check the $\cE$-component of
$(u,v,w)-(w,u,v)$ for $u,v,w\in\cW$.  Now, the $\cE%
$-component of $u(vw)$ is $h(u,v\diamond w)$ which is symmetric in $u,v,w$, as
are the $\cE$-components of $(uv)w=(w(vu))^{*}$ and $(u,v,w)$.
Thus, the $\cE$-component of $(u,v,w)-(w,u,v)$ is $0$.

We next reduce condition (iii) without using the adjoint identity. \ If
$a\in\cE_{+}\subset\cZ(\cA(h,N))$ and $x,y\in
\cA(h,N)_{+}$, then%
\[
D_{x,y}(a)=0,\quad D_{ax,x}=aD_{x,x}=0, \quad D_{a,x}=D_{x,a}=0.
\]
Let $x=(a,v)=a+v\in\cA(h,N)_{+}$ so $x^{2}=a^{2}%
+h(v,v)+2av+v\diamond v$. We see that
\[
D_{x^{2},x}=D_{2av+v\diamond v,v}=D_{v\diamond v,v}=2D_{v^{\natural},v}.
\]
Thus, (iii) reduces to%
\[
D_{v^{\natural},v}(w)=0
\]
for $v,w\in\cW$.

Since $\diamond$ is symmetric, we have $[v^{\natural},v]=h(v^{\natural
},v)-h(v,v^{\natural})=-3(N(v)-N(v)^{*})$,~so
\[
\frac{2}{3}[[v^{\natural},v],w]=-4(N(v)-N(v)^{*})w.
\]
Also, since $(wu)v=h(w\diamond u,v)+h(w,u)v+(w\diamond u)\diamond v$ and since
$h(w\diamond u,v)$ is symmetric in $u,v$, we see that
\begin{align*}
(w,v,&v^{\natural})-(w,v^{\natural},v)  =(wv)v^{\natural
}-(wv^{\natural})v+w[v^{\natural},v]\\
&=h(w,v)v^{\natural}+(w\diamond v)\diamond v^{\natural}
 -h(w,v^{\natural})v-(w\diamond v^{\natural})\diamond v
+3(N(v)-N(v)^{*})w.
\end{align*}
Thus,
\[
D_{v^{\natural},v}(w)=-(N(v)-N(v)^{*})w+h(w,v)v^{\natural}+(w\diamond
v)\diamond v^{\natural}-h(w,v^{\natural})v-(w\diamond v^{\natural})\diamond v
\]
and (iii) is equivalent to
\[
\text{(iii}'\text{) }N(v)w+h(w,v^{\natural})v-(w\diamond v)\diamond
v^{\natural}=N(v)^{*}w+h(w,v)v^{\natural}-(w\diamond v^{\natural})\diamond
v.
\]

If the adjoint identity holds, both sides of (iii$'$) are $0$ by
(ADJ1) and (ADJ5) and $\cA(h,N)$ is a structurable algebra.
Conversely, if $\cA(h,N)$ is a structurable algebra and
$s\in\cE_{-}$ is invertible in $\cE$, $w=v$ in (iii$'$) gives
\begin{itemize}
\item[(iii$''$)]
$4(N(v)v-(v^{\natural})^{\natural
})=N(v)^{*}v+h(v,v)v^{\natural}-(v\diamond v^{\natural})\diamond v$.
\end{itemize}
Replacing $v$ by $sv$ in (iii$''$) multiplies the left side by
$s^{4}$ and the right side by $-s^{4}$. Thus, both sides are $0$ and the
adjoint identity holds.
\end{proof}

We  will also need the following simple fact:

\begin{lemma}
\label{lem:centre}
Suppose that $\cE$ is a commutative associative algebra
with involution, $N:\cW\rightarrow \cE$ is a cubic form, and $h:\cW\times\cW\rightarrow\cE$
is a nondegenerate hermitian form such that $(h,N)$ has an adjoint.
If $\cE$ contains a skew element that is invertible in~$\cE$,
then $\cZ(\cA(h,N)) = \cE_+$.
\end{lemma}

\begin{proof}  Suppose that $s\in \cE_-$ is invertible in $\cE$.
If $x = a+v\in \cZ(\cA(h,N))$, then $a^* = a$ and
$0 = [x,s] = (a+v)s - s(a+v) = -2sv$, which forces $v = 0$. Conversely,
it is easy to check that  $\cE_+ \subset \cZ(\cA(h,N))$.
\end{proof}

We shall now relate $\cA(h,N)$ to two other constructions of
structurable algebras.

We first consider the Cayley-Dickson process for structurable algebras that was
studied in \cite{AF}.  Let $\cJ$ be a commutative algebra
with $1$ over a commutative associative $F$-algebra $K$. In practice,
$\cJ$ will be Jordan, but we do not require that immediately. \ Let
$x\rightarrow x^{\theta}$ be a period two $K$-linear map on $\cJ$ and
let $\mu\in K^{\times}$. The \textit{Cayley-Dickson} process (see
\cite[p.~200]{AF}) gives an algebra with involution
\[\CD(\cJ,\theta,\mu)=\cJ\oplus \cJ,\]
with $(a,b)$ written as $a+s_0b$,
\begin{equation}
\label{eq:CDmult}
(a+s_0b)(c+s_0d)   =(ac+\mu(bd^{\theta})^{\theta})+s_0(a^{\theta
}d+(b^{\theta}c^{\theta})^{\theta}),
\end{equation}
and
\begin{equation}
\label{eq:CDinvol}
(a+s_0b)^{*}  =a-s_0b^{\theta}.
\end{equation}
Note that $\cE=K[s_0]$ with
$s_0^{2}=\mu$ is a subalgebra of $\CD(\cJ,\theta,\mu)$.

\begin{lemma}
\label{lem:CD} Let  $K$ be a commutative and associative algebra which we
regard as an algebra with involution with the trivial involution.  Let
$N_\cV: \cV \to K$ be a cubic form on a $K$-module $\cV$ and let
$T : \cV \times\cV \to \cV$ be a nondegenerate symmetric
$K$-bilinear form such that $(T,N_\cV)$ satisfies the adjoint identity.
Denote the adjoint of $(T,N_\cV)$ by $\#$ and the linearization
of $\#$ by $\times$ (see Remark \ref{rem:classicaladjoint}).
Form the commutative algebra
\[\cJ=\cA(T,N_\cV) = K \oplus \cV\]
(with trivial involution)
and define $\theta : \cJ \to \cJ$ by
$(\al+v)^{\theta}=\al-v$ for $\al\in K,v\in\cV$.
Also form $\cE=K[s_0]$ with $s_0^{2}=\mu\in K^{\times}$
and $K$-linear involution $*$ with $s_0^{*}=-s_0$,
and form the left $\cE$-module $\cW=\cE\ot_{K}\cV$ with
$\cV$ identified with $1\ot\cV$.  Extend  $T$ to a hermitian form
$h:\cW\times \cW\rightarrow\cE$, extend $\times$ to a bi-semilinear
map $\diamond : \cW\times \cW\to \cW$, and set
$w^\natural = \frac 12 w\diamond w$ and $N(w) = \frac 16 h(w,w\diamond w)$
for $w\in \cW$.
Then:

(a)\ms  $h$ is nondegenerate, $N:\cW \to \cE$ is the unique cubic
form extending $N_\cV$, and $(h,N)$ satisfies the adjoint identity with adjoint
$\natural$.

(b)\ms
$\CD(\cJ,\theta,\mu)=\cA(h,N)$ as algebras with  involution.

(c)\ms $\cZ(\CD(\cJ,\theta,\mu)) = K$ and $\CD(\cJ,\theta,\mu)$ is a structurable algebra.
\end{lemma}

\begin{proof}  For (a) it is clear that
$h$ is nondegenerate and that $N$ is a cubic form extending $N_\cV$.
Also, if $1\le n \le 4$, any degree $n$ homogeneous polynomial map
from $\cW$ to $\cE$ that vanishes on $\cV$ is $0$.  This tells
us that the extension of $N_\cV$ is unique.  Furthermore,
since $\partial_{x}N|_{w}=h(x,w^{\natural})$ and $(w^\natural)^\natural = N(w)w$
hold for $w\in \cV$ they hold for $w\in \cW$.
For (b), observe first that
as $K$-modules
\[
\CD(\cJ,\theta,\mu) = \cJ\oplus s_0\cJ=K\oplus\cV\oplus s_0K\oplus
s_0\cV=\cE\oplus\cW = \cA(h,N).
\]
It is also clear that the involutions on $\CD(\cJ,\theta,\mu)$ and $\cA(h,N)$ are identical.
Furthermore, using the multiplication in $\CD(\cJ,\theta,\mu)$, we have $s_0%
(s_0v)=\mu v=s_0^{2}v$, so $\cV\oplus s_0\cV=\cW$ as left
$\cE$-modules. Since $*$ is an involution, we also have
$wb=(b^{*}w)^{*}=b^{*}w$ in $\CD(\cJ,\theta,\mu)$ for
$b\in\cE$ and $w\in\cW$. Also, if $u,v,x,y\in\cV$, we can
compute in $\CD(\cJ,\theta,\mu)$:
\begin{align*}
(u+s_0v)(x+s_0y)  & =ux-\mu(vy)^{\theta}-s_0(uy)+s_0(vx)^{\theta}\\
& =T(u,x)-\mu T(v,y)-s_0T(u,y)+s_0T(v,x)\\
& \qquad +u\times x+\mu(v\times y)-s_0(u\times y)-s_0(v\times x)\\
& =h(u+s_0v,x+s_0y)+(u+s_0v)\diamond(x+s_0y).
\end{align*}
So $\CD(\cJ,\theta,\mu) = \cA(h,N)$ as algebras with involution.
Finally, (c) follows from Lemma \ref{lem:centre} and Theorem~\ref{thm:cubic}.
\end{proof}

Suppose next that $K$ is a commutative associative algebra,
$\cJ$ is a  Jordan algebra over $K$, and  $t: \cJ \to K$ is a $K$-linear
form.  Following \cite[\S~0.10]{BZ}, we say that $t$ is a \emph{trace form} on $\cJ$ if
$t(x(yz)) = t((xy)z)$ for $x,y,z\in \cJ$. Also, we say that $t$ is  \emph{nondegenerate}
if the $K$-bilinear form $(x,y) \mapsto t(xy)$ is nondegenerate.
If $t(1) = 4$, we say that
$t$ satisfies the \emph{Cayley-Hamilton trace identity of degree 4} if
\begin{equation*}
\ch_4(x) := x^4 + q_3(x)x^3 + q_2(x) x^2 + q_1(x) x + q_0(x)1 = 0
\end{equation*}
for $x\in \cJ$, where the coefficients $q_i(x)\in K$ are defined by
\[
\begin{gathered}
q_3(x) = -t(x),\
q_2(x) = \frac 12 (t(x)^2-t(x^2)),\
q_1(x) = \frac 16 (3t(x)t(x^2) - 2t(x^3) - t(x)^3),\\
q_0(x) = \frac 1{24}
(3t(x^2)^2 + 8 t(x)t(x^3) - 6t(x^4) - 6 t(x)^2t(x^2) + t(x)^4).
\end{gathered}
\]
(These expressions are of course familiar from Newton's identities for the coefficients
of the characteristic polynomial of a $4\times 4$-matrix.)

An important classical example of the preceding
occurs when
$K$ is a field, $\cJ$ is a finite dimensional separable degree 4 Jordan algebra over $K$
and $t$ is the generic trace on $\cJ$ \cite[\S VI.3 and VI.6]{J2} .
In that case, it is well known that $t(1) = 4$ and $t$ is a nondegenerate trace form
that satisfies $\ch_4(x) = 0$ (see for example \cite[Prop.~5.1]{AF}).

\begin{corollary}
\label{cor:CDstructurable} Suppose $\cJ$ is a  Jordan algebra,
$K$ is a subalgebra of the centre of $\cJ$,
and  $t: \cJ \to K$ is a nondegenerate trace form such that $t(1) = 4$.
Let $\cV = \{v\in\cJ: t(v)=0\}$, and define a symmetric $K$-bilinear form $T : \cV \times \cV \to K$ and
a cubic form $N_\cV : \cV \to K$~by
\[T(u,v)   =\frac{1}{4}t(uv) \andd N_\cV(v)   =\frac{1}{24}t(v^{3}).\]
Then:

(a)\ms $T$ is nondegenerate, $(T,N_\cV)$ has an adjoint, and
$\cJ = \cA(T,N_\cV)$.

(b)\ms $(T,N_\cV)$ satisfies the adjoint identity if and only if
$t$ satisfies the Cayley-Hamilton identity of degree 4.

(c)\ms If $t$ satisfies the Cayley-Hamilton identity of degree 4, if $\mu\in K^\times$, and
if we define $\theta : \cJ \to \cJ$ by $\theta = \frac{1}{2}t-\id$, then
$\CD(\cJ,\theta,\mu)$ is a structurable algebra with centre~$K$.
\end{corollary}

\begin{proof}  (a): Note that we  have $\cJ = K \oplus \cV$ as $K$-modules.
Hence, since $t$ is nondegenerate,
$T$ is nondegenerate.  Also, since $t$ is a trace form, we have
$\partial_{v}N_\cV|_{u}=\frac{1}{8}t(vu^{2})$ for $u,v\in \cV$.
Thus, $\partial_vN_\cV|_u = T(v,\frac 12 u^2) = T(v,\frac 12 u^2 - \frac 18 t(u^2)),$
so $(T,N_\cV)$ has an adjoint given by
\begin{equation}
\label{eq:usharp}
v^{\#}=\frac 12 (v^{2}-\frac{1}{4}t(v^{2})).
\end{equation}
Next, it follows from
\eqref{eq:usharp} that $uv=T(u,v)+u\times v$, so
$\cJ = \cA(T,N_\cV)$ as algebras (with trivial involution).

(b): Now
$q_i : \cJ \to K$ is a homogeneous polynomial map of degree $4-i$ over $K$
for $0\le i \le 3$, and it is easy to check that
\begin{equation}
\label{eq:diffq}
\partial_{1}q_{i}|_{x} = -(i+1)q_{i+1}(x)
\end{equation}
for $x\in \cJ$ and $1\le i \le 3$, where $q_4(x) := 1$.  Also,
$\ch_4: \cJ \to \cJ$ is a homogeneous polynomial map of degree $4$ over $K$,
and it follows from the definition of $\ch_4$ and \eqref{eq:diffq} that $\partial_{1}\ch_4|_{x} = 0$.
Let $f$ be the full linearization of $\ch_4$.  Then,
$\partial_{1}\ch_4|_{x} = \frac 16 f(1,x,x,x)$, so $f(1,x,x,x) = 0$ for $x\in \cJ$.
Linearization gives $f(1,x_1,x_2,x_3) = 0$ for $x_1,x_2,x_3\in\cJ$.
So, if $a\in K$ and $v\in \cV$, we have
\[\ch_4(a+v) = \frac 1{24} f(a+v,a+v,a+v,a+v) =  \frac 1{24} f(v,v,v,v) = \ch_4(v).\]
Another easy calculation shows that
$8((v^\#)^\# - N_\cV(v) v) = \ch_4(v)$ for $v\in \cV$.
Thus, $\ch_4(a+ v) = 8((v^\#)^\# - N_\cV(v) v)$ for $a\in K$ and  $v\in \cV$.

(c): We have $\al^{\theta}=\al$ and $v^{\theta}=-v$ for $\al\in K$,
$v\in\cV$.  Thus (c) follows from Lemma \ref{lem:CD}(c).
\end{proof}

Corollary \ref{cor:CDstructurable}(c) generalizes \cite[Theorem
6.6]{AF} which dealt with the case when $K$ is a field and $\cJ$ is
finite dimensional over $K$. Corollary \ref{cor:CDstructurable}(c)
can be thought of as a degree 4 analog of the generalized Tits'
construction of Lie algebras described in \cite[Prop.~24]{BZ}, since
the structurable algebra $\CD(\cJ,\theta,\mu)$ can in turn be used
to construct a $\bbZ$-graded Lie algebra using Kantor's construction
\cite{K2,A2}.

Although we will not use the matrix construction of structurable algebras
(see \cite[p.~257]{K2}, \cite[\S~8(v)]{A1} in the finite dimensional case),
we now briefly consider it as a special case of our
construction using cubic forms.  For this purpose,
we require a pair formulation of the adjoint identity. \label{page:paircubic}
Suppose that $K$ is a commutative and associative $F$-algebra, that
$(\cV_{+},\cV_{-})$ is a pair of $K$-modules, that
$N_{\sg}:\cV_{\sg}\rightarrow K$ is a pair of cubic forms, and that
$T:\cV_{+}\times\cV_{-}\rightarrow K$ is a nondegenerate
$K$-bilinear pairing, i.e., for $v_{\sg}\in\cV_{\sg}$ we have
$T(v_{+},\cV_{-})=0$ implies $v_{+}=0$ and $T(\cV_{+},v_{-})=0
$ implies $v_{-}=0$. To allow symmetry of notation, we set $T(v_{-},v_{+})=T(v_{+},v_{-})$.
We say $(T,N_{+},N_{-})$ has an \textit{adjoint} if
for $\sg = \pm$ and for each $u\in\cV_{\sg}$ there is $u^{\#}\in\cV_{-\sg}$
with $\partial_{v}N_{\sg}|_{u}=T(v,u^{\#})$  for all $v\in\cV_{\sg}$.
If $(T,N_+,N_-)$ has an adjoint,
the \textit{adjoint map} $u \mapsto u^{\#}$  from $\cV_\sg$ to $\cV_{-\sg}$ is unique,
and we define $u\times v = (u\times v)^\# - u^\# - v^\#$
for
$u,v\in \cV_\sg$.
We say $(T,N_{+},N_{-})$ satisfies
the \textit{adjoint identity} if it has an adjoint and
$(v^{\#})^{\#}=N_{\sg}(v)v$  for all $v\in\cV_{\sg}$, $\sg = \pm$.

If $(T,N_{+},N_{-})$ has an adjoint, following \cite[p.~148]{A1},
we form the algebra with involution
\[
\cM(T,N_+,N_-)=
\left\{
\left[
\begin{matrix}
\al_{+} & v_{+}\\
v_{-} & \al_{-}
\end{matrix}
\right]  \suchthat \al_{\sg}\in K,v_{\sg}\in\cV_{\sg}\right\}
\]
with
\[
 \left[
\begin{matrix}
\al_{+}\vphantom{\beta_{+}} & v_{+}\\
v_{-} & \al_{-} \vphantom{\beta_{+}}
\end{matrix}
\right]
\left[
\begin{matrix}
\beta_{+} & u_{+}\\
u_{-} & \beta_{-}%
\end{matrix}
\right]  =
\left[
\begin{matrix}
\al_{+}\beta_{+}+T(v_{+},u_{-}) &\mspace{-15mu}\al_{+}u_{+}+\beta_{-}v_{+}+v_{-}\times u_{-}\\
\beta_{+}v_{-}+\al_{-}u_{-}+v_{+}\times u_{+} & \al_{-}\beta_{-}+T(v_{-},u_{+})
\end{matrix}
\right]
\]
and
\[
\left[
\begin{matrix}
\al_{+} & v_{+}\\
v_{-} & \al_{-}%
\end{matrix}
\right]^{*}
=\left[
\begin{matrix}
\al_{-} & v_{+}\\
v_{-} & \al_{+}%
\end{matrix}
\right] .
\]

\bigskip

Part  (c) of the next corollary generalizes \cite[\S8, Example~(v)]{A1}
which dealt with the case when $K$ is a field and the pair $(\cV_+,\cV_-)$ is finite dimensional.

\begin{corollary}
\label{cor:matrix}
Suppose that $K$ is a commutative associative algebra over $F$ and that
the triple $(T,N_{+},N_{-})$, defined on a pair of $K$-modules $(\cV_+,\cV_-)$, has an adjoint.
Let $\cE = K\oplus K$ with the
exchange involution $*$ and let $\cW=\cV_{+}\oplus\cV_{-}$.  For
$\al=(\al_{+},\al_{-})\in\cE$ and $v=(v_{+},v_{-}%
),u=(u_{+},u_{-})\in\cW$, we set $\al v  =(\al_{+}v_{+},\al_{-}v_{-}),$
\[N(v)   =(N_{+}(v_{+}),N_{-}(v_{-})) \andd
h(v,u)  =(T(v_{+},u_{-}),T(v_{-},u_{+})).\]
Then:

(a)\ms $\cW$ is a left $\cE$-module, $N$ is a
cubic form, $h$ is a nondegenerate hermitian form,
$(h,N)$ has an adjoint, and
$\cM(T,N_+,N_-)\simeq\cA(h,N)$ as algebras with involutions.

(b)\ms If $(T,N_+,N_-)$ satisfies the adjoint identity, then
so does $(h,N)$.

(c)\ms If  $(T,N_+,N_-)$ satisfies the adjoint identity then
$\cM(T,N_+,N_-)$ is a structurable algebra with centre $K1$.
\end{corollary}

\begin{proof} For (a), it is clear that $\cW$ is a left $\cE$-module, $N$ is a
cubic form, and $h$ is a nondegenerate hermitian form.  Also,
$(h,N)$ has an adjoint given by
$(v_+,v_-)^\natural =(v_{-}^{\#},v_{+}^{\#})$.
Identifying $(\al_+,\al_-) + (v_+,v_-)$ in $\cA(h,N)$ with
$\left[
\begin{smallmatrix}\al_+ & v_+ \\ v_- & \al_{-} \end{smallmatrix}
\right]$
in $\cM(T,N_+,N_-)$, we see that $\cM(T,N_+,N_-) = \cA(h,N)$ as algebras with involution.
(b) is clear, and (c) follows from (a), (b),
Theorem~\ref{thm:cubic} and Lemma  \ref{lem:centre}.
\end{proof}

\section{The geometry of class III tori}
\label{sec:geometry}

With the general structurable algebra constructions from \S\ref{sec:Cubicforms}
in hand, we now begin our classification of class III structurable tori.
In this section we associate a finite incidence geometry to
any class III torus.  Our analysis of this geometry will
limit the possibilities for class III tori.

Throughout the section we assume that $\cA$ is a class III structurable torus, and  we use the
notation $S$, $S_\sg$, $\Lm_-$, $\cZ$, $\Gm$, $\cE$ and $\cW$ from \S\ref{sec:classes}.

Recall from \S\ref{sec:classes} that $\cA=\cE\oplus\cW$ with $\cE\cE \subset\cE$ and $\cE\cW+\cW\cE\subset\cW$.
Moreover, we have
\[we = e^*w\]
for $e\in \cE$ and $w\in\cW$ \cite[Prop.~8.3(a)]{AY}.  Also, we recall the following from
\cite[Prop.~8.13]{AY}.

\begin{proposition}
\label{prop:fromA&Y}If $\cA$ is a structurable $\Lm$-torus of
class III, $\sg_{0}\in S_{-}$, and $0\neq s_0\in\cA^{\sg_{0}}%
$, then

(a)\ms $\cE$ is a commutative, associative $\Lm_{-}$-torus
with involution.

(b)\ms $\cZ= \cZ(\cE) =\cE_{+}$.
\ Consequently, if $x$ is a homogeneous element of $\cA$, then
$x^{2}\in\cE\Longleftrightarrow x^{2}\in\cZ$.

(c)\ms $\cE_{-}=s_0\cZ$, $s_0^{2}\in\cZ$, and
$\cE=\cZ\oplus s_0\cZ$.

(d)\ms $S_{-}=\sg_{0}+\Gm$, $2\sg_{0}\in\Gm$, and $\Lm_{-}%
=\Gm\cup(\sg_{0}+\Gm)$.

(e)\ms $(\Lm_{-}:\Gm)=2$.

(f)\ms $4\Lm\subset\Gm\subset\Lm_{-}\subset\Lm$.
\end{proposition}

We now let $\cE$ act on $\cW$ by means of the left multiplication
action $(a,w)\mapsto aw$ in $\cA$.  We
define maps $h : \cW \times \cW \to \cE$ and $\diamond : \cW \times \cW \to \cW$
by the equality
\begin{equation*}
w_{1}w_{2}=h(w_{1},w_{2})+w_{1}\diamond w_{2}%
\end{equation*}
for $w_1,w_2\in \cW$.   We also define a map $N:\cW \to \cE$
by
\begin{equation*}
N(w)=\frac{1}{6}h(w,w\diamond w)
\end{equation*}
for $w\in \cW$.

\begin{proposition} \label{prop:ClassIIIcubic}
Suppose that $\cA$ is a structurable torus of class III.
Then

(a)\ms $\cW$ is a graded left $\cE$-module relative to the action  $(a,w)\mapsto aw$.

(b)\ms $h$ is a nondegenerate graded hermitian form on $\cW$ over $\cE$ and
$\diamond$ is an $\cE$-bi-semilinear symmetric graded map.

(c)\ms  $N$ is a graded cubic form on $\cW$ over $\cE$; that is,
$N$ is cubic form and its full linearization is graded.

(d)\ms  The pair  $(h,N)$ satisfies the adjoint identity with
$w^{\natural}=\frac{1}{2}w\diamond w$
for $w\in \cW$.

(e)\ms $\cA = \cA(h,N)$ as $\Lm$-graded algebras,
where $\cA(h,N)$ has the $\Lm$-grading extending the $\Lm$-grading
on $\cE$ and $\cW$.
\end{proposition}

\begin{proof}
(a): For general structurable tori, $\cW$ is a left $\cE$-module
relative to $a\circ w = a^{*}w$ \cite[Prop.~8.4(a)]{AY} . However,
since $\cA$ is of class III, $\cE$ is commutative, so $\cW$ is an $\cE$-module relative to $aw$.
The action $aw$ is clearly graded.

(b): $h$ and $\diamond$ are clearly graded. Also,
in the notation of \cite[\S 8]{AY}, we have $h(w_1,w_2) = \chi(w_2,w_1)$ and
$w_1\diamond w_2 = \xi(w_{2},w_{1})$.  (b) then follows from the corresponding
statements about $\chi$ and $\xi$ \cite[Prop.~8.4 (b)--(e)]{AY}.

(c):  We have $h(w_1,w_2\diamond w_3) = \overline{\chi(w_1,\xi(w_3,w_2)}$ and the right
hand side is symmetric in its arguments by \cite[p.~128]{AY}.  Hence,
$h(w_1,w_2\diamond w_3)$ is symmetric in its arguments, so $N$ is a cubic
form over $\cE$ with full linearization $h(w_1,w_2\diamond w_3)$.

(d) and (e): Since $\partial_u N\mid_v = \frac 12 h(u,v\diamond v)$, the pair $(h,N)$
has an adjoint $v^\natural = \frac 12 v\diamond v$.  Also, we have
\begin{equation*}
(a+v)(b+w)  =ab+h(v,w)+aw+b^{*}v+v\diamond w \andd (a+v)^{*}   =a^{*}+v
\end{equation*}
for $a,b\in\cE$ and $v,w\in\cW$.  Therefore
$\cA = \cA(h,N)$ as graded algebras with involution. Hence $(h,N)$ satisfies the adjoint
identity by Theorem~\ref{thm:cubic}.
\end{proof}

Let \[\bar{\Lm}=\Lm/\Lm_{-} \andd \bal=\al+\Lm_{-}\in\bar{\Lm},\]
for $\al\in \Lm$.  So
\[4\bal=0\]
for $\al\in\Lm$ by Proposition
\ref{prop:fromA&Y}(f).

If $\beta\in\bal$ for $0\neq
\bal\in\bar{S}$, then $\cA^{\beta}=e\cA^{\al}$
for any $0\neq e\in\cE^{\beta-\al}$. Since
$h(w_{1},w_{2}\diamond w_{3})$ is symmetric and trilinear, we see that for
$0\neq\bal_{i}\in\bar{S}$, the condition
$h(\cA^{\al_{1}},\cA^{\al_{2}}\diamond\cA^{\al_{3}})\neq0$ does not
depend on the order of the $\al_{i}$ or on the choice of representatives
for $\bal_{i}$.

We can now define an incidence geometry $\cI$ associated with $\cA$.
The \emph{points} of $\cI$ are the elements
$\bal\in\bLm$ for $\al\in S\smallsetminus\Lm_{-}$.
The \textit{lines} of $\cI$ are the unordered $3$-tuples $[\bal,\bbe,\bgm]$
of points $\bal,\bbe,\bgm$ with
\[h(\cA^{\al},\cA^{\beta}\diamond\cA^{\gm})\neq
0.\]
Note, in this case, that $\bal+\bbe+\bgm=0$.
We also note that since $\cA$
is of class III, $\cI$ has at least one point and one line.

If
$[\bal,\bbe,\bgm]$ is a line, we say that
$\bal$ is \textit{incident} to $[\bal,\bbe,\bgm]$
and we say that $\bal$ and $\bbe$ are
\textit{collinear}.
Since $h$ is nondegenerate, we see for points
$\bal,\bbe$ that $h(\cA^{\al},\cA^{\beta})\neq0$
if and only if $\bal+\bbe=0$. Thus, for points $\bal,\bbe$,
\[\text{$\bal$ and $\bbe$ are collinear} \iff \cA^{\al}\diamond\cA^{\beta}\neq0.\]

We write $\ord{\bal}$ for the order of $\bal\in\bar{\Lm}$, in which case
$\ord{\bal}= 1$, $2$ or $4$. Moreover, if $\bal$ is a point,
$\ord{\bal}= 2$ or $4$.

It is important to note that our definition of collinear allows the possibility
of a point $\bal$ being collinear with itself.  Indeed, by definition, this holds
if and only if $[\bal,\bal,\bgm]$ is a line for some point $\bgm$.  Moreover, we see in part
(a) of the next lemma that this is the case if and only if $\ord{\bal} = 4$.  We also note that
if $\bal$ is a point, then so is $-\bal$ (by (ST2)), but $\bal$ is clearly never collinear with $-\bal$.

\begin{lemma}
\label{lem:geometry} Suppose that $\bal,\bbe,\bgm,\bde$
are points in the geometry $\cI$ for a class III structurable
$\Lm$-torus $\cA$.  Then:

(a)\ms If $\ord{\bal}=4$, then $[\bal,\bal,2\bal]$ is a line. These are the only lines not having
three distinct points and the only lines with points of both orders, $2$ and
$4$.

(b)\ms If $[\bal,\bbe,\bgm]$ is a line, then $\bar{\delta
}$ is collinear with at least one of $\bal,\bbe,\bgm.$

(c)\ms If $\ord{\bde} = 2$ and $\bde$ is collinear with $\bal$, then $\bde$
is collinear with $-\bal$.

(d)\ms If $[\bal,\bbe,\bgm]$ is a line, if $\bde$ is collinear with
both $\bal$ and $\bbe$, and if $\ord{\bde} =2$, then $\bde=\bgm$.
\end{lemma}

\begin{proof}
We let $0\neq a\in\cA^{\al}$, $0\neq b\in\cA^{\beta}$,
$0\neq c\in\cA^{\gm}$, and $0\neq d\in\cA^{\delta}$.

(a):   If $\ord{\bal}=4$, then $h(\cA^{\al},\cA^{\al})=0$,
so $a\diamond a = a^2 \ne 0$ by Remark \ref{rem:powers}.
Hence, $2\bal=-2\bal$ is a point and $[\bar
{\al},\bal,2\bal]$ is a line. \
Moreover, if $[\bar{\al},\bal,\bgm]$ is a line, then $\bgm=-\bal%
-\bal=2\bal$, so $\ord{\bal}=4$. Also, if
$[\bal,\bbe,\bgm]$ is a line with
$\ord{\bal} =4$ and $\ord{\bgm} =2$, then
$\ord{\bal+\bgm} =4$. Thus, $c^{\natural}=\frac{1}%
{2}c\diamond c=0$ and $(a\diamond c)^{\natural}\neq0$. Using (ADJ4), we see
that $h(c,a^{\natural})c=(a\diamond c)^{\natural}\neq0$ so $\bgm%
+2\bal=0$. Therefore, $\bbe= -\bal - \bgm = -\bal+2\bal =\bal$

(b): Suppose $\bde$ is not collinear with
$\bal,\bbe,\bgm$. Since $d\diamond a=d\diamond
b=d\diamond c=0$, (ADJ3) gives $h(a,b\diamond c)d=0$, a contradiction.

(c): We can assume that $\ord{\bal} = 4$.  Since $\bde$ is collinear with $\bal$,
we have $\bde = 2\bal$ by (a).  Therefore, $\bde = 2(-\bal)$, so $\bde$ is collinear with
$-\bal$ by (a).

(d): We know that $\bde$ is collinear with $-\bal$ and $-\bbe$ (by (c)).
So $[\bde,-\bal,\bal-\bde]$ and
$[\bde,-\bbe,\bbe-\bde]$ are lines. Let $0\neq
u\in\cA^{\al-\delta}$ and $0\neq v\in\cA^{\beta-\delta}$.
Thus, $0\neq(d\diamond u)\diamond(d\diamond v)\in\cA^{\al
}\diamond\cA^{\beta}$.  Since $d^{\natural}=0$, (ADJ2) shows that
$h(d,u\diamond v)d\neq0$.  Thus  $h(d,u\diamond v)\ne 0$, so
$\bde=-(\bal-\bde)-(\bbe-\bde)=-\bal-\bbe=\bgm$.
\end{proof}

If a point $\bal$ is collinear with all points $\bbe\neq\bal$
and $\bal$ is incident to each line, we say that $\cI$ is a
\textit{star} with \textit{centre} $\bal$.

\begin{corollary}
\label{cor:geometry} Let  $\cA$ be a structurable $\Lm$-torus of class III.

(a) \ If $\cI$ is a star, then any centre $\bal$ of\/ $\cI$ has order 2.

(b) \ If $\cI$ is a star, the centre  of\/ $\cI$
is unique unless $\cI$ consists of 3 distinct points $\bal$, $\bbe$, $\bgm$
of order 2 and one line $[\bal,\bbe,\bgm]$.

(c) \ If some point $\bal$ is collinear with all points $\bbe \ne \pm \bal$, then
$\cI$ is a star.

(d) \ $\cI$ is a star if and only if there exists a homogeneous element $x\in \cW$ such that
$L_x$ is invertible.
\end{corollary}

\begin{proof}
(a): Suppose that $\cI$ is a star with centre $\bal$.
Since $\bal$ is not collinear with $-\bal$,
we have $-\bal = \bal$, so
$\ord{\bal} = 2$.

(b):  Suppose that $\bal \ne \bbe$ are centres of $\cI$.
Then there is a line $[\bal,\bbe,\bgm]$ and $\ord{\bal} = \ord{\bbe} = 2$, so
$\ord{\bgm} = \ord{-\bal-\bbe} = 2$.  Let $\bde$ be a point not equal to $\bal$ or $\bbe$.
If $\ord{\bde} = 4$, then Lemma \ref{lem:geometry}(a) tells us that
$\bal  = 2\bde= \bbe$, a contradiction.
So $\ord{\bde} = 2$, and thus $\bde = \bgm$ by Lemma \ref{lem:geometry}(d).
Therefore, $\bal$, $\bbe$ and $\bgm$ are the only points and $[\bal,\bbe,\bgm]$ is the only line.

(c): If $\ord{\bal} =2$, then
$\bal$ is incident to
each line by Lemma  \ref{lem:geometry}(d), so $\cI$ is a star with
centre $\bal$. Suppose then that $\ord{\bal} =4$. If
$\bbe$ is any point of order 2, then $\bbe$ is collinear with
$\bal$, so $\bbe=2\bal$ by Lemma \ref{lem:geometry}(a).
Thus, if $\bbe$ is any point not equal to $2\bal$, we have
$\ord{\bbe} = 4$ and hence $\ord{2\bbe} = 2$.  But then, as we've just seen,
$2\bbe=2\bal$, so $2\bal$ is collinear with $\bbe$  by Lemma \ref{lem:geometry}(a).
Thus, $2\bal$ is collinear with all points $\bbe \neq 2\bal$,
and the first case tells us that $\cI$ is a star with centre $2\bal$.

(d):  Suppose first that $x\in\cW^\al$ with $L_x$ invertible.
Then, $\bal$ is a point, and
$\cA^\al \diamond \cA^\beta = L_x \cA^\beta \ne 0$
for all points $\bbe \ne -\bal$.  So $\cI$ is a star by (c). Conversely,
suppose that $\cI$ is a star with centre $\bal$.  We choose $0\ne x \in \cW^\al$
and show that $L_x$ is invertible.
Since $\cA$ is finely graded, it is sufficient to show that $L_x \cA^\beta \ne 0$
for all $\beta\in S$.  But if $\beta \in \Lm_-$, this is clear from the multiplication
in $\cA = \cA(h,N)$.  So we can assume that $\bbe$ is a point.  If $\bbe \ne \bal$,
then $\bbe$ and $\bal$ are collinear, so $L_x \cA^\beta \ne 0$.  Finally,
if $\bbe = \bal$, then
$\cA^\beta = \cE^{\beta-\al}x$,
so $L_x \cA^\beta  = \cE^{\beta-\al}x^2\ne 0$ by Remark \ref{rem:powers}.
\end{proof}

It is natural and convenient now to subdivide class III structurable
tori into three subclasses:

\begin{definition}
\label{def:trichotomy}
If $\cA$ is a class III structurable $\Lm$-torus, we
say that $\cA$ has \textit{class} III(a), III(b), or III(c), if
the corresponding condition below holds:
\begin{mylist}
\item[III(a):] $\cI$ is a star.
\item[III(b):] $\cI$ is not a star and all points have order $2$.
\item[III(c):] $\cI$ is not a star and there is a point of order $4$.
\end{mylist}
\end{definition}

In view of Lemma \ref{lem:geometry}(a), condition III(b) can be equivalently stated as:
$\cI$ is not a star and each line of $\cI$ has 3 distinct points. Similarly,
condition III(c) can be stated as:
$\cI$ is not a star and there is a line of $\cI$ with a repeated point.  So,
the trichotomy in Definition \ref{def:trichotomy} is purely based on the properties of the incidence
 geometry~$\cI$.

\begin{example}
To  provide some intuition, we now give examples of
incidence geometries that are associated with structurable tori of class III(a), III(b) and III(c).
In the diagrams, solid circles represent points of order 2 and open circles represent
points of order 4.

\bigskip
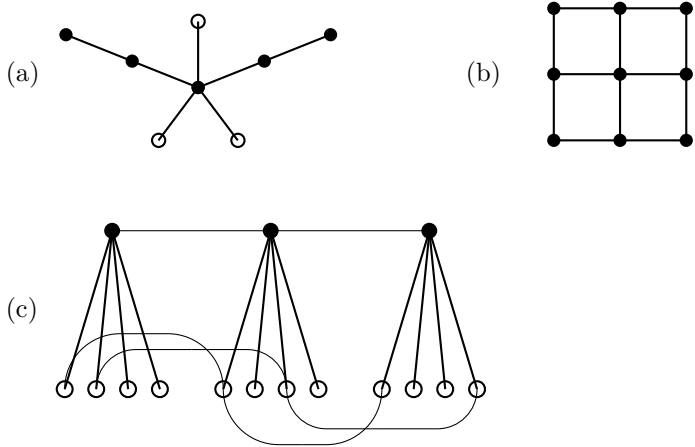
\newcommand\ptwo{\circle*{1}} 
\newcommand\pfour{\circle{1}} 
\newcommand\shower 
{\begin{picture}(6,10)\thicklines
    \put(0,0){\Line(3,10)}
    \put(2,0){\Line(1,10)}
    \put(4,0){\Line(-1,10)}
    \put(6,0){\Line(-3,10)}
    \put(0,0){\pfour}
    \put(2,0){\pfour}
    \put(4,0){\pfour}
    \put(6,0){\pfour}
    \put(3,10){\ptwo}
    \end{picture}}

(a)\quad
\parbox{1.5truein}
{
\setlength{\unitlength}{5pt}
\begin{picture}(20,10)(-10,-4)\thicklines
\put(0,0){\Line(10,4)}
\put(0,0){\Line(-10,4)}
\put(0,0){\Line(0,5)}
\put(0,0){\Line(-3,-4)}
\put(0,0){\Line(3,-4)}
\multiput(0,0)(5,2){3}{\ptwo}
\multiput(0,0)(-5,2){3}{\ptwo}
\put(0,5){\pfour}
\put(-3,-4){\pfour}
\put(3,-4){\pfour}
\end{picture}
}
\hspace{.5truein}
(b)\qquad
\parbox{1truein}
{
\setlength{\unitlength}{5pt}
\begin{picture}(10,10)\thicklines
\put(0,0){\Line(0,10)}
\put(5,0){\Line(0,10)}
\put(10,0){\Line(0,10)}
\put(0,0){\Line(10,0)}
\put(0,5){\Line(10,0)}
\put(0,10){\Line(10,0)}
\multiput(0,0)(0,5){3}{\ptwo}
\multiput(5,0)(0,5){3}{\ptwo}
\multiput(10,0)(0,5){3}{\ptwo}
\end{picture}
}

\bigskip\bigskip
(c)\quad
\parbox{2.5truein}
{
\setlength{\unitlength}{6pt}
\begin{picture}(26,10)
\put(0,0){\shower}
\put(10,0){\shower}
\put(20,0){\shower}
\put(3,10){\Line(20,0)}
\put(8,0){\oval(12,5)[t]}
\put(20,0){\oval(12,5)[b]}
\put(5,0){\oval(10,7)[t]}
\put(15,0){\oval(10,7)[b]}
\end{picture}
}

\bigskip\bigskip\bigskip\noindent
In example (c), we have  drawn only 2 of the 16 lines that are
made up of points from the 3 clusters at the bottom of the picture.
Indeed, if we label the points in each of these clusters by the elements
of the group $\bbZ_2\oplus\bbZ_2$, there are 16 lines of the form $[\bal_1,\bal_2,\bal_3]$,
where the points $\bal_1$, $\bal_2$, $\bal_3$ are chosen one from each cluster with
labels summing to~0.  (See Lemma \ref{lem:IIIcgeometry} below in the case when
all  $k_i = 1$.)
\end{example}

\begin{remark}
If $\cA$ is of class III(b), then, using Lemma \ref{lem:geometry},
$\cI$ can be seen to be
a generalized quadrangle $\cQ$ of order $(2,t)$ (see \cite[\S 1.1]{PT})
imbedded in the
projective space $\cP$ of the $\bbZ_2$-vector space
$\Lm/\Lm_-$.  Furthermore, one can classify such pairs $(\cQ,\cP)$.
Methods similar to those in \cite{F1} should then classify the pairs
$(h,N)$ leading to class III(b) tori. Although we have been
guided in our research by these facts, in the end
we adopted a different approach to the
classification of class III(b) tori (see \S\ref{sec:IIIb}).
\end{remark}

\section{Tori of class III$(a)$}
\label{sec:IIIa}

In this section we obtain the classification of structurable tori
of class III(a).

Recall that in Example \ref{ex:herm}, we constructed
the structurable torus $\cA(\hermk)$ associated with a
diagonal graded hermitian form $\hermk$ over
an associative torus with involution $\cB$.
We saw in Theorem \ref{thm:classII}
that if we take $\cB$ of class I in this construction
we obtain all class II structurable tori.
On the other hand, there are,
up to isograded-isomorphism,
precisely two associative tori with involution that are not
of class~I:
$\cP(r)$ and $\Cay_*(2) \ot \cP(r)$
(see Remark \ref{rem:assoctori}).
If we use the first of these for $\cB$ in the construction, we obtain
Jordan tori \cite[\S~5]{Y1}.
If we use the second, we now see that we
obtain all class III(a) structurable tori.

\begin{theorem}
\label{thm:classIIIa}
If  $M$ is a subgroup of
$\Lm$ such that $2\Lm \subset M$,
if $\cB$ is an associative $M$-torus with involution that is isograded-isomorphic
to
$\Cay_*(2) \ot \cP(r)$
for some $r\ge 0$,
and if $\hermk$ is a diagonal $\Lm$-graded hermitian form
over $\cB$, then the structurable $\Lm$-torus $\cA(\hermk)$ associated
with $\hermk$ (as in Example \ref{ex:herm}) is of class III(a).
Moreover, any structurable $\Lm$-torus of class III(a) is graded-isomorphic
to a structurable $\Lm$-torus $\cA(\hermk)$  obtained in this way.
\end{theorem}

\begin{proof} For the first statement,
assume that $\cA = \cA(\hermk) = \cB \oplus \cX$,
where $M, \cB$, $\cX$ and $\hermk$ are as in  Example \ref{ex:herm}, and
$\cB$ is isograded isomorphic to $\Cay_*(2) \ot \cP(r)$.
We use the notation
$S$, $S_\sg$, $\Lm_-$, $\cZ$, $\Gm$, $\cE$ and $\cW$ from \S\ref{sec:classes}
for $\cA$.  Now the centre $\cZ$ of the algebra with involution
$\cA$ is also the centre of  the algebra with involution
$\cB$, so we have $\Gm = \Gm(\cB)$. Thus, since
$\cB$ is isograded-isomorphic to $\Cay_*(2) \ot \cP(r)$,
we may choose a basis $\set{\mu_1,\mu_2,\dots,\mu_{r+2}}$ of $M$ and nonzero elements $a_1\in \cB^{\mu_1}$,
$a_2\in \cB^{\mu_2}$ such that
$\Gm = \langle 2\mu_1,2\mu_2,\mu_3,\dots,\mu_{r+2}\rangle,$
\[\cB = \cZ \oplus \cZ a_1 \oplus \cZ a_2 \oplus \cZ a_1a_2,\]
$a_i^2\in \cZ$, $a_2a_1 = -a_1a_2$, $a_1^* = a_1$ and $a_2^* = a_2$.
Then $S_- = \Gm + \mu_1+\mu_2$, so
$\Lm_- = \langle \Gm, \mu_1+\mu_2\rangle$,
\[\cE = \cZ \oplus \cZ a_1a_2 \andd \cW = \cZ a_1\oplus \cZ a_2 \oplus \cX.\]
Now $\cW \cW \supset \cX \cX = \hermk(\cX,\cX) = \cB$ and hence $\cW \cW \not \subset \cE$.  Thus,
$\cA$ is of class III.
Moreover, $L_{a_1}$ is invertible, so, by Corollary \ref{cor:geometry}(d),
$\cI$ is a star. (This can also  easily be seen directly.)  So, $\cA$ is of class III(a).

For the second statement, suppose
that $\cA$ is a structurable torus of class III(a), and we use
the notation from \S\ref{sec:classes} and \S\ref{sec:geometry} for $\cA$.
We use the recognition theorem \cite[Prop.~9.8]{AY}.
Indeed, according to that result, it is sufficient
to show that  there exists a subgroup $M$ of $\Lm$
satisfying the following conditions:
\begin{itemize}
\item[(a)] $\Lm_- \subset M \subset S$ and $(M:\Lm_-) = 2$
\item[(b)] $\cA^{S\setminus M}\cA^{S\setminus M} \subset \cA^M$
\item[(c)] For each $\sg\in S\setminus M$ there exists $\tau\in S\setminus M$ so
that $\cA^\sg \cA^\tau \not\subset \cE$.
\end{itemize}
To do this, let $\bde$ be a centre of the geometry $\cI$ associated with $\cA$.
Then,
by Corollary \ref{cor:geometry}(a), $\ord{\bde} = 2$.
Let $M=\left\langle \Lm_{-},\delta\right\rangle = \Lm_- \cup (\delta + \Lm_-)$,  in which
case we have condition~(a).  To show (b), it suffices to check that
$\cA^{\al}\cA^{\beta} \subset \cA^M$ for points $\bal,\bbe$ distinct from~$\bde$.
If $\bal+\bbe=0$,
then $\cA^{\al}\cA^{\beta}=h(\cA^{\al},\cA^{\beta})\subset\cE$. \ If $\bal+\bbe%
\neq0$ and $\bal$ and $\bbe$ are not collinear, then
$\cA^{\al}\cA^{\beta}=\cA^{\al}\diamond
\cA^{\beta}=0$. \ If $[\bal,\bbe,\bgm]$ is a
line, then $\bde=\bgm$ (since $\bde$ is a centre of $\cI$)
and $\cA^{\al}\cA^{\beta}=\cA^{\al}\diamond\cA^{\beta}\subset\cA^{M}$
(since $\bar{\al}+\bbe=-\bde$).  Thus, we have (b).
Finally, if $\sg\in S\smallsetminus M$, then $[\bar{\sg},\bde,\bar{\tau
}]$ is a line for $\tau=-\sg-\delta\notin M$, so
$\cA^{\sg}\cA^{\tau}=\cA^{\sg}\diamond\cA^{\tau}=\cA^{-\delta}\not \subset \cE$.
\end{proof}

\section{Tori of class III$(b)$}
\label{sec:IIIb}

In this section we show that each structurable torus of class III(b) can be constructed
by doubling a Jordan torus of degree 4 using the
Cayley-Dickson process (see \S\ref{sec:Cubicforms}).

We will need the following simple fact from Jordan theory:

\begin{lemma}
\label{lem:Jordfact} Suppose that $\cJ$ is a finite
dimensional central simple Jordan algebra  over $F$.  If $x$
is an element of $\cJ$ such that $x^2\in F$ and $L_x$ is invertible,
then $x\in F$.
\end{lemma}

\begin{proof}  By extending the base field and replacing $x$ by a scalar multiple of $x$, we can
assume that $x^2 = 1$.  Then $e := \frac 12(1-x)$ is an idempotent
and we have the Peirce decomposition $\cJ = \cJ_0 \oplus \cJ_1
\oplus \cJ_{\frac 12}$ of $\cJ$ relative to $e$ \cite[\S III.1]{J2}
Since $x = 1-2e$, we have $L_x \cJ_{\frac 12} = 0$, so $\cJ_{\frac
12} = 0$ and $\cJ = \cJ_0 \oplus \cJ_1$ where $\cJ_0$ and $\cJ_1$
are ideals of $\cJ$.  Thus, by simplicity, $e=0$ or $e=1$, so $x =
\pm 1$.
\end{proof}

Suppose that $\cJ$ is a Jordan $\Lm$-torus with centre $\cZ$ and
central grading group $\Gm$. Then $\cJ$ is a free $\cZ$-module
\cite[Lemma 3.9(ii)]{Y1}. Indeed the rank (possibly infinite) of
this module is easily seen to be the cardinality $\ord{S(\cJ)/\Gm}$
of the set $S(\cJ)/\Gm$, where $S(\cJ)$ is the support of $\cJ$ and
\[S(\cJ)/\Gm = \set{\al + \Gm \suchthat \al \in S(\cJ)}.\]
Also, $\cZ$ is an integral domain
\cite[Lemma 3.6(i)]{Y1}, so we can form the field $\tcZ$ of fractions of
$\cZ$.  The Jordan algebra
\[\tcJ =\tcZ\ot_\cZ\cJ\]
over $\tcZ$ is called the \textit{central closure} of $\cJ$. Since
$\cJ$ is a free $\cZ$-module, it follows that the $\cZ$-algebra
$\cJ$ embeds  in $\tcJ$ as $1\ot_\cZ \cJ$. Moreover, each element of
$\tcJ$ can be expressed in the form $z^{-1} \ot_\cZ x$, where
$0\ne z\in \cZ$ and $x\in \cJ$.  Again since $\cJ$ is a free
$\cZ$-module, $\cJ$ is finitely generated as $\cZ$-module if and
only if $\tcJ$ is finite dimensional over $\tcZ$. In fact if $\cJ$
is finitely generated as a $\cZ$-module, its central closure $\tcJ$
is a  finite dimensional central division algebra over $\tcZ$
\cite[3.6, 2.6 and 2.10]{Y1}, and thus the degree of $\tcJ$ over $\tcZ$
is defined as the degree of the generic minimum polynomial of each
element of $\tcJ$ \cite[\S VI.3]{J2}. Following
\cite[Definition~6.4]{Y1}, we say that a Jordan $\Lm$-torus $\cJ$
has \textit{degree} $n$ (or central degree $n$) if $\cJ$ is finitely
generated as a $\cZ$-module and the degree of $\tcJ$ (in the above
sense) over $\tcZ$ is~$n$.

\begin{proposition}\label{prop:CDtorus}
(a)\ms
Let $M$ be a subgroup of $\Lm$ such that $(\Lm:M)=2$; let $\cJ$ be a
degree $4$ Jordan $M$-torus with centre $\cZ$ and central grading
group $\Gm$ satisfying $2M\subset \Gm$; let $\sg_0$ be an element of
$\Lm$ such that
\begin{equation*}
\Lm = \langle M,\sg_0\rangle \andd \quad 2\sg_{0}\in\Gm;
\end{equation*}
and let $0\neq\mu\in\cZ^{2\sg_{0}}$.
Define $\theta : \cJ \to \cJ$ by
\begin{equation}
\label{eq:deftheta}
 \theta(x) = \left\{
\begin{array}{ll}
    \hskip-.5em\hphantom{-}x, & \text{if } x\in\cJ^\Gm \\
    \hskip-.5em -x, & \text{if } x\in\cJ^{M\setminus\Gm}\\
\end{array},
\right.
\end{equation}
and give the algebra with involution $\CD(\cJ,\theta,\mu)$ the $\Lm$-grading defined~by
\begin{equation}
\label{eq:CDgrading} \CD(\cJ,\theta,\mu)^{\tau}=\cJ^{\tau} \andd
\CD(\cJ,\theta,\mu)^{\sg_{0}+\tau}=s_0\cJ^{\tau}.
\end{equation}
for $\tau\in M$. Then, $\CD(\cJ,\theta,\mu)$ is a class III(b)
structurable $\Lm$-torus with centre $\cZ$.

(b)\ms Any class III(b) structurable $\Lm$-torus is graded
isomorphic to a structurable torus $\CD(\cJ,\theta,\mu)$ constructed
from some $M$, $\cJ$, $\sg_0$ and $\mu$ as in (a).
\end{proposition}

\begin{proof}
(a):
Let  $\tcZ$ be the field of fractions of  $\cZ$, let
$\tcJ =\tcZ\ot_\cZ\cJ$, and let
$\ttrace : \tcJ \to \tcZ$ be the generic trace on $\tcJ$.

We first claim that the restriction of $\ttrace$ to
$\cJ^{M\setminus \Gm}$ is zero.%
\footnote{It is not difficult to give a
direct argument for this without using
\cite[Prop.~4.9]{NY}, a result that uses the classification of Jordan tori.}
Indeed, by \cite[Prop.~4.9]{NY},
we have $\cJ = \cZ \oplus (\cJ,\cJ,\cJ)$, where
$(\cJ,\cJ,\cJ)$ is the $F$-span of the associators $(x,y,z)$, $x,y,z\in\cJ$.
Thus, since $(\cJ,\cJ,\cJ)$ is a graded subspace of $\cJ$, it follows
that $\cJ^{M\setminus \Gm} = (\cJ,\cJ,\cJ)$.  Then the claim follows from the
fact that the generic trace $\ttrace$ is a trace form on $\tcJ$.

Now, by assumption, $\tcJ$ is a finite dimensional central division algebra of degree $4$
over $\tcJ$ and hence it is a finite dimensional separable
algebra of degree 4 over $\tcJ$.  So, by Corollary \ref{cor:CDstructurable}(c) (and
the paragraph preceding Corollary \ref{cor:CDstructurable}),
$\CD( \tcJ,\tilde\theta,\mu)$ is a
structurable algebra, where $\tilde\theta : \tcJ \to \tcJ$
is defined by $\tilde\theta=\frac{1}{2}\ttrace-\id$.
But since $\ttrace$ is $\tcZ$-linear, $\ttrace(1) = 4$,
and $\ttrace$ is zero on $\cJ^{M\setminus \Gm}$, it follows that
$\tilde\theta$ restricts to $\theta$  on $\cJ$.
Hence, setting  $\cA=\CD(\cJ,\theta,\mu)$, we see that
$\cA$ is a subalgebra with involution of $\CD(\tcJ,\tilde\theta,\mu)$,
so  $\cA$ is structurable.

It is clear that $\cA$ is a finely $\Lm$-graded structurable algebra
using the grading given by \eqref{eq:CDgrading}. Moreover,
$\supp(\cA)=\supp(\cJ)\cup(\sg_{0}+\supp(\cJ))$ generates $\Lm$. If
$0\neq x\in\cJ^{\tau}$, where $\tau\in M$, then $x$ has a Jordan
inverse $y\in\cJ^{-\tau}$, in which case $xy=yx=1$ and
$[L_{x},L_{y}]\cJ=0$. Then, since $\theta(x)=-x$ if and only if
$\theta(y)=-y$, we have
\[
[ L_{x},L_{y}](s_0\cJ)=
s_0([L_{\theta(x)},L_{\theta(y)}]\cJ)=
s_0([L_{x},L_{y}]\cJ)=0,
\]
so $y$ is the inverse of $x$ in $\cA$. Also, since $L_{s_0}$ is invertible,
it follows using Remark~\ref{rem:invertible}(c) that $s_0x$ is invertible.
Therefore $\cA$ is a structurable torus.  Also, by Corollary \ref{cor:CDstructurable}(c),
$\cZ(\cA) = \cZ$ and hence $\Gm(\cA) =~\Gm$.

Next, using the notation from \S\ref{sec:classes} and \S\ref{sec:geometry} for $\cA$, we have
$S_- = \sg_0 + \Gm$, so $\Lm_- = \langle \Gm, \sg_0 \rangle$,
$\cE = \cZ \oplus s_0\cZ$ and $\cW = \cJ^{M\setminus \Gm} \oplus s_0 \cJ^{M\setminus \Gm}$.
But, since $\tcJ$ is of degree 4, we  have
$\cJ^{M\setminus \Gm}  \cJ^{M\setminus \Gm} \not\subset \cZ$, so
$\cW \cW \not\subset \cE$ and $\cA$ is of
class~III.

It remains to show that $\cA$ is of class III(b).
Since $2M\subset \Gm$, we have $2\Lm \subset \Lm_-$, so
every point in the geometry $\cI$ associated with~$\cA$ has order 2.
Suppose finally for contradiction that $\cI$ is a star.  So, by
Corollary \ref{cor:geometry}(d), $L_x: \cA \to \cA$ is invertible
for some homogeneous $x\in \cW$.  It follows then from the multiplication and grading
in $\CD(\cJ,\theta,\mu)$ that $L_x: \cJ \to \cJ$ is invertible for some homogeneous
$x\in \cJ^{M\setminus \Gm}$.  Moreover, since $2M\subset \Gm$, we have $x^2 \in \cZ$.
This contradicts Lemma \ref{lem:Jordfact} applied to the finite dimensional central simple
algebra $\tcJ$ over $\tcZ$.

(b): Suppose that $\cA$ is a structurable $\Lm$-torus of class III(b).
We use the notation of \S\ref{sec:classes} and \S\ref{sec:geometry} for $\cA$.
As in Proposition \ref{prop:fromA&Y}, we fix a choice of $\sg_0\in S_-$, in which case
$\sg_0\notin \Gm$ and $2\sg_0\in \Gm$.  We also fix a choice of $0\ne s_0\in \cA^{\sg_0}$ and we let
\[\mu = s_0^2 \in \cZ^{2\sg_0}.\]

Now, since every point of $\cI$ has order $2$, we have $2(S\setminus \Lm_-) \subset \Lm_-$.
So $2S \subset \Lm_-$.  But $2S \subset S_+$ \cite[Prop.~7.1]{AY}, so
$2S \subset \Lm_-\cap S_+ = \Gm$.  Therefore $2\Lm \subset \Gm$.

So $\Lm/\Gm$ is a $\bbZ_2$-vector space and $\sg_0\notin \Gm$,
and therefore we can find a subgroup $M$ of $\Lm$ such
that
\[\Gm \subset M,\quad (\Lm:M)=2 \andd \sg_0\notin M.\]
We fix a choice of $M$ with these properties.
Then certainly $\Lm = \langle M, \sg_0\rangle$.

Let
\[\cJ=\cA^{M}.\]
Since $M\cap S_{-}=\emptyset$, we see that
$v^{*}=v$ for all $v\in\cJ$.  Thus $\cJ$ is a structurable algebra with trivial involution,
so $\cJ$ is a Jordan algebra \cite[\S 1]{A1}.
Since $\Lm=M\cup(\sg_{0}+M)$ and $\sg_{0}+S =S$, we have $S=(M\cap S)\cup(\sg_{0}+(M\cap S))$. Since $\Lm
=\left\langle S\right\rangle $, $\sg_0\notin M$ and $2\sg_{0}\in\Gm\subset M\cap S$, we
see that $M=\left\langle M\cap S\right\rangle $.  Therefore $\cJ$ is a Jordan
$M$-torus.

Now
\[
\cA=\cA^{M}\oplus\cA^{\sg_{0}+M}=\cJ\oplus
s_0\cJ \andd \cJ=\cZ\oplus\cV,
\]
where
\[\cV=\cJ^{M\smallsetminus\Gm}.\]
Thus, $\cW=\cV\oplus
s_0\cV\simeq\cE\ot\cV$ as $\cE$-modules.
Let $T$ and $N_{\cV}$ be the restrictions of
$h$ and $N$ respectively to $\cV$ (where $h$
and $N$ are defined as in \S\ref{sec:geometry}).
Then $T$ and $N_\cV$ take values in
$\cE^{M}=\cZ$.  Now, by Proposition \ref{prop:ClassIIIcubic} (d) and (e),
$(h,N)$ satisfies the adjoint identity and
$\cA = \cA(h,N)$.  Thus,
$(T, N_{\cV})$ satisfies the adjoint identity and
$\cJ=\cA(T, N_{\cV})$.
So, by Lemma \ref{lem:CD}~(a) and (b) (with $K = \cE$), we have
\[\cA = \CD(\cJ,\theta,\mu)\]
as algebras with involution, where
$\theta: \cJ \to \cJ$ is defined by
\eqref{eq:deftheta}.  Also, if we give $\CD(\cJ,\theta,\mu)$ a  $\Lm$-grading
by \eqref{eq:CDgrading}, this is an equality of $\Lm$-graded algebras with involution.

Next, the centre $\cZ$ of $\cA$ is clearly  contained in $\cZ(\cJ)$. On the other hand,
if $x$ is a nonzero homogeneous element of $\cZ(\cJ)$, then
$L_x : \cJ \to \cJ$ is invertible, so $L_x : \cA \to \cA$ is invertible.
Since $\cI$ is not a star, this implies that $x\in \cZ$ by Corollary~\ref{cor:geometry}(d).  Thus,
$\cZ(\cJ) = \cZ$ and therefore $\Gm(\cJ) = \Gm$.  Also, since $2\Lm \subset \Gm$, we have
$2M \subset \Gm$.

Now $\cJ$ is a free $\cZ$-module of rank $\ord{S(\cJ)/\Gm}$, and
$S(\cJ)/\Gm \subset M/\Gm$ is finite. Hence, the degree $n$ of the
$M$-torus $\cJ$ is defined. But $(h,N)$ satisfies the adjoint
identity by Proposition 7.2, so $(T,N_\cV)$ satisfies the adjoint
identity $(v^\#)^\# = N_\cV(v) v$ for $v\in \cV$, where $\#$ is the
restriction of $\natural$ to $\cV$. Hence, since $v^\# = \frac 12
(v^2 - T(v,v))$, we see that each $v\in \cV$ is a root of a monic
polynomial of degree $4$ over $\cZ$.  Thus, the same is true for
each element $x\in \cJ$.  So each element of $\tcJ$ is a root of a
monic polynomial of degree $4$ over $\tcZ$. Consequently, $n\le 4$
\cite[p.~224]{J2} (since $\tcZ$ is infinite). Now if $n = 3$, then,
by \cite[Prop.~6.7]{Y1}, $3M \subset \Gm$, which contradicts
$2M\subset \Gm$ and $M \ne \Gm$. Also, if $\cJ$ has degree 2, then
$\tcJ$ is the Jordan algebra of a symmetric bilinear form over
$\tcZ$ \cite[p.~207]{J2}.  Thus, $\cJ^{M\setminus \Gm}
\cJ^{M\setminus \Gm} \subset \cZ$ \cite[p.~153]{Y1}, so $\cW \cW
\subset \cE$, contradicting the fact that $\cA$ is of class~III.  Of
course, $n\ne 1$ since $\cJ \ne F$. So $\cJ$ is a Jordan $M$-torus
of degree~$4$.
\end{proof}

We call the structurable $\Lm$-torus $\CD(\cJ,\theta,\mu)$ constructed in
Theorem \ref{prop:CDtorus}(a)
a \emph{Cayley-Dickson torus obtained from $\cJ$}.

\begin{remark} It is interesting to note that in  the proof of Proposition \ref{prop:CDtorus}(a)
we did not use the assumption $2\Lm \subset M$ to show that
$\CD(\cJ,\theta,\mu)$ is a class III structurable torus.
That assumption was only used to show that $\CD(\cJ,\theta,\mu)$ is of class III(b).
In fact, if we drop the assumption that $2\Lm \subset M$,
we do obtain an extra torus that is of class III(c), but it is not possible
to obtain all class III(c) tori in this way.
\end{remark}

To explain the dependence of the Cayley-Dickson torus
$\CD(\cJ,\theta,\mu)$ on the scalar $\mu$, we prove the following:

\begin{lemma}
\label{lem:mu}
Suppose that $\CD(\cJ,\theta,\mu)$ and $\CD(\cJ',\theta',\mu')$
are Cayley-Dickson tori constructed as in Theorem
\ref{prop:CDtorus}(a) (with the assumptions therein).  Suppose that there exists
an isograded algebra isomorphism $\ph : \cJ \to \cJ'$ such that
\begin{equation}
\label{eq:mumu'}
\ph(\mu) = {z'}^2 \mu'
\end{equation}
for some nonzero homogeneous $z'\in \cZ'(\cJ')$.
Then $\CD(\cJ,\theta,\mu) \ig \CD(\cJ',\theta',\mu')$.
\end{lemma}

\begin{proof}
Let $\cA = \CD(\cJ,\theta,\mu)$ and $\cA' = \CD(\cJ',\theta',\mu')$.
We use the notation $\Lm$, $M$, $\cJ$, $\cZ$, $\Gm$, $\sg_0$, $\mu$ and $s_0$ for $\cA$
as in Theorem \ref{prop:CDtorus}(a), and we use corresponding primed notation
for $\cA'$.

Now by assumption we have a group isomorphism $\ep : M \to M'$ such that
\[\ph(\cJ^\tau) = {\cJ'}^{\ep(\tau)}\]
for $\tau\in M$.  Since $\ph$ is an algebra isomorphism it follows that
$\ep(\Gm) = \Gm'$, so $\ep(M\setminus \Gm) = M'\setminus \Gm'$. Thus,
\begin{equation}
\label{eq:phtheta}
\ph\theta = \theta' \ph.
\end{equation}

Next let $\gm'\in \Gm'$ be the degree of $z'\in {\cZ'}$.  Then, since
$\mu\in \cZ^{2\sg_0}$ and $\mu'\in {\cZ'}^{2\sg_0'}$, it follows from
\eqref{eq:mumu'} that
\begin{equation}
\label{eq:sgsg'}
\ep(2\sg_0) = 2\gm' + 2\sg_0'.
\end{equation}
We now define $\ep_\Lm : \Lm \to \Lm'$ by
\[\ep_\Lm(\tau) = \ep(\tau) \andd
\ep_\Lm(\sg_0 + \tau) = \sg_0' + \ep(\tau) + \gm'\]
for $\tau\in M$.  Then, one checks using
\eqref{eq:sgsg'} that $\ep_\Lm$ is an isomorphism of groups.

Next we define $\psi : \cA \to \cA'$ by
\[\psi(a+s_0b) = \ph(a) + s_0'(z'\ph(b))\]
for $a,b\in \cJ$.  Then, $\psi$ is a linear bijection, and
one checks directly using the definitions of $\psi$, $\ep_\Lm$
and the gradings on $\cA$ and $\cA'$ (see \eqref{eq:CDgrading}) that
\[\psi(\cA^\lm) = {\cA'}^{\ep_\Lm(\lm)}.\]
for $\lm\in \Lm$.  Also, since $z'\in \cZ'$, we have
\[z'{a}'^{\theta'} = (z'a')^{\theta'}\]
for $a'\in\cJ'$.  Using this fact and \eqref{eq:phtheta}, as well as the definitions
of the involutions and products, one checks directly that
$\psi$ preserves the products and  involutions.
\end{proof}

As mentioned previously, Yoshii  has classified Jordan tori in \cite{Y1}.
In view of Theorem \ref{prop:CDtorus}, it important for us to identity in his list
the Jordan $\Lm$-tori  $\cJ$ of degree~4 that satisfy $2\Lm \subset \Gm(\cJ)$.
(Of course we will then apply this with $\Lm$ replaced by $M$.)  To do this,
we recall some standard terminology from Jordan theory.

If $\cA$ is an associative algebra, the \emph{plus algebra} of $\cA$ is the algebra $\cA^+$ with underlying
vector space $\cA$ and product $x\cdot y = \frac 12 (xy+yx)$.  In that case, $\cA^+$ is a Jordan algebra.
If $\cA$ is an associative algebra with involution $*$  then
\[\cH(\cA) := \cA_+ = \set{ x\in \cA \suchthat x^* = x}\]
is a subalgebra of the Jordan algebra $\cA^+$.
If we wish to emphasize the role of the involution $*$, we write
$\cH(\cA)$ as $\cH(\cA,*)$.  If $\cA$ is a $\Lm$-graded algebra, then $\cA^+$ is a $\Lm$-graded algebra (with the
same grading).  Also, if $\cA$ is a $\Lm$-graded algebra with involution, then $\cH(\cA)$ is a $\Lm$-graded
subalgebra of $\cA^+$.  Finally, if $\cA$ is an associative $\Lm$-torus (without involution), then $\cA^+$ is a Jordan
$\Lm$-torus; and if $\cA$ is an associative $\Lm$-torus with involution such that $\langle S_+(\cA)\rangle = \Lm$,
then $\cH(\cA)$ is a Jordan $\Lm$-torus

If $E$ is a split 2-dimensional composition algebra (with its canonical involution)
over $F$ and $\cA$ is a
$\Lm$-graded algebra with involution $*$, then it is well known that
\begin{equation}
\label{eq:splitE}
\cH(\cA\ot E) \simeq_\Lm \cA^+.
\end{equation}
(If we identify $E = F\oplus F$, an isomorphism from right to left
in \eqref{eq:splitE} is $a \mapsto a\ot (1,0) + a^* \ot (0,1)$.)

\begin{proposition}
\label{prop:Jord4}
Suppose that $\cJ$ is a $\Lm$-graded Jordan algebra.
Then,
\begin{equation}
\label{eq:Jord4}
\textit{$\cJ$ is a Jordan $\Lm$-torus of degree 4 satisfying $2 \Lm \subset \Gm(\cJ)$}
\end{equation}
if and only if\/
$\cJ$ is isograded  isomorphic to
\begin{equation*}
\label{eq:Jord4list}
\cH(\Cay(2)\ot\Cay(2)\ot \cT \ot \cP(r)),
\end{equation*}
where $r\ge 0$, $\cT = \Cay(0)$, $\Cay(1)$, $\Cay(2)$ or $E$, and
$E$ is a 2-dimensional composition algebra over $F$.
\end{proposition}

\begin{proof}
``$\Leftarrow$'' We can assume that
$\cJ = \cH(\cA)$, where
\[\cA = \cA_1\ot\cA_2\ot\cA_3\ot\cA_4,\]
with $\cA_1 = \Cay(2)$, $\cA_2 = \Cay(2)$,
\begin{equation}
\label{eq:A3choice}
\cA_3 = \Cay(0),\ \Cay(1),\ \Cay(2) \orr E,
\end{equation}
and $\cA_4 = \cP(r)$ with $r\ge 0$. Then $\cA$ is a $\bbZ^n$-graded
associative algebra with involution, where $n = 2+2+k+r$, with $k =
0$, $1$, $2$ or $0$ in the cases covered by \eqref{eq:A3choice} (in
order). (Here we are making the obvious identification of
$\bbZ^2\oplus \bbZ^2 \oplus \bbZ^k \oplus \bbZ^r$ with  $\bbZ^n$.)
So $\cJ$ is a Jordan $\Zn$-torus by Remark \ref{rem:tensor}(a).  (To
see this in the case when $\cT = E$, we can extend the base field,
assume that $E$ is split, and use  \eqref{eq:splitE}.) Next
let $\cZ_i = \cZ(\cA_i)$ for $1\le i \le 4$, and let $\cZ =
\cZ(\cJ)$.  Using the fact that $\cZ$ is a graded subspace of $\cJ$,
it is easy to check that $\cZ = \cZ_1\ot \cZ_2\ot  \cZ_3\ot \cZ_4$.
Thus $\Gm(\cJ) = 2\bbZ^2\oplus 2\bbZ^2 \oplus 2\bbZ^k\oplus \bbZ^r$,
so $2\Zn \subset \Gm(\cJ)$. Finally, let $\tcZ$ be the quotient
field of~$\cZ$. Then, by Lemma \ref{lem:baseringtensor} (extended to
more than two factors), we have
\[
(\tcZ\ot_{\cZ_1} \cA_1)
\ot_\tcZ
(\tcZ\ot_{\cZ_2} \cA_2)
\ot_\tcZ
(\tcZ\ot_{\cZ_3} \cA_3)
\ot_\tcZ
(\tcZ\ot_{\cZ_4} \cA_4)
\simeq
\tcZ\ot_\cZ \cA
\]
as algebras with involution over $\tcZ$.
So, since $\cA_4 = \cZ_4$, we have
\[\tcZ\ot_\cZ \cJ = \tcZ\ot_\cZ \cH(\cA)= \cH(\tcZ\ot_\cZ \cA)
\simeq
\cH((\tcZ\ot_{\cZ_1} \cA_1)\ot_\tcZ(\tcZ\ot_{\cZ_2} \cA_2)\ot_\tcZ(\tcZ\ot_{\cZ_3} \cA_3)).
 \]
as algebras over $\tcZ$.
Finally, since $\tcZ\ot_{\cZ_1} \cA_1$ and $\tcZ\ot_{\cZ_2} \cA_2$ are
quaternion algebras with canonical involution over $\tcZ$ and
$\tcZ\ot_{\cZ_3} \cA_3$ is a composition algebra with canonical involution
of dimension $1$, $2$, $4$ or $2$ over $\tcZ$, it is well known that the righthand side of this
isomorphism is a degree 4 Jordan algebra over $\tcZ$ \cite[\S V.7, Theorem 11]{J2}.

``$\Rightarrow$'' Suppose that \eqref{eq:Jord4} holds, and let $\cZ = \cZ(\cJ)$ and $\Gm = \Gm(\cJ)$.
Then the central closure $\tcJ = \tcZ \ot_{\cJ} \cJ$
of $\cJ$ is a finite dimensional central simple Jordan division algebra of degree $4$ over $\tcJ$.
Thus, the dimension $\tcJ$ over $\tcZ$ is $10$, $16$ or $28$ \ibid.  So
\[d := \ord{S(\cJ)/\Gm} =10, 16 \text{ or } 28.\]

We now use the classification of Jordan tori of degree $\ge 4$.
Indeed, since $\tcJ$ has degree $\ge 4$, \cite[Theorem 7.1]{Y1} tells us that
one of the following holds:
\begin{mylist}
\item[(i)] $\cJ \simeq_\Lm \cA^+$, where $\cA$ is an associative $\Lm$-torus,
\item[(ii)] $\cJ \simeq_\Lm \cH(\cA)$, where $\cA$ is an associative $\Lm$-torus with involution,
\item[(iii)] $\cJ \simeq_\Lm \cH(\cB,\sg)$, where
$\cB$ is an associative $\Lm$-torus over $E$, $E/F$ is a
quadratic field extension,
and $\sg$ is a graded $\sg_E$-semilinear involution of $\cB$.
\end{mylist}
(In fact \cite[Theorem 7.1]{Y1} says more in each case, but this is all we need.)
We consider the cases (i), (ii) and (iii) separately.

(i)\ms Suppose that $\cJ = \cA^+$, where $\cA$ is an associative $\Lm$-torus.
Then $S(\cJ) = \Lm$, so $d = \ord{\Lm/\Gm}$ is a power of $2$. Thus $d = 16$.

Since $2\Lm\subset \Gm$, we may choose a basis $\lm_1,\ldots,\lm_r$
for $\Lm$ such that $n_1\lm_1,\ldots,n_r\lm_r$ is a basis for $\Gm$
with $n_i = 1$ or $2$. We then choose $0\ne x_i\in \cA^{\lm_i}$, in
which case the elements $x_1^{\pm 1},\dots,x_n^{\pm 1}$ generate the
algebra $\cA$ and satisfy $x_jx_i = q_{ij}x_j x_i$, where $\bq =
(q_{ij}) \in F_{n\times n}$ is a \emph{quantum matrix} (that is
$q_{ii} = 1$ and $q_{ij} = q_{ji}^{-1}$).  In other words, $\cA$ is the
\textit{quantum torus} $F_\bq$ (see for example \cite[p.129]{Y1}).
Now, since $\cJ$ is a Jordan torus, the centre $\cZ$ of the Jordan
algebra $\cJ=\cA^+$ is also the centre of the associative algebra
$\cA$ (without involution) \cite[Lemma 3.6(i) and 2.5]{Y1}. So the
elements  $x_j^{n_j}\in \cZ$ commute with the elements $x_i$ of
$\cA$. But $x_j^{n_j}x_i=q_{ij}^{n_j}x_ix_j^{n_j}$, so $q_{ij}^{n_j}
= 1$.  Thus, $q_{ij} = \pm 1$; that is, $\bq$ is an
\emph{elementary} quantum matrix.  Therefore $\cA$ has a graded
involution $*$ such that $x_i^* = x_i$ \cite[Example 4.3]{Y1}. So,
using the classification of associative tori \textit{with}
involution (see Remark \ref{rem:assoctori}), we conclude that we
have
\[
\cA \ig \cA_1\ot\ldots\ot\cA_p \ot \cP(r),\] as graded algebras
\emph{without} involution, where $p,r \ge 0$ and $\cA_i = \Cay(2)$
for each~$i$. (Here we have used the fact that $\Cay_*(2) \ig
\Cay(2)$ and $\Cay(1) \ig \cP(1)$ as graded algebras without
involution.) But $\Gm$ is the support of the centre of the
associative algebra $\cA$ (without involution), and therefore the
same is true for $\cA_1\ot\ldots\ot\cA_p \ot \cP(r)$. Thus $d =
\ord{\Lm/\Gm} =  4^p$.  So, $p = 2$, and hence $\cA \ig
\cA_1\ot\cA_2\ot \cP(r)$ as graded algebras without involution. Thus
$\cJ \ig \left(\Cay(2)\ot\Cay(2) \ot \cP(r)\right)^+$. So, by
\eqref{eq:splitE},
 $\cJ \ig \cH(\Cay(2)\ot\Cay(2)\ot E\ot \cP(r))$,
where $E$ is the split 2-dimensional composition algebra.

(ii) Suppose $\cJ = \cH(\cA)$, where $\cA$ is an associative $\Lm$-torus with involution.
By Remark \ref{rem:assoctori}, we may assume that $\Lm = \Zn$ and
\[
\cA = \cA_{1}\ot\ldots\ot\cA_{l}\ot\cA_{l+1}\ot \cP(r)
\]
where $\ell, r\ge 0$, $\cA_{i}=\Cay(2)$ for $1\leq i\leq l$,  and
$\cA_{l+1}$ is $\Cay(0)$, $\Cay(1)$, or $\Cay_*(2)$. From this it is
easy to check that $\cZ = \cZ_1\ot \ldots\ot\cZ_{\ell+1}\ot \cP(r)$,
where $\cZ_i$ is the centre of $\cA_i$ for
$1\le i \le \ell+1$.  If $\cC$ and
$\cD$ are algebras with involution, then
\[
\cH(\cC\otimes\cD)=(\cC\otimes\cD)_{+}
=(\cC_{+}\otimes\cD_{+})\oplus(\cC_{-}\otimes\cD_{-}).
\]
Thus, we can use induction on $l$ to compute
$d=\left\vert S(\cJ)/\Gamma\right\vert $ getting
\[
\begin{tabular}
[c]{|c|c|c|c|}\hline
& $\cA_{l+1}=\cC(0)$ & $\cC(1)$ & $\cC_{\ast
}(2)$\\\hline
$l=0$ & $1$ & $1$ & $3$\\\hline
$1$ & $1$ & $4$ & $6$\\\hline
$2$ & $10$ & $16$ & $36$\\\hline
$3$ & $28$ & $64$ & \\\hline
$4$ & $136$ &  & \\\hline
\end{tabular}
.
\]
Since $d$ increases with increasing $l$ and since $d=10,16$, or $28$, we
 have only the following
possibilities: $\ell = 2$ and $\cA_{\ell+1} = F$; $\ell = 3$ and
$\cA_{\ell+1} = F$; or $\ell = 2$ and $\cA_{\ell+1} = \Cay(1)$.
Therefore, $\cJ \ig \cH(\Cay(2)\ot\Cay(2)\ot\Cay(k)\ot\cP(r))$,
where $0\le k \le 2$.

(iii)\ms Suppose that $\cJ = \cH(\cB,\sg)$, with
$\cB$, $E/F$ and $\sg$ as above in (iii).
Now $E\ot \cJ \simeq \cB^+$ as graded  Jordan algebras \cite[Example 4.3(3)]{Y1}, and hence
$\cB^+$ is a Jordan $\Lm$-torus over $E$ that satisfies the condition
\eqref{eq:Jord4}.   Thus,  by the argument in (i) (with $F$ replaced by $E$
and $\cA$ replaced by $\cB$), we may identify
\begin{equation*}
\label{eq:bident}
\cB = \cB_1\ot_E \cB_2\ot_E \cB_3,
\end{equation*}
as $\Lm$-graded algebras without involution over $E$, where
$\Lm = \Lm_1\oplus\Lm_2\oplus\Lm_3$,
$\cB_i$ is a $\Lm_i$-graded associative algebra without involution over $E$ for $1\le i \le 3$,
\begin{equation*}
\label{eq:bcompig}
\cB_1 \ig E\ot \Cay(2),\quad \cB_2 \ig E\ot \Cay(2) \andd
\cB_3  \ig E\ot \cP(r)
\end{equation*}
as graded associative algebras without involution for some $r\ge 0$.
Since $\sg$ is graded, it has the form $\sg = \sg_1\ot\sg_2\ot \sg_3$,
where $\sg_i$ is a $\sg_E$-semilinear involution of $\cB_i$ for $1\le i \le 3$.
But for $\lm_i\in \Lm_i$ we have $\cB_i^{\lambda_i} \ne 0$ and
therefore (since $\sg_i|_{\cB_i^{\lambda_i}}$ is a graded $\sg_E$-semilinear map of period 2)
$\cB_i$ contains an nonzero element fixed
by $\sg_i$ and a nonzero element antifixed by~$\sg_i$.
Hence, rechoosing canonical generators, we have
$(\cB_i,\sg_i) \ig (E\ot \Cay(2),\sg_E\ot \natural)$
for $i = 1,2$, and
$(\cB_3,\sg_3) \ig (E\ot \cP(r),\sg_E\ot 1)$
as graded algebras with involution over $F$.
Thus, as graded algebras with involution over $F$, we
have
\begin{align*}
(\cB,\sg)
&\ig \big(E\ot (\Cay(2) \ot  \Cay(2) \ot \cP(r)), \sg_E\ot (\natural\ot \natural \ot 1)\big)\\
&\ig (\Cay(2) \ot \Cay(2) \ot E \ot \cP(r), \natural \ot \natural \ot \sg_E \ot 1)
\end{align*}
So  $\cJ \ig \cH(\Cay(2)\ot\Cay(2)\ot E \ot \cP(r))$.
\end{proof}

We now put together our results to give a precise description of structurable
tori of class III(b).

\begin{theorem}
\label{thm:IIIb}
(a)\ms Let
 \[\cJ = \cH(\Cay(2)\ot\Cay(2)  \ot \cT \ot \cP(r)),\]
where $r\ge 1$,
\begin{equation}
\label{eq:choiceT}
\cT = \Cay(0),\ \Cay(1),\ \Cay(2)\ \text{or } E,
\end{equation}
and $E$ is a 2-dimensional composition algebra over $F$.
Set $k = 0$, $1$, $2$  or  $0$ in the cases covered by \eqref{eq:choiceT} in order; and
set $q = 4+k$ and $n = q +r$.
Let $\Lm$ be a free abelian group with
basis $\lm_1,\dots,\lm_n$, and let
\[M = \langle \lm_1,\dots,\lm_{q},2\lm_{q+1},\lm_{q+2},\dots,\lm_n\rangle.\]
Decompose $M$ as $M = M_1\oplus M_2 \oplus M_3 \oplus M_4$, where
\[M_1 = \langle \lm_1,\lm_2\rangle,\ M_2 = \langle \lm_3,\lm_4\rangle,\
M_3 = \langle \lm_5,\dots,\lm_q\rangle,\ M_4 = \langle 2\lm_{q+1},\lm_{q+1},\dots,\lm_n\rangle.\]
(So $M_3 = 0$ if $\cT = \Cay(0)$ or $E$.)
Give the associative algebra with involution $\Cay(2)\ot\Cay(2)  \ot \cT \ot \cP(r)$ the tensor product
$M$-grading, where the $M_1$, $M_2$, $M_3$ and $M_4$ gradings on the factors
$\Cay(2)$, $\Cay(2)$, $\cT$ and $\cP(r)$ are determined respectively by the
bases $\set{\lm_1,\lm_2}$, $\set{\lm_3,\lm_4}$,
$\set{\lm_5,\dots,\lm_q}$ and $\set{2\lm_{q+1},\lm_{q+1},\dots,\lm_n}$ (see
Remark \ref{rem:altgrading}).
Give the Jordan algebra $\cJ$ an $M$-grading as a
graded subalgebra of $\left(\Cay(2)\ot\Cay(2)  \ot \cT \ot \cP(r)\right)^+$.
Set
\[\sg_0 = \lm_{q+1} \in \Lm \andd \mu = 1\ot 1 \ot 1 \ot t_1 \in \cZ(\cJ)^{2\sg_0}.\]
Define $\theta : \cJ \to \cJ$ by \eqref{eq:deftheta} and construct the $\Lm$-graded
algebra with involution
\[\CD(\cJ,\theta,\mu) = \cJ  \oplus s_0 \cJ,\]
with multiplication, involution and grading given by  \eqref{eq:CDmult}, \eqref{eq:CDinvol}
and \eqref{eq:CDgrading}.  Then, $\CD(\cJ,\theta,\mu)$ is a structurable $\Lm$-torus of class III(b) with centre
$\cZ(\cJ)$.

(b)\ms  Any structurable $\Lm$-torus of class III(b) is isograded isomorphic to a Cayley-Dickson torus
$\CD(\cJ,\theta,\mu)$ constructed as in (a).
\end{theorem}

\begin{proof} (a): By Proposition \ref{prop:Jord4},
$\cJ$ is an $M$-graded Jordan torus  of degree 4 with $2M \subset \Gm(\cJ)$.
Then the conclusion follows from Theorem \ref{prop:CDtorus}(a).

(b): Suppose that $\cA$ is a structurable $\Lm$-torus of class III(b).  We use the notation of
\S\ref{sec:classes} and \S\ref{sec:geometry} for $\cA$.  Choose $\sg_0\in S_-$ and, as in the proof
of Theorem \ref{prop:CDtorus}(b), choose a subgroup $M$ of $\Lm$ such that $\Gm \subset M$,
$(\Lm:M) = 2$ and $\sg_0\notin M$. Choose $0\ne s_0\in \cA^{\sg_0}$, let $\mu = s_0^2 \in \cZ^{2\sg_0}$,
and let $\cJ = \cA^M$.  Then, by the proof
of Theorem~\ref{prop:CDtorus}(b), $\cJ$ is a Jordan $M$-torus of degree 4 with
$\cZ(\cJ) = \cZ$, $\Gm(\cJ) = \Gm$ and $2M \subset \Gm$.  Moreover,
$\cA = \CD(\cJ,\theta,\mu)$ as $\Lm$-graded algebras with involution, where
$\theta: \cJ \to \cJ$ is defined by
\eqref{eq:deftheta} and the $\Lm$-grading on $\cA$ is defined by
\eqref{eq:CDgrading}.

Now, by Proposition
\ref{prop:Jord4}, $\cJ$ is isograded isomorphic to
$\cH(\Cay(2)\ot\Cay(2)  \ot \cT \ot \cP(r))$, where $r\ge 0$,
\begin{equation*}
\label{eq:choiceT'}
\cT = \Cay(0),\ \Cay(1),\ \Cay(2)\ \text{or } E,
\end{equation*}
and $E$ is a 2-dimensional composition algebra.
Let $k = 0$, $1$, $2$ or $0$ respectively; and let  $q = 4+k$ and $n=q+r$.
Then, we may identify
$M =\bbZ^{2+2+k+r} =  \bbZ^n$ and
\[\cJ = \cH(\Cay(2)\ot\Cay(2)  \ot \cT \ot \cP(r)).\]

Next, as we saw in the proof of Proposition \ref{prop:Jord4}(a),
we have $\cZ = \cP(2) \ot \cP(2) \ot \cZ(\cT)\ot \cP(r)$, with
$\cZ(\cT) =  \cP(k)$ if $\cT =\cA(k)$ and $\cZ(\cT) = F$ if $\cT = E$.
So, letting $\ep_1,\ldots,\ep_n$ be the standard basis of $M = \bbZ^n$, we have
\[\Gm =  \sum_{i=1}^q 2\bbZ \ep_i + \sum_{i=q+1}^n \bbZ \ep_i.\]

Also $\sg_0 \notin M$, so $2\sg_0 \notin 2M$.  Thus $2M \ne \Gm$, so $r\ge 1$.

Now let $\lm_1,\dots,\lm_n$ (in $(\frac 12 \bbZ)^n$) be defined by
$\lm_i = \ep_i$ if $i\ne q+1$ and $\lm_{q+1} = \frac 12 \ep_{q+1}$.
Let $\Lm' = \langle \lm_1,\dots,\lm_n\rangle$, in which case
$M = \langle \lm_1,\dots,\lm_{q},2\lm_{q+1},\lm_{q+2},\dots,\lm_n\rangle$ in $\Lm'$.
Let
\[\sg_0' = \lm_{q+1}\in \Lm' \andd \mu' = 1\ot 1 \ot 1\ot t_1 \in \bbZ(\cJ)^{2\sg_0'}.\]
Construct the Cayley-Dickson $\Lm'$-torus $\CD(\cJ,\theta,\mu')$ as in (a)
(using $\Lm'$, $\sg_0'$ and $\mu'$ rather that $\Lm$, $\sg_0$ and $\mu$).
To complete the proof of (b), we show that
$\CD(\cJ,\theta,\mu') \ig \CD(\cJ,\theta,\mu)$.  To do this it suffices, by
Lemma \ref{lem:mu}, to show that
there exists an isograded  automorphism
$\ph : \cJ \to \cJ$ such that
\[\ph(\mu') = z^2 \mu\]
for some nonzero homogeneous $z\in \cZ(\cJ)$.

Now, since $2\sg_0\in \Gm$, we can write
\[2\sg_0 = \sum_{i=1}^q 2a_i \ep_i + \sum_{i=q+1}^n a_i \ep_i,\]
where $a_i\in \bbZ$ for all $i$.  Moreover, since $2\sg_0\notin 2M$, $a_{q+j}$ is odd for some $j\ge 1$.
Replacing $\mu$ by $z^2\mu$ for some $z\in \cZ$, we can assume that $a_{q+j} = 1$.
Next, there exists an isograded automorphism
$\psi$ of $\cP(r)$ that permutes the canonical generators and exchanges $t_1$
and $t_j$.  Replacing $\mu$ by $(1\ot 1\ot 1\ot \psi)(\mu)$, we can assume that
$a_{q+1} = 1$.  So we have
\begin{equation}
\label{eq:zlocation}
w := \mu(\mu')^{-1} \in \cZ^{\sum_{i\ne q+1}\bbZ \ep_i}.
\end{equation}
Define $\ph : \cJ \to \cJ$ by $\ph(x) = w^m x$ for $x\in  \cJ^{m\ep_{q+1} + \sum_{i\ne q+1}\bbZ \ep_i}$,
$m\in \bbZ$.  Then, one checks directly that $\ph$ is an algebra homomorphism, and
(using \eqref{eq:zlocation}) that
$\ph$ is invertible with inverse defined by
$x \mapsto  w^{-m} x$ for $x\in  \cJ^{m\ep_{q+1} + \sum_{i\ne q+1}\bbZ \ep_i}$,
$m\in \bbZ$.   Finally, it is clear that $\ph(\mu') = \mu$ and that
$\ph$ is isograded.
\end{proof}

\begin{remark}
\label{rem:splitE} In Theorem \ref{thm:IIIb}(a), suppose
that $\cT = E$.  If $E$ split (which holds automatically
when $F$ is algebraically closed) then instead of choosing
$\cJ = \cH(\Cay(2)\ot\Cay(2)\ot E \ot \cP(r))$ we may, by \eqref{eq:splitE}, choose
\[\cJ = \left(\Cay(2)\ot \Cay(2)\ot\cP(r)\right)^+\]
with the tensor product grading by $M = M_1\oplus M_2\oplus M_4$.
In that case our choice of $\mu$ is $\mu = 1\ot 1 \ot t_1$.
\end{remark}

\section{The space of hermitian matrices $\Herm$}
\label{sec:hermitian}

In the  next section we will construct structurable tori of class III(c)
using a cubic form defined on the space of hermitian matrices $\Herm$
over a composition algebra $\cC$.  To prepare for that we prove two
lemmas about $\Herm$ in this section, a coordinatization lemma and a lemma about gradings.

Throughout the section, we assume that \emph{$K$ is a commutative associative
algebra over $F$.}

If $\cC$ is a composition algebra over $K$ with canonical involution
$\barinv$, we  let $\Herm$ denote the $K$-module of all
$3\times3$-hermitian matrices
\[
x=\left[
\begin{array}
[c]{ccc}%
a_{1} & x_{3} & \bar{x}_{2}\\
\bar{x}_{3} & a_{2} & x_{1}\\
x_{2} & \bar{x}_{1} & a_{3}%
\end{array}
\right]  .
\]
over $\cC$, where $a_{i}\in K$, $x_{i}\in\cC$. We can also
write the above element $x\in\Herm$ as
\[
x={\sum_{i}}a_{i}[ii]+{\sum_{(i,j,k)\cyc}}x_{i}[jk],
\]
where $(i,j,k)\cyc$ means that $(i,j,k)$ is a cyclic permutation of
$(1,2,3)$. For this $x$ and for $y=\textstyle{\sum_{i}}b_{i}[ii]+{\sum
_{(i,j,k)\cyc}}y_{i}[jk]$, we set
\begin{align}
T(x,y) &  =\sum_{i}a_{i}b_{i}+\sum_{i}n(x_{i},y_{i})\in K,\label{eq:trace}\\
N(x) &  =a_{1}a_{2}a_{3}-a_{1}n(x_{1})-a_{2}n(x_{2})-a_{3}n(x_{3}%
)+t(x_{1}x_{2}x_{3})\in K,\label{eq:norm}
\end{align}
where $n$ and $t$ denote the norm and trace on $\cC$ respectively.
Then, $T$ is a nondegenerate $K$-bilinear form on $\Herm$
and $N$ is a cubic form on $\Herm$ over $K$. We
call $T$ and $N$ the \textit{standard trace and norm} on $\Herm$.
It is well known \cite[p.~488]{Mc} that the pair $(T,N)$
satisfies the adjoint identity with adjoint
\[
x^{\#}={\sum_{(i,j,k)\cyc}}(a_{i}a_{j}-n(x_{k}))[kk]+{\sum
_{(i,j,k)\cyc}}(\overline{x_{i}x_{j}}-a_{k}x_{k})[ij].
\]
Hence, for $x$ and $y$ as above in $\Herm$, we have
\begin{equation}
\label{eq:adjointlin} x\times y = {\sum_{(i,j,k)\cyc}}
(a_ib_j + b_ia_j -n(x_k,y_k))[kk]+ {\sum_{(i,j,k)\cyc}}
(\overline{x_iy_j}+\overline{y_ix_j}-a_ky_k - b_kx_k)[ij].
\end{equation}

We now prove a coordinatization result for a pair $(\tilde{T},\tilde{N})$
satisfying the adjoint identity.

\begin{lemma}
\label{lem:normH}
For $i=1,2,3$, let $\cW_{i}$ be a $K$-module with a  quadratic form
$n_{i}:\cW_{i}\rightarrow K$ and let $\tau:\cW_{1}
\times\cW_{2}\times\cW_{3}\rightarrow K$ be trilinear. On
$\cW=K^{3}\oplus\cW_{1}\oplus\cW_{2}\oplus
\cW_{3}$, define%
\begin{align*}
\tilde{T}(x,y)  & =\sum_{i}a_{i}b_{i}+\sum_{i}n_{i}(x_{i},y_{i}),\\
\tilde{N}(x)  & =a_{1}a_{2}a_{3}-a_{1}n_{1}(x_{1})-a_{2}n_{2}(x_{2}%
)-a_{3}n_{3}(x_{3})+\tau(x_{1},x_{2},x_{3}),
\end{align*}
for
$x=(a_{1},a_{2},a_{3},x_{1},x_{2},x_{3})$ and $y=(b_{1},b_{2},b_{3}%
,y_{1},y_{2},y_{3})$.
Suppose $u_{i}\in\cW_{i}$ with $n_{i}(u_{i})=1$ for
$i=1,2$,  $\tilde T$ is nondegenerate, and
the pair $(\tilde{T},\tilde{N})$ satisfies the adjoint
identity. Then there is a composition algebra $\cC$ and $K$-linear
isomorphisms $\eta_{i}:\cC\rightarrow\cW_{i}$ such that
$\eta_{i}(1)=u_{i}$, where $u_{3}=u_{1}\times u_{2}$, and such that the $K$-linear
isomorphism $\eta:\Herm\rightarrow\cW$ given by
\begin{equation}
\label{eq:defeta}
\eta({\sum_{i}}a_{i}[ii]+{\sum_{(i,j,k)\cyc}}x_{i}[jk])=(a_{1},a_{2},a_{3},\eta_{1}%
(x_{1}),\eta_{2}(x_{2}),\eta_{3}(x_{3})),
\end{equation}
satisfies $\tilde{T}(\eta(x),\eta(y))=T(x,y)$ and $\tilde{N}(\eta(x))=N(x) $
for all $x,y\in\Herm$.
\end{lemma}

\begin{proof}
Since $\tilde{T}$ is nondegenerate, each $n_{i}$ is nondegenerate. Using
$\partial_{y}\tilde{N}|_{x}=\tilde{T}(y,x^{\#})$ to compute $x^{\#}$
for $x=(a_{1},a_{2},a_{3},x_{1},x_{2},x_{3})$, we get
\begin{gather*}
x^{\#}   = (c_{1},c_{2},c_{3},z_{1},z_{2},z_{3}),\quad\text{where}\\
c_{k}   = a_{i}a_{j}-n_{k}(x_{k})\in K, \quad z_{k} = x_{i}\times x_{j}-a_{k}x_{k} \in\cW_{k}
\end{gather*}
for $(i,j,k)\cyc$. In particular, if $e_{1},e_{2},e_{3}$ is the
standard basis of $K^{3}$ and $x_{k}\in\cW_{k}$ for
$\{i,j,k\}=\{1,2,3\}$, then
\[e_{i}\times e_{j}   =e_{k},\quad x_{k}^{\#}  =-n_{k}(x_{k})e_{k},\quad  e_{k}\times x_{k}   =-x_{k},\]
so $\tilde{T}(x_{i},x_{j}^{\#})=\tilde{T}(x_{i},-n_{j}(x_{j})e_{j})=0$.
Thus, (ADJ4) for $(\tilde{T},\tilde{N})$ gives $(x_{i}\times x_{j})^{\#}+x_{i}^{\#}\times x_{j}^{\#}=0$;
that is
\begin{equation}
n_{k}(x_{i}\times x_{j})=n_{i}(x_{i})n_{j}(x_{j}).\label{eq:normcomp}%
\end{equation}
Since $n_{i}(u_{i})=1$ and $u_{i}^{\#}=-e_{i}$ for $i=1,2$, we can set
$u_{3}=u_{1}\times u_{2}\in\cW_{3}$ to get $n_{3}(u_{3})=1$ and
$u_{3}^{\#}=-e_{3}$.

If $x_{j},y_{j}\in\cW_{j}$ and $(i,j,k)\cyc$, then $n_{k}(u_{i}\times
x_{j},u_{i}\times y_{j})= n_i(u_i)n_j(x_j,y_y) = n_{j}(x_{j},y_{j})$ by the polarization of
(\ref{eq:normcomp}).  Since $n_{j}$ is nondegenerate, we see that
$y_{j}\rightarrow u_{i}\times y_{j}$ is a linear monomorphism of $\cW_j$ into $\cW_k$.
Polarizing (ADJ5) for $(\tilde{T},\tilde{N})$, we get%
\[
(u\times(v\times w))\times v+(u\times v^{\#})\times w=\tilde{T}(w,v^{\#}%
)u+\tilde{T}(u,w)v^{\#}+\tilde{T}(u,v)(v\times w)\text{.}%
\]
Taking $u=v=u_{i}$ and $w=x_{j}\in\cW_{j}$ with $i\neq j$, we get
\begin{align*}
u_{i}\times(u_{i}\times(u_{i}\times x_{j})))  & =-(u_{i}\times u_{i}%
^{\#})\times x_{j}+\tilde{T}(x_{j},u_{i}^{\#})u_{i}\\
&\qquad\quad +\tilde{T}(u_{i},x_{j})u_{i}^{\#}+\tilde{T}(u_{i},u_{i})(u_{i}\times
x_{j})\\
& =-u_{i}\times x_{j}+2(u_{i}\times x_{j})=u_{i}\times x_{j},
\end{align*}
since $\tilde{T}(x_{j},u_{i}^{\#})=\tilde{T}(u_{i},x_{j})=0$, $-(u_{i}\times
u_{i}^{\#})=u_{i}\times e_{i}=-u_{i}$, and $\tilde{T}(u_{i},u_{i}%
)=2n_{i}(u_{i})=2$.  Since $y_{j}\rightarrow u_{i}\times y_{j}$ is a linear
monomorphism, we have%
\begin{equation}
u_{i}\times(u_{i}\times x_{j})=x_{j}.\label{eq:ux}%
\end{equation}
Since $u_{3}=u_{1}\times u_{2}$ by definition, (\ref{eq:ux}) shows that
$u_{i}\times u_{j}=u_{k}$ for all $\{i,j,k\}=\{1,2,3\}$. \ Let $\cC%
=\cW_{3}$ and define a bilinear product on $\cC$ by%
\[
xy=(x\times u_{1})\times(u_{2}\times y)
\]
for $x,y\in\cC$.  Clearly, $1_{\cC}:=u_{3}$ is the identity
element for $\cC$ by (\ref{eq:ux}).    Moreover, $n:=n_{3}$ has
$n(1_{\cC})=1$ and $n(xy)=n(x)n(y)$ by (\ref{eq:normcomp}).  Thus,
$\cC$ is a composition algebra.
We denote the norm, trace and canonical involution
on $\cC$ by $n$, $t$ and $\barinv$ respectively.
Define
$\eta_{i}:\cC\rightarrow\cW_{i}$ by $\eta_{1}(x)=u_{2}\times\bar{x}$,
$\eta_{2}(x)=u_{1}\times\bar{x}$, and $\eta_{3}=\id$.
By (\ref{eq:ux}), $\eta_{i}$
is a $K$-linear isomorphism satisfying $\eta_i(1_\cC) = u_i$.
Now define $\eta$ by \eqref{eq:defeta}.
Then, if $x\in\cC$,  using \eqref{eq:normcomp} we have
$n_{1}(\eta(x[23])) = n_1(\eta_1(x)) = n_1(u_2\times \bar x)
= n_2(u_2)n_3(\bar x) = n(\bar x) = n(x)$, and similarly
$n_{2}(\eta(x[31])) = n(x)$.  Thus
\begin{equation}
\label{eq:nin}
n_{1}(\eta(x[23])=n_{2}(\eta(x[31]))= n_{3}(\eta(x[12])) =n(x)
\end{equation}
for $x\in \cC$.
Next, if $x_i\in \cW_i$, we have
$\tau(x_1,x_2,x_3) = \partial_{x_1} \partial_{x_2} \tilde{N}|_{x_3} =
\tilde{T}(x_1\times x_2,x_3)= n_3(x_1\times x_2,x_3)$.
So for
$x,y,z\in\cC$,
\begin{align}
\tau(\eta(x[23]),\eta(y[31]),\eta(z[12]))  & =n_{3}(\eta(x[23])\times
\eta(y[31]),\eta(z[12]))\notag \\
& =n_{3}((\bar{y}\times u_{1})\times(u_{2}\times\bar{x}),z)\notag \\
& =n(\bar{y}\bar{x},z)=n(\overline{xy},z)=t(xyz)\text{.} \label{eq:taut}
\end{align}
Thus for $x,y\in\Herm$, we see, using \eqref{eq:nin} and
\eqref{eq:taut}, that $\tilde{N}%
(\eta(x))=N(x)$ and $\tilde{T}(\eta(x),\eta(y))=T(x,y)$.
\end{proof}

We next characterize certain gradings of $\Herm$.

\begin{lemma}
\label{lem:gradH}  Suppose that $K$ is $\Lm$-graded as an algebra,
and suppose $\cC$ is a composition algebra over $K$.

(a) If $\cC$ is $\Lm$-graded
as an algebra with involution and as a $K$-module and if
$\lambda_{i}\in\Lm$ with $\lambda_{1}+\lambda_{2}+\lambda_{3}=0$,
then letting
\begin{equation}
\label{eq:gradH}
a[ij]\in\Herm\quad\text{have degree}\quad\lambda
_{i}+\lambda_{j}+\mu,
\end{equation}
for $a\in\cC^{\mu}$, $\mu\in\Lm$, $(i,j,k)\cyc$, and
for $a\in K^{\mu}$, $\mu\in\Lm$, $i=j$, defines a $\Lm$-grading of
$\Herm$ as a $K$-module for which $N$ is a graded
cubic form.

(b) Conversely, suppose that we are given a $\Lm$-grading of $\Herm$
as a $K$-module with the properties that $N$ is a graded cubic form,
$1[ij]$ is homogeneous for $(i,j,k)\cyc$ or $i=j$, and
each $\cC[ij]$ with $(i,j,k)\cyc$ is a graded submodule of $\Herm$.
Then, there exists a $\Lm$-grading on $\cC$ as an algebra with involution and as a $K$-module
so that  the grading on $\Herm$ is given by $\eqref{eq:gradH}$ where
$-\lm_i$ is the degree of  $1[jk]$ for $(i,j,k)\cyc$.
\end{lemma}

\begin{proof}
(a): Clearly, (\ref{eq:gradH}) gives a $\Lm$-grading of
$\Herm$ as an $K$-module. \ Also, $N\ $is a graded cubic form.
Indeed, each term in the full linearization of $N$ involves one factor from
each row and one from each column, so the contribution to the degree from the
matrix locations is $2(\lambda_{1}+\lambda_{2}+\lambda_{3})=0$.

(b): Given a $\Lm$-grading of $\Herm$ as specified in
the converse, let $-\lambda_{i}$ be the degree of $1[jk]$ for
$(i,j,k)\cyc$ and let $\rho_{i}$ be the degree of $1[ii]$. \ We
can obtain three gradings of $\cC$ as a $K$-module by translating the
gradings of the submodules $\cC[jk]$. \ Specifically, let
$\cC_{(i)}=\cC$ with the grading given by $\cC%
_{(i)}^{\mu}[jk]=(\cC[jk])^{-\lambda_{i}+\mu}$ for
$(i,j,k)\cyc$. \ Thus, $1\in\cC_{(i)}^{0}$. \ If $x_{i}%
\in\cC_{(i)}^{\mu_{i}}$, we have%
\[
t(x_{1}x_{2}x_{3})=T(x_{1}[23],x_{2}[31]\times x_{3}[12])\in K^{-(\lambda
_{1}+\lambda_{2}+\lambda_{3})+(\mu_{1}+\mu_{2}+\mu_{3})}.
\]
Taking each $x_{i}=1$ shows $\lambda_{1}+\lambda_{2}+\lambda_{3}=0$, so the
trilinear map $\cC_{(1)}\times\cC_{(2)}\times\cC%
_{(3)}\rightarrow K$ with $(x_{1},x_{2},x_{3})\rightarrow t(x_{1}x_{2}x_{3}) $
is graded. \ Also, if $x\in\cC_{(i)}^{\lambda} $ and $y\in
\cC_{(i)}^{\mu}$, then%
\[
n(x,y)=-T(1[ii],x[jk]\times y[jk])\in K^{\rho_{i}-2\lambda_{i}+\lambda+\mu
}\text{.}%
\]
Again $x=y=1$ shows $\rho_{i}=2\lambda_{i}$, so $n:\cC_{(i)}%
\times\cC_{(i)}\rightarrow K$ and $t:\cC_{(i)}\rightarrow K$
are graded. \ Thus, $\bar{x}=t(x)1-x\in\cC_{(i)}^{\lambda}$, so $-$ is
graded on $\cC_{(i)}$. \ Since $n$ is nondegenerate and graded on
$\cC_{(i)}$, it is easy to see that $\cC_{(i)}^{\mu}%
=\{x\in\cC\mid n(\cC_{(i)}^{\lambda},x)\in K^{\lambda+\mu}\}$.
\ Now if $x\in\cC_{(j)}^{\mu}\cC_{(k)}^{\nu}$, then%
\[
n(\cC_{(i)}^{\lambda},x)=n(\overline{\cC_{(i)}^{\lambda}%
},x)\subset t(\cC_{(i)}^{\lambda}\cC_{(j)}^{\mu}%
\cC_{(k)}^{\nu})\subset K^{\lambda+\mu+\nu}%
\]
so $\cC_{(j)}^{\mu}\cC_{(k)}^{\nu}\subset\cC%
_{(i)}^{\mu+\nu}$. \ In particular, $\cC_{(1)}^{\lambda}%
=1\cC_{(3)}^{\lambda}=\cC_{(2)}^{\lambda}1$, so the three
gradings coincide. \ Now writing $\cC^{\lambda}=\cC%
_{(i)}^{\lambda}$, we see that $\cC$ is $\Lm$-graded as an algebra
with involution. \ Since $K^{\mu}[ii]=K^{\mu}(1[ii])\subset\cH%
(\cC_{3})^{2\lambda_{i}+\mu}$ and $\cC^{\mu}[jk]=
(\cC[jk])^{-\lambda_{i}+\mu}\subset\Herm^{-\lambda
_{j}-\lambda_{k}+\mu}$, we have the desired grading of
$\Herm$.
\end{proof}

\section{Construction of tori of class III(c)}
\label{sec:constructionIIIc}

In  this section  we construct structurable tori of class III(c) using the construction from
\S \ref{sec:Cubicforms} of a structurable algebra $\cA(h,N)$ from a pair
$(h,N)$ satisfying the adjoint identity. These tori will be
used in the next section to classify all tori of class III(c).

We will use the following proposition which can be viewed as a converse to   Proposition
\ref{prop:ClassIIIcubic}.

\begin{proposition}
\label{prop:torusIII}
Suppose $\cE$ is a $\Lm$-graded associative commutative
algebra with involution $*$ and centre $\cZ$; suppose $\cW$ is a $\Lm$-graded
left $\cE$-module;
suppose $h : \cW \times \cW \to  \cE$ is a non-degenerate $\Lm$-graded hermitian form; and
suppose $N : \cW \to \cE$ is a $\Lm$-graded cubic form over $\cE$
such
that $(h,N)$ satisfies the adjoint identity.
Let
\[\cA(h,N) = \cE \oplus \cW\]
be the algebra with involution
constructed from the pair $(h,N)$ as in \S\ref{sec:Cubicforms},
and give $\cA(h,N)$ a $\Lm$-grading by extending the given $\Lm$-gradings
on $\cE$ and $\cW$.  Then $\cA(h,N)$ is a $\Lm$-graded structurable algebra.
Moreover, if

(i)\ms the support $\Lm_-$ of $\cE$ is a
subgroup of $\Lm$,
and $\cE$ is a commutative associative $\Lm_{-}$-torus with involution,

(ii)\ms $\cW$ is finely graded,
$\supp(\cW)\cap \Lm_- = \emptyset$,
and  $\Lm = \langle \Lm_-, \supp(\cW)\rangle$,

(iii)\ms $4\Lm\subset\Gm(\cE)$,

(iv)\ms $0\neq x\in\cW^{\al}$ with $2\al\notin\Gm(\cE)$ implies
$x^\natural \neq0,$

\noindent then $\cA$ is a structurable $\Lm$-torus.
If, in addition,

(v)\ $N\ne 0$ and the involution $*$ on $\cE$ is nontrivial,

\noindent then $\cA$ is a structurable torus of class III and centre $\cZ$.

\end{proposition}

\begin{proof}
Let  $\cA=\cA(h,N)$. Then,
by Theorem \ref{thm:cubic}, $\cA$ is a structurable
algebra. \ To see that $\cA$ is $\Lm$-graded as an algebra with
involution, it suffices to show that $\cW^{\alpha}\diamond
\cW^{\beta}\subset\cW^{\alpha+\beta}$ for $\alpha,\beta
\in\Lm$. \ This follows from $h(\cW^{\lambda},\cW^{\alpha
}\diamond\cW^{\beta})\subset\mathcal{E}^{\lambda+\alpha+\beta}$, which
holds since $N$ is graded, and from $\cW^{\mu}=\{x\in\cW\mid
h(\cW^{\lambda},x)\in\mathcal{E}^{\lambda+\mu}$ for all $\lambda
\in\Lm\}$, which holds since $h$ is nondegenerate and  graded.

Now suppose that
$\Lm_-$ is a subgroup of $\Lm$ satisfying (i)--(iv), and let $\Gm = \Gm(\cE) = \supp(\cZ)$.
By (i)--(ii), $\cA$ is finely graded and $\Lm = \langle\supp(\cA)\rangle$.

To show that $\cA$ is a structurable torus,  it remains to show that each homogenous
element of $\cA$ is invertible.  To do this, it suffices, by Remark \ref{rem:invertible}(i),
to show that for $0\neq x\in\cA^{\al}$
there is $y\in\cA^{-\al}$ with $xy=1$ and $[L_{x},L_{y}]=0$.
If $\al\in\Lm_{-}$, we can take $y=x^{-1}$ in
$\cE$ and use the fact that $\cW$ is a left $\cE$-module.
So we may suppose that $\al\notin\Lm_{-}$.
Then $h(x,\cW^{\beta})\neq0$ for some
$\beta\in-\al+\Lm_{-}$, since $h$ is nondegenerate.  However,
$\cE^{\gm}\cW^{\beta}=\cW^{\gm+\beta}$, so
$h(x,\cW^{\beta})\neq0$ for all $\beta\in-\al+\Lm_{-}$. In
particular,
\[h(x,y)=1\]
for some $y\in\cW^{-\al}$.  Since $0\notin\supp(\cW)$,
we also have $x\diamond\cW^{-\al}=0$, so
\[xy=1.\]
If
$2\al\in\Gm$, then $x=zy$ for some $z\in\cE^{2\al}$, so
$[L_{x},L_{y}]=L_z[L_{y},L_{y}]=0$.
Thus, we can assume that $2\al\notin\Gm$.
Now, for $e\in\cE$, we have
$x(ye)=h(x,e^{*}y)+x\diamond(e^{*}y)=e$
and similarly $y(xe)=e$, so $[L_{x},L_{y}]\cE=0$.  Hence, it remains to show that
\[[L_{x},L_{y}]w=0\]
for $w\in\cW$. Since $2\al\notin\Gm$,
we have $x^{\natural}\neq 0$ (by (iv)) and similarly $y^{\natural}\neq0$.
Now $x\diamond y^{\natural}\in\cW^{-\al}$ so
\[x\diamond y^{\natural}=ay\]
with
$a\in F$. Moreover, $a= a^* = h(x,ay)=h(x,x\diamond y^{\natural})
= h(y^\natural,x\diamond x) =
2h(y^{\natural},x^{\natural})\neq0$.
By (ADJ2), we  have
\[
(w\diamond y^{\natural})\diamond x^{\natural} + (w\diamond x)\diamond(y^{\natural}\diamond x)
=h(y^{\natural
},x^{\natural})w + h(w,y^{\natural}\diamond x)x+h(w,x^{\natural})y^{\natural},
\]
so
\[
(w\diamond y^{\natural})\diamond x^{\natural} + a(w\diamond x)\diamond y
=aw + ah(w,y)x+h(w,x^{\natural})y^{\natural}.
\]
Thus
\begin{equation}
\label{eq:dp}
y \diamond(x\diamond w) -h(w,y)x
=w +a^{-1}h(w,x^{\natural})y^{\natural}-a^{-1}(w\diamond y^{\natural})\diamond x^{\natural}.
\end{equation}
Since $x^{\natural}=zy^{\natural}$ for some $z\in\cZ(\cE)^{4\al}$ (by (iii)), we
see that the right side of (\ref{eq:dp}) is symmetric in $x$ and $y$.  Thus,
the left side of (\ref{eq:dp}) or equivalently $x\diamond(y\diamond
w)+h(w,y)x$ is symmetric in $x$ and $y$.  Since
\[
x(yw)=h(w,y)x+h(x,y\diamond w)+x\diamond(y\diamond w),
\]
we see that $[L_{x},L_{y}]w=0$. So $\cA$ is a structurable torus.

Finally, suppose that (v) also holds.  Then we may choose  $\sg_0\in S_-(\cE)$,
in which case $2\sg_0\in \Gm$, $S_+(\cE) = \Gm$, and $S_-(\cE) = \sg_0 + \Gm$.  Thus, $S_-(\cE)$
generates $\Lm_-$, so our notation
$\Lm_-$, $\cE$, $\cW$ agrees with the notation
from \S\ref{sec:classes}.  Hence, since the map $\diamond : \cW \times \cW \to \cW$ is nonzero, it
follows that $\cA$ is of class III.  Also, the centre of $\cA$ is
$\cZ$ by Lemma \ref{lem:centre}.
\end{proof}

We will also need the following simple lemma.

\begin{lemma}
\label{lem:psi1}  Let $\mathcal{E}$ be a commutative, associative algebra
with involution $*$ and let $\cW$ be a left $\mathcal{E}$-module with
cubic form $N$, hermitian form $h$, and symmetric bilinear form $T$. If
$\psi:\cW\rightarrow\cW$ is a $*$-semilinear map satisfying $N(\psi(x)) = N(x)^*$ and
$h(x,y)=T(x,\psi(y))$, then $N$ is nondegenerate and $(h,N)$ satisfies the adjoint identity if and only
if $T$ is nondegenerate and $(T,N)$ satisfies the adjoint identity. In that case the adjoints
($\natural$ for $(h,N)$ and $\#$ for $(T,N)$) are related by
\[y^{\natural}=\psi^{-1}(y^{\#})=\psi(y)^{\#}.\]
\end{lemma}

\begin{proof}
Clearly $h$ is nondegenerate and $y^{\natural}$ is an adjoint for $(h,N)$
if and only if $T$ is nondegenerate and $y^{\#}=\psi(y^{\natural})$ is an
adjoint for $(T,N)$. Also if those conditions hold, then,
since $N=\ast\circ N\circ\psi$ and $\psi^{\ast
}=\psi$, we have
\begin{align*}
T(x,y^{\#})  & =\partial_{x}(\ast\circ N\circ\psi)|_{y}=(\partial_{\psi
(x)}N|_{\psi(y)})^{\ast}\\
& =(T(\psi(x),\psi(y)^{\#}))^{\ast}=T(x,\psi(\psi(y)^{\#})),
\end{align*}
so $y^{\#}=\psi(\psi(y)^{\#})$, and hence $y^{\natural}=\psi(y)^{\#}$. \ Thus,%
\[
(y^{\natural})^{\natural}=(\psi^{-1}(y^{\#}))^{\natural}=(\psi(\psi
^{-1}(y^{\#})))^{\#}=(y^{\#})^{\#},
\]
proving (a).
\end{proof}

Our goal  now is to construct a pair $(h,N)$ satisfying the assumptions
of Proposition \ref{prop:torusIII}, where $h$ and $N$ are defined
on the space of $3\times 3$-hermitian matrices over an alternative torus with involution $\cC$.
Ultimately we will use two choices for $\cC$, a quaternion torus $\cQ$
and an octonion torus $\cO$.  However, it will only be at the very end of the proof
of the classification theorem (Theorem \ref{thm:IIIc})
that we will see that these are the only choices.
For this reason, we begin with three choices for $\cC$,
a quaternion torus $\cQ$ and \emph{two} octonion tori $\cO$ and $\cO'$.

We now describe these initial choices for $\cC$.
In each case $\cC$ is a finely $\Lm$-graded
algebra with involution $\barinv$, where $\Lm$ is a
finitely generated free abelian group.
Also in each case, we let
\[M = \supp(\cC),\quad \cE = \cZ(\cC) \andd \Lm_- = \supp(\cE).\]

\medskip
\textbf{Case 1:}\
Let $\cC = \cQ$, where
\[\cQ=\CD(F[t_1^{\pm1},t_2^{\pm1},s^{\pm1}],t_1,t_2) \]
with canonical involution $\barinv$ and canonical generators $v_1,v_2$
(see \eqref{eq:CDprocess}).
Let
\[\Lm = \langle\lm_1,\lm_2,\sg\rangle\]
with basis
$\lm_1,\lm_2,\sg$, and give $\cC$ the $\Lm$-grading
so that
\[\text{$v_1,v_2,s$ have degrees  $2\lm_1,2\lm_2,\sg$ respectively.}\]
Then,
\[M = \langle2\lm_1,2\lm_2,\sg\rangle,\quad \cE = F[t_1^{\pm1},t_2^{\pm1},s^{\pm1}] \andd
\Lm_- = \langle 4\lm_1,4\lm_2,\sg\rangle.\]

\medskip
\textbf{Case 2:}\ Let $\cC = \cO$, where
\[\cO=\CD(F[t_1^{\pm1},t_2^{\pm1},s^{\pm1}],t_1,t_2,s)\]
with canonical involution $\barinv$ and canonical generators $v_1,v_2,v$.
Let
\[\Lm = \langle\lm_1,\lm_2,\lm\rangle\]
with basis
$\lm_1,\lm_2,\lm$, and give $\cC$ the $\Lm$-grading
so that
\[\text{$v_1,v_2,v$ have degrees  $2\lm_1,2\lm_2,\lm$ respectively.}\]
Then,
\[M = \langle2\lm_1,2\lm_2,\lm\rangle,\quad
\cE = F[t_1^{\pm1},t_2^{\pm1},s^{\pm1}] \andd \Lm_- = \langle 4\lm_1,4\lm_2,2\lm\rangle.\]

\medskip
\textbf{Case 3:}\ Let $\cC = \cO'$, where
\[\cO'=\CD(F[t_1^{\pm1},t_2^{\pm1},s^{\pm1},r^{\pm1}],t_1,t_2,r)\]
with canonical involution $\barinv$ and canonical generators
$v_1,v_2,v$. Let
\[\Lm = \langle\lm_1,\lm_2,\sigma,\lm\rangle\]
with
basis $\lm_1,\lm_2,\sigma,\lm$, and give $\cC$ the $\Lm$-grading so
that
\[\text{$v_1,v_2,s,v$ have degrees  $2\lm_1,2\lm_2,\sigma,\lm$ respectively.}\]
Then,
\[M = \langle2\lm_1,2\lm_2,\sigma,\lm\rangle,\quad
\cE = F[t_1^{\pm1},t_2^{\pm1},s^{\pm1},r^{\pm1}] \andd\Lm_- = \langle 4\lm_1,4\lm_2,\sigma,2\lm\rangle.\]

\medskip

For the remainder of this section
\emph{we assume that $\cC = \cQ$, $\cO$ or $\cO'$ and we use the notation
introduced in Cases 1, 2 and 3.}

Note that in each case,  $\cC$ is a composition algebra over $\cE$
with norm $n:\cE \to \cE$ defined by $n(a) = a\bar a$ for $a\in \cC$
(see \S \ref{sec:composition}).  We denote the trace on $\cC$
by $t:\cC \to \cE$.  Note also that $\cC$ is an alternative $M$-torus with involution
since, regarding $\cC$ as an $M$-graded algebra with involution, we have
$\cC \ig \cA(2)\otimes \cP(1)$,
$\cA(3)$ or $\cA(3)\otimes \cP(1)$ in Cases 1, 2 and 3 respectively.

In each case we define an involution $*$ on $\cE$
such that
\begin{equation}
\label{eq:def*}
t_i^{*}=t_i,\quad  s^*=-s \quad\text{and, when $\cC = \cO'$,}\quad  r^* = r.
\end{equation}
We then let
\[\cZ = \set{a\in \cE \suchthat a^* = a}\]
be the centre of the algebra with involution $\cE$, and we let
\[\Gm = \supp(\cZ).\]
Then $\Gm = \langle 4\lm_1,4\lm_2,2\sg\rangle$, $\langle 4\lm_1,4\lm_2,4\lm\rangle$
or $\langle 4\lm_1,4\lm_2,2\sigma,2\lm\rangle$ in Cases 1, 2 and 3 respectively.
Thus in all cases we have
\begin{equation}
\label{eq:4Lm}
4\Lm \subset \Gm.
\end{equation}

In some cases, we can introduce an important $*$-semilinear automorphism $\theta$
of~$\cC$:

\begin{lemma}
\label{lem:theta}
Suppose that either $\cC=\cQ$ or there exists $\iota\in F$ with $\iota^{2}=-1$ and
$\cC=\cO$.
Then there exists a unique $*$-semilinear $F$-algebra automorphism $\theta$ of $\cC$
which fixes $v_1$ and $v_2$, and, if $\cC = \cO$, satisfies $\theta(v)=\iota v$.
\end{lemma}

\begin{proof}   Uniqueness is clear.
If $\cC=\cQ = \CD(\cE,t_1,t_2)$, there is clearly
an $F$-algebra automorphism $\theta$ of $\cC$
extending $*$ on $\cE$ and fixing $v_i$. On the other hand, if
$\cC = \cO = \CD(\cE,t_1,t_2,s)$, we can
view $\cO$ as $\CD(\cE,t_1,t_2,-s)$ with canonical generators $v_1,v_2,\iota v$,
so there exists an $F$-algebra automorphism $\theta$ of $\cC$ extending
$*$ on $\cE$ and mapping $v_i \mapsto v_i$ and $v\mapsto \iota v$.
\end{proof}

As in \S \ref{sec:hermitian} (with $K = \cE$), let $\Herm$ be the $\cE$-module
of $3\times 3$-hermitian matrices over $\cE$, and let $T$ and $N$ denote
the standard trace and norm on $\Herm$ respectively.
Let
\[\lambda_{3}=-\lambda_{1}-\lambda_{2},\]
and we use Lemma \ref{lem:gradH}(a) to grade $\Herm$
so that $N$ is a graded cubic form.  Thus,
\[a[ij]\in\Herm\quad\text{has degree}\quad\lambda
_{i}+\lambda_{j}+\mu,\]
for $a\in\cC^{\mu}$, $\mu\in M$, $(i,j,k)\cyc$, and
for $a\in K^{\mu}$, $\mu\in\Lm_-$, $i=j$.

Although $N$ is graded,
$T$ is not
graded and needs to be modified to give a graded hermitian form. We do this in
the next lemma using a $\ast$-semilinear vector space isomorphism
$\psi:\Herm\rightarrow\Herm$.
Such an isomorphism is said to be \textit{semi-norm preserving} if
$N(\psi(x))=N(x)^{\ast}$ for all $x\in\Herm$.

In
general, if $\cU$ is an $\mathcal{E}$-module with a nondegenerate
symmetric bilinear form $f:\cU\times\cU\rightarrow\mathcal{E}$
and $\psi:\cU\rightarrow\cU$ is an $\alpha$-semilinear vector
space isomorphism where $\alpha$ is an $F$-algebra automorphism of $\mathcal{E}$, we write
$\psi^{\ast}=\varphi$ (relative to $f$) if
$f(\psi(x),y)=\alpha(f(x,\varphi(y))$ for $x,y\in\cU$.
We note that if $\varphi$
exists, it is unique and is an $\alpha^{-1}$-semilinear vector space
isomorphism. \ Moreover, if $\psi_{1}^{\ast}$ and $\psi_{2}^{\ast}$ exist,
then $(\psi_{1}\psi_{2})^{\ast}=\psi_{2}^{\ast}\psi_{1}^{\ast}$. \ Also, if
$\alpha=\ast$, the involution on $\mathcal{E}$, then $h(x,y)=f(x,\psi(y))$ is
a nondegenerate sesquilinear form, which is hermitian if and only if
$\psi^{\ast}=\psi$ (relative to $f$).

A key lemma in our construction and classification of structurable
tori of class III(c) is the following.

\begin{lemma}
\label{lem:psi2}  Suppose that $\psi:\Herm\rightarrow\Herm$
is a $*$-semilinear vector space isomorphism
 and set $h(x,y):=T(x,\psi(y))$ for $x,y\in\Herm$.  Then
$\psi$ is semi-norm preserving and $h$ is a graded
hermitian form if and only if the following conditions hold:
\begin{itemize}
\item[(i)] Either $\cC=\mathcal{Q}$ or there exists $\iota\in F$ with
$\iota^{2}=-1$ and $\cC=\mathcal{O}$,

\item[(ii)] $\psi: \sum_{i} a_{i}[ii]+\sum_{(i,j,k)\cyc}
x_{i}[jk]\rightarrow\sum_{i} a_{i}^{*}n(w_{i})[ii]+ \sum
_{(i,j,k)\cyc} (w_{j}\theta(x_{i}))\bar{w}_{k}[jk],$ where
$w_{i}\in\cC^{2\lambda_{i}}$ with $w_{1}w_{2}w_{3}=1$ and $\theta$ is
the $*$-semilinear $F$-algebra automorphism of $\cC$ which fixes
$v_{1}$ and $v_{2}$, and satisfies $\theta(v)=\iota v$ in the case
$\cC=\mathcal{O}$.
\end{itemize}
\end{lemma}

\begin{proof}
Before beginning the proof, we need some preliminaries about triality
triples. (See for example \cite{J1} for
a slightly more general concept in the finite dimensional case.)

Let $\alpha$ be an automorphism of $\mathcal{E}$. \ A \emph{triality triple}
of $\alpha$-semisimilarities of $\cC$ is a triple
$(\varphi_{1},\varphi_{2},\varphi_{3})$ of $\alpha$-semilinear bijections from
$\cC$ to $\cC$ such that for some invertible $c_{i}%
\in\mathcal{E}$ we have
\begin{align}
n(\varphi_{i}(x_{i}))  & =c_{i}\alpha(n(x_{i})),\label{eq:phi1}\\
t(\varphi_{1}(x_{1})\varphi_{2}(x_{2})\varphi_{3}(x_{3}))  & =\alpha
(t(x_{1}x_{2}x_{3}))\label{eq:phi2}%
\end{align}
for $x_{i}\in\cC$. In that case, the elements $c_{i}$ are called the
\textit{multipliers} of the triple.  If $\alpha=id$, we say that
$\varphi_{i}$ is a  \emph{similarity}. Since $t(abc)$ is invariant under cyclic
permutations, any cyclic permutation of a triality triple is a triality
triple.

If $\varphi:\Herm\rightarrow\Herm$
is an $\alpha$-semilinear bijection which stabilizes the
spaces $\mathcal{E}[ii]$, $\cC[jk]$ for $(i,j,k)\cyc$, we can
write
\begin{equation}
\varphi:\sum_{i}a_{i}[ii]+\sum_{(i,j,k)\cyc}x_{i}[jk]\rightarrow
\sum_{i}\alpha(a_{i})c_{i}^{-1}[ii]+\sum_{(i,j,k)\cyc}\varphi
_{i}(x_{i})[jk]\label{eq:phi3}%
\end{equation}
for some $\alpha$-semilinear bijections $\varphi_{i}:\cC%
\rightarrow\cC$ and some $0\neq c_{i}\in\mathcal{E}$.
We claim that $N(\varphi(x)) = \alpha(N(x))$ for $x\in \Herm$ if and only if $(\varphi_{1},\varphi
_{2},\varphi_{3})$ is a triality triple with multipliers $c_{i}$. \ Indeed,
(\ref{eq:phi1}) and (\ref{eq:phi2}) follow from $N(\varphi(x))=\alpha(N(x))$
for $x=1[ii]+x_{i}[jk]$ and $x={\sum_{(i,j,k)\cyc}}
x_{i}[jk]$, respectively.  The converse is immediate.

Triality triples
of similarities can be created using invertible elements of $\cC$.
Indeed, it is obvious that $(L_{a},L_{b},L_{c})$ is a triality triple of similarities
if $a,b,c\in\mathcal{E}$ with $abc=1$ (see \S \ref{sec:composition} for this notation). Also, if
$x\in\cC$ is invertible, then the right (or left) Moufang identity
and the associativity of $t(abc)$ shows that
$t((ax)(x^{-1}bx^{-1})(xc))=t(abc)$,
so $(R_{x},L_{x^{-1}}R_{x^{-1}},L_{x})$ is a triality triple.
Moreover, if $xyz=1$ in $\cC$, then
$yz=n(x)^{-1}\bar{x}$ and we can compute
\begin{eqnarray*}
\bar{z}^{-1}(by)\bar{z}^{-1} &=&n(z)^{-2}(zb)(yz)=n(x)^{-1}n(z)^{-2}(zb)\bar{%
x}, \\
\bar{z}(y^{-1}cy^{-1}) &=&n(y^{-2})((\bar{z}\bar{y})c)\bar{y}%
=n(x)^{-1}n(y)^{-2}(xc)\bar{y}.
\end{eqnarray*}%
Thus,
\begin{multline}
(\id,L_{n(x)n(z)^{2}},L_{n(x)n(y)^{2}})(R_{\bar{z}},L_{\bar{z}^{-1}}
R_{\bar{z}^{-1}},L_{\bar{z}})(L_{y},R_{y},L_{y^{-1}}R_{y^{-1}})\\
=(R_{\bar{z}}L_{y},R_{\bar{x}}L_{z},R_{\bar{y}}L_{x})
\end{multline}
is a triality triple of similarities.

To begin  the proof, we recall that
$T$ is not graded using the given grading on $\Herm$. But we obtain a second
grading on $\Herm$ by letting
\[
a[ij]\in\Herm\quad\text{have degree}\quad-\lambda
_{i}-\lambda_{j}+\mu,
\]
for $a\in\cC^{\mu}$, $\mu\in M$, $(i,j,k)\cyc$, and
for $a\in K^{\mu}$, $\mu\in\Lm_-$, $i=j$.
Denoting $\Herm$ with this grading
by $\Herm^{\prime}$, we see that $T:\Herm\times\Herm^{\prime}\rightarrow
\mathcal{E}$ is graded. Since $T$ is nondegenerate, we see that $h$ is graded
if and only if $\psi:\Herm\rightarrow\cH(\cC_{3})^{\prime}$ is graded; that is
\[
\psi:\Herm^{\lambda_{i}+\lambda_{j}+\mu}\rightarrow
\Herm^{\prime\lambda_{i}+\lambda_{j}+\mu}
=\Herm^{\prime-\lambda_{i}-\lambda_{j}+\mu'
}=\Herm^{\lambda_{i}+\lambda_{j}+\mu'}
\]
for $\mu\in M$, $(i,j,k)\cyc$, and
for $\mu\in\Lm_-$, $i=j$, where $\mu'=2(\lambda_{i}+\lambda_{j})+\mu$. In
particular, this implies that $\psi$ stabilizes $\mathcal{E}[ii]$ and
$\cC[ij]$ for $(i,j,k)\cyc$.   Thus, $\psi$ is semi-norm
preserving and $h$ is a graded sesquilinear form if and only if $\psi$ has the
form (\ref{eq:phi3}) for a triality triple $(\psi_{1},\psi_{2},\psi_{3})$ of
$\ast$-semisimilarities with $\psi_{i}(\cC^{\mu})=\cC%
^{-2\lambda_{i}+\mu}$ and multipliers $c_{i}\in\mathcal{E}^{-4\lambda_{i}}$.
If these conditions hold, let $u_{i}=\psi_{i}(1)\in\cC%
^{-2\lambda_{i}}$. We see that $\bar{u}_{i}=-u_{i}$ and $n(u_{i})=c_{i}$ by
(\ref{eq:phi1}).  Also, since $(u_{1}u_{2})u_{3}\in\cC^{0}=F1$,
(\ref{eq:phi2}) shows that $u_{1}u_{2}u_{3}=1$. \ Let $w_{i}=\bar{u}_{i}%
^{-1}=-u_{i}^{-1}\in\cC^{2\lambda_{i}}$ so $w_{1}w_{2}w_{3}=1$ and
$\bar{w}_{i}=-w_{i}$.  We know that $\eta_{i}=-R_{w_{k}}L_{w_{j}}$ for
$(i,j,k)\cyc$  defines a
triality triple $(\eta_1,\eta_2,\eta_3)$ of similarities with multipliers $n(w_{j})n(w_{k}%
)=n(w_{i}^{-1})=c_{i}$ and hence a norm preserving map $\eta$. \ We see that
$\eta_{i}(1)=-w_{j}w_{k}=-w_{i}^{-1}=u_{i}$. \ Thus, $\eta^{-1}\psi$ is
semi-norm preserving and given by the triality triple $(\theta_{1},\theta
_{2},\theta_{3})$ where $\theta_{i}=\eta_{i}^{-1}\psi_{i}$. $\ $%
Clearly, $\theta_{i}(1)=1$, so $\theta_{i}$ has multiplier $1$.  Thus,
$t(\theta_{i}(x))=t(x)$ and $\theta_{i}$ commutes with $-$.  Now
\begin{align*}
n(ab,c)^{\ast} &  =t(ab\bar{c})^{\ast}=t(\theta_{1}(a)\theta_{1}(b)\theta
_{3}(\bar{c}))=t(\theta_{1}(a)\theta_{1}(b)\overline{\theta_{3}(c)})\\
&  =n(\theta_{1}(a)\theta_{2}(b),\theta_{3}(c))=n(\theta_{3}^{-1}(\theta
_{1}(a)\theta_{2}(b)),c)^{\ast},
\end{align*}
so $\theta_{3}(ab)=\theta_{1}(a)\theta_{2}(b)$. \ Since $\theta_{i}(1)=1$, we
see that $\theta_{1}=\theta_{2}=\theta_{3}$ is a $*$-semilinear automorphism
$\theta$ of $\cC$. \ Since $\eta_{i}(\cC^{\mu})=\cC%
^{-2\lambda_{i}+\mu}=\psi_{i}(\cC^{\mu})$, we see that $\theta$ is
graded. We have shown that if $\psi$ is semi-norm preserving and $h$ is a
graded sesquilinear form, then $\psi$ is of the form (ii) for some graded
$*$-semilinear automorphism $\theta$
and some $w_i\in \cC^{2\lm_i}$ with $w_1w_2w_3 = 1$.
Conversely, if $\psi$ is of the form (ii)
for such a choice of $\theta$ and $w_i$, then $(\psi_{1},\psi_{2},\psi_{3})$ is a triality triple
of $\ast$-semisimilarities with $\psi_{i}(\cC^{\mu})=\cC%
^{-2\lambda_{i}+\mu}$ and multipliers $c_{i}\in\mathcal{E}^{-4\lambda_{i}}$,
so $\psi$ is semi-norm preserving and $h$ is a graded sesquilinear form.

So  we assume for the rest of the proof that $\psi$ has the form (ii) for some graded
$*$-semilinear automorphism $\theta$
and some $w_i\in \cC^{2\lm_i}$ with $w_1w_2w_3 = 1$.
We then have
\[\psi_i = \eta_i\theta,\quad \text{where } \eta_{i}=-R_{w_{k}}L_{w_{j}}\]
for
$(i,j,k)\cyc$, and we set $u_{i}=\psi_{i}(1)= -w_jw_k = -w_i^{-1}$
as above. It remains to show that under these conditions $\psi^{\ast}=\psi$ if and only
if (i) holds and $\theta$ is as in (ii).

First note, that $n(w_{i})\in\mathcal{E}^{4\lambda_{i}}=Ft_{i}$ for $i=1,2$, and
$n(w_{3})\in\mathcal{E}^{4\lambda_{3}}=F(t_{1}t_{2})^{-1}$.  Thus, $n(w_{i})^{\ast}=n(w_{i})$,
so $\psi^{\ast}=\psi$ (relative to $T$)
if and only if $\psi_{i}^{\ast}=\psi_{i}$ (relative to $n$)
for $i=1,2,3$.  Now $n(\theta(x),y)=n(x,\theta^{-1}(y))^{\ast}$ so
$\theta^{\ast}=\theta^{-1}$.   Also, $L_{w}^{\ast}=L_{\bar{w}}$ and
$R_{w}^{\ast}=R_{\bar{w}}$, so
$\eta_{i}^{\ast}=-L_{w_{j}}R_{w_{k}}$.
Hence
\begin{align*}
\psi_i^* = \psi_i &\iff \theta^* \eta_i^* = \eta_i \theta \iff \theta^{-1} \eta_i^* = \eta_i \theta
\\
&\iff \theta^2 = \theta\eta_i^{-1}\theta^{-1}\eta_i^*
\iff \theta^2 = L_{\theta(w_j)}^{-1}R_{\theta(w_k)}^{-1}L_{w_j}R_{w_k}.
\end{align*}
Now
$\eta_{i}^{\ast}(1)=-w_{j}w_{k}=u_{i}=\eta_{i}(1)$.  So if $\psi_{i}^{\ast}=\psi_{i}$,
then $\theta^{-1} \eta_i^*(1) = \eta_i \theta(1)$, which implies that $\theta^{-1}(u_i) = u_i$
and hence $\theta(w_i) = w_i$.
Thus setting
\[\beta_i = L_{w_j^{-1}}R_{w_k^{-1}}L_{w_j}R_{w_k}\]
for $(i,j,k)\cyc$, we see that $\psi^* = \psi$ if and only if
$\theta(w_i) = w_i$ and $\theta^2 = \beta_i$ for  $i=1,2,3$.

If $\cC=Q$ and $(i,j,k)\cyc$, then
$\beta_i(x)=w_{j}^{-1}w_jxw_{k}w_{k}^{-1}=x$ for $x\in \cC$. If
$\cC=\mathcal{O}$ or $\mathcal{O}^{\prime}$, we have
$\cC=\CD(\mathcal{E},t_{1},t_{2})\oplus v\CD(\mathcal{E},t_{1},t_{2})$,
with multiplication given in \eqref{eq:CDdef}
(with $v^{2}=s$ or $r$ respectively).
Moreover, each $w_i$ is in $\CD(\mathcal{E},t_{1},t_{2})$
since $w_{i}\in\cC^{2\lambda_{i}}=Fv_{i}$, and hence,
for $x\in\CD(\mathcal{E},t_{1},t_{2})$, we have
$\beta_i(x)=x$ and
\[
\beta_i(vx)= v(w_j^{-1}w_k^{-1}w_jw_kx) = - v(w_j^{-1}w_k^{-1}w_kw_jx) = -vx.
\]

Thus, $\beta_i=\id$ if $\mathcal{C=Q}$, whereas $\beta_i$ is the automorphism with
$\beta_i(v_1)=v_1$, $\beta_i(v_2)=v_2$ and $\beta_i(v)=-v$ if $\cC=\mathcal{O}$ or
$\mathcal{O}^{\prime}$.  In particular,
$\beta := \beta_i$ is independent of $i$. Since $Fw_{i}=Fv_{i}$, we have shown that
$\psi^{\ast}=\psi$ if and only if $\theta(v_{i})=v_{i}$, for $i=1,2$, and
$\theta^{2}=\beta$. Since $\theta$ is graded so $\theta(v)\in Fv$, this is
equivalent to $\theta(v_{i})=v_{i}$ and $\theta(v)=\iota v$ with $\iota
^{2}=-1$ if $\cC=\mathcal{O}$ or $\mathcal{O}^{\prime}$. However, if
$\cC=\mathcal{O}^{\prime}$, then $v^{2}=r\in\mathcal{E}$, so $r^{\ast
}=\theta(v^{2})=-r$, a   contradiction.
\end{proof}

We can now use Proposition \ref{prop:torusIII} and Lemma \ref{lem:psi2}
to give a construction of structurable tori of class
III(c).  In this construction we make a choice of $w_i\in \cC^{2\lm_i}$ with
$w_1w_2w_3 = 1$.  We note that there always exists such a choice since we may
take $w_i = v_i$, where $v_3 = (v_1v_2)^{-1}$.

\begin{theorem}
\label{thm:constructIIIc}
Let $\cE = F[t_1^{\pm1},t_2^{\pm1},s^{\pm1}]$ with the involution $*$
such that $t_i^* = t_i$ and $s^* = -s$, and let $\cZ$
be the centre of this algebra with involution.
Suppose that either
\begin{itemize}
\item[(1)] $\cC = \cQ = \CD(\cE,t_1,t_2)$ with canonical generators $v_1,v_2$ and
canonical involution $\barinv$, and
$\cC$ is $\Lm$-graded with $v_1,v_2,s$ of degrees $2\lm_1,2\lm_2,\sg$ respectively,
where $\Lm$ has basis $\lm_1,\lm_2,\sg$;
or
\item[(2)] There exists $\iota\in F$ with $\iota^{2}=-1$,
$\cC = \cO = \CD(\cE,t_1,t_2,s)$ with canonical generators $v_1,v_2,v$
and canonical involution $\barinv$, and $\cC$ is
$\Lm$-graded with $v_1,v_2,v$ of degrees $2\lm_1,2\lm_2,\lm$ respectively,
where $\Lm$ has basis $\lm_1,\lm_2,\lm$.
\end{itemize}
Let
$\Herm$ be the $\cE$-module of hermitian $3\times 3$-matrices over $\cC$,
set $\lm_3 = -\lm_1-\lm_2$,  and
grade $\Herm$ by $\Lm$ with $a[ij]$ of degree
$\lm_i +\lm_j+\mu$ if $a\in\cC^{\mu}$.
Choose $w_{i}\in\cC^{2\lm_{i}}$ with $w_{1}w_{2}w_{3}=1$
and define $\psi: \Herm \to \Herm$ by
\[
\psi \Big(\sum_i a_{i}[ii]+\sum_{(i,j,k)\cyc} x_{i}[jk]\Big)
=
\sum_i
a_{i}^{*}n(w_{i})[ii]+
\sum_{(i,j,k)\cyc}
(w_{j}\theta(x_{i}))\bar{w}_{k}[jk],\]
where $\theta$ is the $*$-semilinear $F$-algebra automorphism of $\cC$
which fixes $v_1$ and $v_2$, and satisfies $\theta(v)=\iota v$ in the case
$\cC=\cO$.  Let $T$ and $N$ be the standard trace and norm forms on $\Herm$
and define $h : \Herm \times \Herm \to \cE$ by
\[h(x,y)=T(x,\psi(y)).\]
Let $\cA(h,N) = \cE \oplus \Herm$ with product and involution given
(as in \eqref{eq:defAhN})
by
\[(a,x)(b,y) = (ab + h(x,y), ay + b^*x + x\diamond y)
\andd (a,u)^* = (a^*,u),\] where $x\diamond y =
\psi(x)\times\psi(y)$ and $\times$ is defined by
\eqref{eq:adjointlin}. Give $\cA(h,N)$ the $\Lm$-grading extending
the gradings on $\cE$ and $\Herm$. Then $\cA(h,N)$ is a
structurable $\Lm$-torus of class III(c) with centre $\cZ$.
Moreover, up to isograded isomorphism, $\cA(h,N)$ depends only on
$\cC$ and not on the choice of $w_{i}$ or $\iota$.
\end{theorem}

\begin{proof} Let $\cA = \cA(h,N) = \cE \oplus \cW$,
where $\cW = \Herm$.
By Lemmas \ref{lem:psi1} and \ref{lem:psi2},
all of the assumptions in the first sentence
of Proposition \ref{prop:torusIII} are satisfied.  So
$\cA$ is a $\Lm$-graded structurable algebra.  To show
that $\cA$ is a structurable torus of class III, we check conditions (i)--(v)
of Proposition \ref{prop:torusIII}.
We note first that $\cW$
is finely graded by $\Lm$ and the support of $\cW$ is disjoint from $\Lm_-$.
Indeed, $\cC$ is finely graded by
$M$, so the three spaces $\cC[ij]$ are finely graded by the three
cosets $-\lm_{k}+M$ where $(i,j,k)\cyc$, and the three spaces
$\cE[ii]$ are finely graded by the three cosets
$2\lm_{i}+\Lm_{-}$ which lie in $M$.
Also, since $\lm_i\in \supp(\cW)$ and (if $\cC = \cO$)
$\lm_i + \lm \in \supp(\cW)$, we have
$\Lm = \langle \Lm_-, \supp(\cW)\rangle$.
Next, we have seen in \eqref{eq:4Lm}
that $4\Lm\subset\Gm$. Also, if $0\neq x\in\cW^{\al}$
with $2\al\notin\Gm$, then $\al\in-\lm_{k}+M$ for some $k$ and
$x=a[ij]$ where $(i,j,k)\cyc$. Thus,
$x^{\natural}=\psi^{-1}(x^{\#})=-n(a)^{*}n(w_{k})^{-1}%
[kk]\neq0\text{.}
$
The other conditions in (i)--(v) of Proposition \ref{prop:torusIII}
clearly hold, so $\cA$ is a
structurable torus of class III with centre $\cZ$.

To see that $\cA$ is of class III(c), let $\bLm = \Lm/\Lm_-$ and
$\cI$ be as in \S\ref{sec:geometry}. The points in $\cI$ of order 2
are $2\blm_i$, $1\le i \le 3$, and the points in $\cI$ of order
4 are $-\blm_i + \bmu$, $1\le i \le 3$, $\bmu\in \bar M$.
In particular, $\cI$ has a point of order 4, so $\cA$ is not of class III(b).
Assume for contradiction that $\cI$ is of class III(a).  Then, by Corollary
\ref{cor:geometry}(a), $\cI$ is a star with centre $2\blm_i$ for some $i$.
But if $(i,j,k)\cyc$,  then $2\blm_i$ is not collinear with
$-\blm_k$ since $1[ii]\in \cW^{2\lm_i}$, $1[ij]\in \cW^{-\lm_k}$
and $(1(ii])\diamond (1[ij]) \in
\cE[ii]\times \cC[ij] = 0$ (by \eqref{eq:adjointlin}).
This contradiction shows that $\cA$ is of class III(c).

For the last statement, suppose we use the $v_{i}$ and (if $\cC=\cO$) a fixed choice of $\iota$  to
construct $\theta$, $\psi$, $h$ and $\cA = \cA(h,N)$.
Suppose that $\theta'$, $\psi^{\prime}$, $h'$ and
$\cA' = \cA(h',N)$
are constructed from some
$w_{i}$  and (if $\cC=\cO$) the same $\iota$.  Now there is a graded
automorphism $\ph$ of $\cC$ with $\ph(v_{i})=w_{i}$ and (if $\cC=\cO$) $\ph(v)=v$.
Also $\theta'=\theta$, so $\ph\theta = \theta'\ph$. Hence,
$\ph R_{\bar v_k}L_{v_j}\theta=R_{\bar w_k}L_{w_j}\theta'\ph$.
Also, $n(\ph(a))=\ph(n(a)) $ for $a\in\cO$.  Thus, if we define
$\ph$ on $x\in\Herm$ by letting $\ph$ act on the entries of $x$, we have
$\ph\circ\psi=\psi^{\prime}\circ\ph$.
So $N(\ph(x))=\ph(N(x))$, $T(\ph(x),\ph(y))=\ph(T(x,y))$ and
$h^{\prime}(\ph(x),\ph(y))=\ph(h(x,y))$.  Hence,
the map $(a,x) \mapsto (\ph(a),\ph(x))$ is a graded isomorphism of
$\cA$ onto $\cA'$.

Now suppose $\cC=\cO$, suppose
$\theta$, $\psi$, $h$ and $\cA$ are constructed from
$v_i$ and $\iota$ as above,
and suppose $\theta'$, $\psi'$, $h'$ and $\cA'$ are constructed
from $v_i$ and $-\iota$.
Let $\ph$ be the isograded automorphism of
$\cO$ with $\ph(v_{i})=v_{i}$ and $\ph(v)=v^{-1}$ associated with
the automorphism of $\Lm$ fixing $\lm_{i}$ and with $\lm\rightarrow-\lm$.
Now $\theta^{\prime}(v)=-\iota v$ implies
$\theta^{\prime}(v^{-1})=\iota v^{-1}$, so
$\ph\theta
=\theta'\ph$.
Hence,
$\ph R_{\bar v_k}L_{v_j}\theta=R_{\bar v_k}L_{v_j}\theta'\ph$.
The other calculations are the same as before and
the map $(a,x) \mapsto (\ph(a),\ph(x))$ is an isograded isomorphism of
$\cA$ onto $\cA'$.
\end{proof}

\begin{definition}
\label{def:constructIIIc}
Suppose that we have the assumptions of Theorem
\ref{thm:constructIIIc}.  Let
\[\cA(\Herm):=\cA(h,N)\]
be the structurable torus of class III(c) constructed in the theorem
using the choice $w_i = v_i$ for $1\le i \le 3$, and (if $\cC = \cO$)
a fixed choice of $\iota$. We call $\cA(\Herm)$ the
\textit{torus based on}  $\Herm$.
\end{definition}

We have thus constructed two
structurable tori of class III(c): $\cA(\cH(\cQ_3))$
and, if $-1\in (F^\times)^2$,  $\cA(\cH(\cO_3))$.
Note that under the inclusion map we  have
\[\cA(\cH(\cQ_3)) \ig \cA(\cH(\cO_3))^{\langle \lm_1,\lm_2,2\lm\rangle}.\]

\section{Classification of tori of class III$(c)$}
\label{sec:classificationIIIc}

In this final section, we classify structurable tori of class III(c).

\begin{lemma}
\label{lem:basis} Let $\Lm_{1}$ be a subgroup of a free abelian group
$\Lm$ of rank $n$ and let $\Lm_{2}$ be a subgroup of $\Lm_{1}$
with $[\Lm_{1}:\Lm_{2}]=2$ and $4\Lm\subset\Lm_{2}$.
If $v_1,\dots,v_r\in \Lm/\Lm_1$
and $2v_{1},\ldots,2v_{r}$ is a basis for the $\mathbb{Z}_{2}$-vector space
$2(\Lm/\Lm_{1})$, then there is a basis $\lm_{1},\ldots
,\lm_{n}$ for $\Lm$ such that
\begin{itemize}
\item[(i)] $v_{i}=\lm_{i}+\Lm_{1}$, $1\leq i\leq r$,
\item[(ii)] $\Lm_{1}=\left\langle 4\lm_{1},\ldots
,4\lm_{r},k_{r+1}\lm_{r+1},\ldots,k_{n}\lm_{n}\right\rangle$,
\item[(iii)] $\Lm_{2}=\left\langle 4\lm_{1},\ldots,
4\lm_{r},2k_{r+1}\lm_{r+1},k_{r+2}\lm_{r+2}\ldots,k_{n}\lm_{n}\right\rangle$
\end{itemize}
with $k_{i}=1$ or $2$ for $r+1 \le i\leq n$.
\end{lemma}

\begin{proof}
Since $4\Lm\subset\Lm_{1}$, it follows from the fundamental theorem for finitely generated
abelian groups that there is a basis $\lm_{1},\ldots,\lm_{n}$ for $\Lm$ with
\[\Lm_{1}=\left\langle k_{1}\lm_{1},\ldots,k_{n}\lm_{n}\right\rangle\]
where $k_i=1$ , $2$, or $4$.
Assume that $l$ is an integer such that $1\le l \le r$
and $v_i=\lm_i+\Lm_1$ for $i<l$.  Since
$4\lm_{i}\in\Lm_{1}$, we can write $v_{l}=\sum_{i=1}^{n}m_{i}%
\lm_{i}+\Lm_{1}$ with $-1\leq m_{i}\leq2$. \ Now
\[
2v_{l}=\sum_{i=1}^{l-1}m_{i}2v_{i}+(\sum_{i=l}^{n}2m_{i}\lm_{i}%
+\Lm_{1})
\]
so $\sum_{i=l}^{n}2m_{i}\lm_{i}\notin\Lm_{1}$. \ We can rearrange
$\lm_{l},\ldots,\lm_{n}$ to assume that $2m_{l}\lm_{l}%
\notin\Lm_{1}$; i.e., $m_{l}=\pm1$. \ Thus, $2\lm_{l}\notin\Lm_{1}
$ and $k_{l}=4$. \ Replacing $\lm_{l}$ by $\lm_{l}^{\prime}=\sum
_{i=1}^{n}m_{i}\lm_{i}$ gives a basis $B^{\prime}$ for $\Lm$ with
$v_{l}=\lm_{l}^{\prime}+\Lm_{1}$. \ Since $4\lm_{l}\equiv
\pm4\lm_{l}^{\prime}$ modulo the subgroup generated by $k_{i}\lm_{i}$
with $i\neq l$, we still have $\Lm_{1}=\left\langle k_{1}\lm
_{1},\ldots,k_{l}\lm_{l}^{\prime},\ldots,k_{n}\lm_{n}\right\rangle $
using the basis $B^{\prime}$. \ By induction on $l$, we see that there is a
basis for $\Lm$ satisfying (i) and $\Lm_{1}=\left\langle k_{1}%
\lm_{1},\ldots,k_{n}\lm_{n}\right\rangle $. \ Since $2v_{1}%
,\ldots,2v_{r}$ is a basis for $2(\Lm/\Lm_{1})$, we have also have (ii).

Clearly, $k_{i}\lm_{i}=4\lm_{i}\in\Lm_{2}$ if $i\leq r$. \ Since
some $k_{i}\lm_{i}\notin\Lm_{2}$, we can rearrange $\lm
_{r+1},\ldots,\lm_{n}$ so that $k_{r+1}\lm_{r+1}\notin\Lm_{2}$ and
$k_{r+1}\geq k_{j}$ if $k_{j}\lm_{j}\notin\Lm_{2}$. \ Suppose
$k_{j}\lm_{j}\notin\Lm_{2}$ with $j>r+1$. \ Set $\lm_{j}^{\prime
}=k_{r+1}\lm_{r+1}+\lm_{j}$ if $k_{j}=1$, and set $\lm_{j}%
^{\prime}=\lm_{r+1}+\lm_{j}$ if $k_{j}=2$ (so $k_{r+1}=2$). \ In
either case, replacing $\lm_{j}$ by $\lm_{j}^{\prime}$ gives a basis
$B^{\prime}$ for $\Lm$ and $k_{j}\lm_{j}^{\prime}=k_{r+1}\lm
_{r+1}+k_{j}\lm_{j}$. \ Moreover, $k_{j}\lm_{j}^{\prime}\in
(\Lm_{1}\backslash\Lm_{2})+(\Lm_{1}\backslash\Lm_{2}%
)\subset\Lm_{2}$, since$[\Lm_{1}:\Lm_{2}]=2$. \ We see that we can
choose a basis for $\Lm$ satisfying (i) and (ii) with $k_{i}\lm_{i}%
\in\Lm_{2}$ for all $i\neq r+1$. \ Now $\Lm_{2}^{\prime}:=\left\langle
4\lm_{1},\ldots,4\lm_{r},2k_{r+1}\lm_{r+1},k_{r+2}\lm
_{r+2}\ldots,k_{n}\lm_{n}\right\rangle \subset\Lm_{2}$ has
$[\Lm_{1}:\Lm_{2}^{\prime}]=2$, so $\Lm_{2}=\Lm_{2}^{\prime}$,
showing (iii).
\end{proof}

Assume for the rest of the section that $\cA$ is a structurable
$\Lm$-torus of class III(c).  We use the notation of \S\ref{sec:classes}
and \S\ref{sec:geometry}.  In particular, $\cI$ denotes the incidence geometry
associated with $\cA$.

\begin{lemma}
\label{lem:IIIcgeometry}
If $\cA$ is a structurable $\Lm$-torus
of class III(c), then there is a basis
$\lm_{1},\lm_{2},\eta_{0},\eta_{1},\ldots,\eta_{n}$ of $\Lm$ with $n\ge 0$ such that
\begin{itemize}
\item[(i)] $\Lm_{-}=\left\langle 4\lm_{1},4\lm_{2},k_{0}\eta_{0},
k_{1}\eta_{1},\ldots,k_{n}\eta_{n}\right\rangle $ with $k_{i}=1$ or $2$,
\item[(ii)] $\Gm=\left\langle 4\lm_{1},4\lm_{2},2k_{0}\eta_{0},k_{1}\eta_{1},\ldots,k_{n}\eta_{n}\right\rangle$,
\item[(iii)] if we set $\lm_3 = -\lm_1-\lm_2$ and
$M=\left\langle2\lm_{1},2\lm_{2},\eta_{0},\eta_{1},\ldots,\eta_{n}\right\rangle $,
then the points of $\cI$ are $2\blm_{i} $ and all
$\bal_{i}$ with $\al_{i}\in\lm_{i}+M$, and the lines are
$[2\blm_{1},2\blm_{2},2\blm_{3}]$, $[\bar{\al
}_{i},\bal_{i},2\blm_{i}]$ and $[\bal_{1},\bar
{\al}_{2},\bal_{3}]$ with $\al_{1}+\al_{2}+\al_{3}=0$.
\end{itemize}
\end{lemma}

\begin{proof}
Let $\blm_{1}$ be a point of order $4$. By Lemma
\ref{lem:geometry} (b) and (a),  if
$\bde\neq2\blm_{1}$ is any point of order $2$, then
$\bde$ is collinear with some point on $[\blm_{1}%
,\blm_{1},2\blm_{1}]$ which must be $2\blm_{1}$.
Since $\cI$ is not a star, Corollary \ref{cor:geometry}(c) shows that
there is some $\blm_{2}$ of order $4$ which is not collinear with
$2\blm_{1}$. Since $\blm_{2}$ is collinear with
$2\blm_{2}$ by Lemma \ref{lem:geometry}(a), we see that $2\blm_{2}\neq2\blm_{1}$.
Since $\ord{2\blm_{2}} =2$, this
shows that $2\blm_{2}$ is collinear with $2\blm_{1}$. Thus,
$[2\blm_{1},2\blm_{2},2\blm_{3}]$ is a line where
$\lm_{3}= -\lm_1-\lm_2$.

If $\bal$ is a point of order $4$, then $\bal$ is collinear
with some $2\blm_{i}$ by Lemma \ref{lem:geometry}(b),
and thus $2\bal=2\blm_{i}$; i.e.,
$\bal=\blm_{i}+\bar{\mu}$ with $2\bar{\mu}=0$ \ Define
\[
M_{i}=\{\mu \in \Lm :\blm_{i}+\bar{\mu}\text{ is a point and }2\bar{\mu}=0\},
\]
so the points of order $4$ are all $\bal$ with
$\al\in\lm_{i}+M_{i}$ for some $1\leq i\leq3$. Clearly,
\[2M_{i}\subset\Lm_{-}\subset M_{i}.\]
If $\al\in\lm_{i}+M_{i}$ and $\beta\in\lm_{j}+M_{j}$ with
$i\neq j$, then $2\bal=2\blm_{i}\neq2\blm_{j}%
=2\bbe$, so $\bbe$ must be collinear with $\bal$ on
$[\bal,\bal,2\bal]$. This shows that if $\bar
{\al}_{i}\in\blm_{i}+M_{i}$, $i=1,2,3$, and $\al_{1}+\al
_{2}+\al_{3}=0$, then $[\bal_{1},\bal_{2},\bal%
_{3}]$ is a line. In particular, $[\blm_{1},\blm_{2}%
,\blm_{3}]$ is a line. If $\bde$ is a point of order $2$,
then $\bde$ is collinear with some $\blm_{i}$, so $\bar
{\delta}=2\blm_{i}$. Thus, the points are $2\blm_{i}$ and
all $\bal_{i}$ with $\al_{i}\in\lm_{i}+M_{i}$.

Next, if
$[\bal,\bbe,\bgm]$ is a line with $\al,\beta\in\lm_{i}+M_{i}$, then
$2\bgm=-2(\bal+\bbe%
)=-4\blm_{i}=0$. \ Since $\bgm$ has order $2$, we see that
$\bgm=2\bal=2\blm_{i}$. \ Also, $\bal%
+\bbe=-\bgm=2\bal$ shows that $\bbe=\bar{\al}$.
Thus, the only lines are $[2\blm_{1},2\blm_{2},2\blm_{3}]$,
$[\bal_{i},\bal_{i},2\blm_{i}]$ and $[\bal_{1},\bal_{2},\bal_{3}] $ with
$\al_{i}\in\lm_{i}+M_{i}$ and $\al_{1}+\al_{2}+\al_{3}=0$.

If $\al\in\lm_{i}+M_{i}$ and $\beta\in\lm_{j}+M_{j}$ with $i\neq
j$, then $-(\al+\beta)\in\lm_{k}+M_{k}$ where $\{i,j,k\}=\{1,2,3\}$.
Thus, $M_{i}+M_{j}\subset-M_{k}$. \ Since $0\in M_{j}$, we see $M_{i}%
\subset-M_{k}$ so $M_{i}=-M_{k}=M_{j}$ is a subgroup $M$.

Since $\bar{\Lm}=\Lm/\Lm_{-}$ is generated by the points of
$\cI$ (by (ST3)), we see $\bar{\Lm}=\left\langle \blm_{1}%
,\blm_{2},\bar{M}\right\rangle $. Now $2M\subset\Lm_{-}$ and
$2\blm_{2}\neq2\blm_{1}$ shows that $v_{i}=2\bar{\lm
}_{i}$, $i=1,2,$ is a basis for $2(\Lm/\Lm_{-})$. Thus, we may apply
Lemma \ref{lem:basis} to $\Gm\subset\Lm_{-}\subset\Lm$ to
rechoose $\lm_1,\lm_2$ and get a basis $\lm_1,\lm_2,\eta_0,\dots,\eta_n$
for $\Lm$ with $n\ge 0$ satisfying (i) and (ii).

It remains to show that
\begin{equation}
\label{eq:basisM}
M = \langle 2\lm_1,2\lm_2,\eta_0,\dots,\eta_n\rangle.
\end{equation}
Since $-\blm_{i}$ is a point, we have $2\lm_{i}\in M$.  Set
$N=\left\langle \eta_{0},\eta_{1},\ldots,\eta_{n}\right\rangle $.  If $\mu\in
M$, write $\mu=l_{1}\lm_{1}+l_{2}\lm_{2}+\eta$ with $\eta\in N$.
Since $2\mu\in\Lm_{-}$, we see that $l_{i}$ is even, so $l_{i}\lm_{i}, \eta\in M$.
Thus, $M=\left\langle 2\lm_{1},2\lm
_{2},M\cap N\right\rangle $. On the other hand, $\Lm_{-}\subset M$ and
$\bar{\Lm}=\left\langle \blm_{1},\blm_{2},\bar
{M}\right\rangle $ show $\Lm=\left\langle \lm_{1},\lm
_{2},M\right\rangle =\left\langle \lm_{1},\lm_{2},M\cap N\right\rangle
$. \ Thus, $M\cap N=N$ giving \eqref{eq:basisM}.
\end{proof}

\begin{theorem}
\label{thm:IIIc} Suppose that $\cA$ is a structurable $\Lm$-torus
of class III(c). Then $\cA$ is isograded isomorphic to either
\[
\cA(\cH(\mathcal{Q}_{3}))\otimes\mathcal{P}(r)\quad\text{ or,
if $-1\in(F^{\times})^{2}$,}\quad\cA(\cH(\mathcal{O}%
_{3}))\otimes\mathcal{P}(r)
\]
for some $r\geq0$. (See Theorem \ref{thm:constructIIIc} and Definition
\ref{def:constructIIIc}.)
\end{theorem}

\begin{proof}
Recall  that we are using the notation from
\S\ref{sec:classes}
and \S\ref{sec:geometry}.  We make one exception to this and denote the hermitian form and
cubic form defined in \S\ref{sec:geometry} by $\th$ and $\tN$ respectively,
rather that $h$ and $N$.  Thus, by Proposition \ref{prop:ClassIIIcubic},
$\th : \cW \times \cW \to \cE$
is a $\Lm$-graded hermitian form, $\tN : \cW \to \cE$ is a $\Lm$-graded cubic form,
$(\th,\tN)$ satisfies the adjoint identity, and
\[\cA = \cA(\th,\tN).\]
as $\Lm$-graded algebras.

Assume now that we have chosen a basis  $\lm_1,\lm_2,\eta_0,\dots,\eta_n$ for
$\Lm$ satisfying the properties in Lemma \ref{lem:IIIcgeometry}, and, as in that lemma, we set
$\lm_3 = -\lm_1-\lm_2$ and
$M=\langle2\lm_{1},2\lm_{2},\eta_{0},\eta_{1},\ldots,\eta_{n}\rangle$.
We can rearrange this basis for $\Lm$ so
that $k_{i}=2$ for $1\leq i\leq n^{\prime}$ and $k_{i}=1$ for $n^{\prime
}<i\leq n$. Let $\Lm^{\prime}=\left\langle \lambda_{1},\lambda_{2}%
,\eta_{0},\eta_{1},\ldots,\eta_{n^{\prime}}\right\rangle $ and $\Lm
^{\prime\prime}=\left\langle \eta_{n^{\prime}+1},\ldots,\eta_{n}\right\rangle
\subset\Gamma$. Thus, by Lemma \ref{lem:tensor}, $\cA\simeq_{\Lm
}\cA^{\Lm^{\prime}}\otimes\cA^{\Lm^{\prime\prime}%
}\simeq_{ig}\cA^{\Lm^{\prime}}\otimes\mathcal{P}(r)$, where
$r=n-n^{\prime}$. Replacing $\cA$ by $\cA^{\Lm^{\prime}}$,
we can assume that $k_{i}=2$ for $1\leq i\leq n$. Thus, $\Lm$ has basis
$\lambda_{1},\lambda_{2},\eta_{0},\eta_{1},\dots,\eta_{n}$,
\begin{gather}
M=\langle2\lambda_{1},2\lambda_{2},\eta_{0},\eta_{1},\dots,\eta_{n}%
\rangle,\quad\Lm_{-}=\langle4\lambda_{1},4\lambda_{2},k_{0}\eta_{0}%
,2\eta_{1},\dots,2\eta_{n}\rangle,\label{eq:ML}\\
\Gamma=\langle4\lambda_{1},4\lambda_{2},2k_{0}\eta_{0},2\eta_{1},\dots
,2\eta_{n}\rangle,\label{eq:G}%
\end{gather}
where $n\geq0$ and $k_{0}=1$ or $2$.

Note    that since $[2\bar\lm_i,2\bar\lm_j,2\bar\lm_k]$ and
$[-\bar\lm_i,-\bar\lm_i,2\bar\lm_i]$ are lines by Lemma \ref{lem:IIIcgeometry}(iii),
we have
\begin{equation}
\label{eq:lines}
\cW^{-2\lm_i} \diamond \cW^{-2\lm_j} \ne 0 \andd \cW^{-\lm_i}\diamond \cW^{-\lm_i} \ne 0
\end{equation}
for  $\set{i,j,k} = \set{1,2,3}$.

Choose $0\neq u_{i}\in\cW^{-\lambda_{i}}$ for $i=1,2$.  Let
$f_{i}=-u_{i}^{\natural}\in\cW^{-2\lambda_{i}}$ for $i=1,2$.
By \eqref{eq:lines}, we have $f_1\diamond f_2 \ne 0$, and so we may choose
$f_{3}\in\cW^{-2\lambda_{3}}$ with $\tilde{h}(f_{3},f_{1}\diamond
f_{2})=1$. Let $e_{i}=f_{j}\diamond f_{k}\in \cW^{2\lm_i}$
for $\{i,j,k\}=\{1,2,3\}$ in which case we have  $h(f_i,e_i) = 1$.
\ Since $4\lambda_{i}\in\Lm_{-}$, we have $f_{i}^{\natural}=0$. \ Thus,
(ADJ2) gives%
\[
e_{i}\diamond e_{j}=(f_{j}\diamond f_{k})\diamond(f_{i}\diamond f_{k}%
)=\tilde{h}(f_{k},f_{i}\diamond f_{j})f_{k}=f_{k}%
\]
so
\[
\tilde{h}(e_{1},e_{2}\diamond e_{3})=\tilde{h}(f_{2}\diamond f_{3}%
,f_{1})=1^{\ast}=1.
\]
Since $\cW^{2\lambda_{i}+\Lm_{-}}=\mathcal{E}e_{i}$, the description of points
in Lemma \ref{lem:IIIcgeometry}(iii) shows that
\[
\mathcal{W=E}e_{1}\oplus\mathcal{E}e_{2}\oplus\mathcal{E}e_{3}\oplus
\cW_{1}\oplus\cW_{2}\oplus\cW_{3}
\]
where $\cW_{i}:=\cW^{-\lambda_{i}+M}$.
Moreover, since $\tilde h(\cW^{\al_1},\cW^{\al_2},\cW^{\al_3}) = 0$
unless $[\bar\al_1,\bar\al_2,\bar\al_3]$ is a line,
the description of lines in Lemma \ref{lem:IIIcgeometry}(iii) shows that, if
$x={\sum_{i}}a_{i}e_{i}+{\sum_{i}}x_{i}$
with $a_{i}\in\mathcal{E}$ and $x_{i}\in\cW_{i}$, then
\[
\tilde{N}(x)=\frac{1}{6}\tilde{h}(x,x\diamond x)=a_{1}a_{2}a_{3}+%
{\sum_{i}}
a_{i}\tilde{h}(e_{i},x_{i}^{\natural})+\tilde{h}(x_{1},x_{2}\diamond
x_{3})\text{.}%
\]
Since $e_{i}\diamond x_{j}=0$ for $i\neq j$, (ADJ3) yields%
\begin{align}
(e_{j}\diamond x_{j})\diamond(e_{k}\diamond x_{k})  & =-(e_{j}\diamond
e_{k})\diamond(x_{j}\diamond x_{k}),\label{eq:ex1}\\
(e_{j}\diamond e_{k})\diamond(e_{i}\diamond x_{i})  & =\tilde{h}(e_{i}%
,e_{j}\diamond e_{k})x_{i}=x_{i}.\label{eq:ex2}%
\end{align}

We define $\ph:\cW\rightarrow\cW$ by $\ph(ae_{i})=a^{\ast
}f_{i}$ for $a\in\mathcal{E}$ and $\ph(x_{i})=-e_{i}\diamond x_{i}$ for
$x_{i}\in\cW_{i}$. \ Since $\cW^{2\lambda_{i}+\Lm_{-}%
}=\cW^{-2\lambda_{i}+\Lm_{-}}=\mathcal{E}f_{i}$ and
$\cW^{-\lambda_{i}+M}=\cW^{\lambda_{i}+M}=\cW_{i}$, we see $\ph$
is a $\ast$-semilinear vector space isomorphism by (\ref{eq:ex2}). We next
show that
\begin{equation}
\label{eq:tNpreserved}
\tilde{N}(\ph(x))=\tilde{N}(x)^{\ast}.
\end{equation}
for $x\in \cW$.
Clearly
\[
\tilde{h}(a_{1}^{\ast}f_{1},a_{2}^{\ast}f_{2}\diamond a_{3}^{\ast}%
f_{3})=(a_{1}a_{2}a_{3})^{\ast}.
\]
Since $(e_{i}\diamond x_{i})^{\natural}=\tilde{h}(e_{i},x_{i}^{\natural})e_{i}$
by (ADJ4), we also have
\[
\tilde{h}(a_{i}^{\ast}f_{i},(-e_{i}\diamond x_{i})^{\natural})
 =\tilde{h}(a_{i}^{\ast}f_{i},\tilde{h}(e_{i},x_{i}^{\natural})e_{i})
=(a_{i}\tilde
{h}(e_{i},x_{i}^{\natural}))^{\ast}.
\]
Finally using (\ref{eq:ex1}) and (\ref{eq:ex2}), we have%
\begin{align*}
\tilde{h}(-e_{1}\diamond x_{1},&(-e_{2}\diamond x_{2})\diamond(-e_{3}\diamond
x_{3}))   =\tilde{h}(e_{1}\diamond x_{1},(e_{2}\diamond e_{3})\diamond
(x_{2}\diamond x_{3}))\\
& =\tilde{h}(x_{2}\diamond x_{3},(e_{2}\diamond e_{3})\diamond(e_{1}\diamond
x_{1}))  =\tilde{h}(x_{2}\diamond x_{3},x_{1})=(\tilde{h}(x_{1},x_{2}\diamond
x_{3}))^{\ast}.
\end{align*}
This shows \eqref{eq:tNpreserved} as claimed.

Define $\tilde{T}(x,y):=\tilde{h}(x,\ph(y))$. \ If $y=%
{\sum_{i}}
b_{i}e_{i}+{\sum_{i}}
y_{i}$ with $b_{i}\in\mathcal{E}$ and $y_{i}\in\cW_{i}$,
then%
\[
\tilde{T}(x,y)  ={\sum_{i}}\tilde{h}(a_{i}e_{i},b_{i}^{\ast}f_{i})+{\sum_{i}}
\tilde{h}(x_{i},-e_{i}\diamond y_{i})
={\sum_{i}}a_{i}b_{i}-{\sum_{i}}\tilde{h}(e_{i},x_{i}\diamond y_{i}),
\]
so $\tilde{T}$ is symmetric. Now Lemma \ref{lem:psi1} (with $\psi=\ph^{-1}$)
shows that $\tilde{T}$ is nondegenerate and
$(\tilde{T},\tilde{N})$ satisfies the adjoint identity with
$x^{\#}=\ph^{-1}(x^{\natural})=\ph(x)^{\natural}$. For $x_{i}%
\in\cW_{i}$, set
\[n_{i}(x_{i})=-\tilde{h}(e_{i},x_{i}^{\natural})\in\mathcal{E}.\]
We now have
\begin{align*}
\tilde{T}(x,y)  & ={\sum_{i}}a_{i}b_{i}+{\sum_{i}}n_{i}(x_{i},y_{i}),\\
\tilde{N}(x)  & =a_{1}a_{2}a_{3}-
{\sum_{i}}
a_{i}n_{i}(x_{i})+\tilde{h}(x_{1},x_{2}\diamond x_{3}).
\end{align*}
Also  $n_{i}(u_{i})= -\tilde h(e_i,u_i^\natural)
= \tilde{h}(e_{i},f_{i})=1$ for $i=1,2$.  So we may apply Lemma \ref{lem:normH}
(with the
identification of $\mathcal{E}e_{1}\oplus\mathcal{E}e_{2}\oplus\mathcal{E}e_{3}$
with $\mathcal{E}^{3}$ and with $\tau(x_{1},x_{2},x_{3})=\tilde
{h}(x_{1},x_{2}\diamond x_{3})$) to get a
composition algebra $\cC$ over $\cE$ and $\cE$-linear isomorphisms $\eta_{i}:\cC\rightarrow\cW_{i}$
such that $\eta_i(1) = u_i$
(with $u_3 = u_1\times u_2$)
and such that the $\cE$-linear
isomorphism $\eta:\Herm\rightarrow\cW$ given by
\[\eta(\sum_i a_i[ii] + \sum_{(i,j,k)\cyc}x_i[jk])
= \sum_i a_i e_i + \eta_i(x_i),\]
satisfies
\begin{equation}
\label{eq:etapreserves}
\tilde T(\eta(x),\eta(y)) = T(x,y) \andd
\tilde N(\eta(x)) = T(x,y).
\end{equation}

We next use $\eta^{-1}$ to transfer the grading of
$\cW$ to a grading of $\Herm$, in which case $\eta$ is a graded map.
Then, by Lemma \ref{lem:gradH}(a),
$\cC$ has a $\Lm$-grading as an algebra with involution
over $\mathcal{E}$ and as an $\cE$-module
so that $\Herm$ has the grading given by~(\ref{eq:gradH}).

Next, since $\eta_{3}(\cC%
)=\cW_{3}=\cW^{-\lambda_{3}+M}$, we see that $\cC$ is
finely graded with support $M$. \ Also, $-n(\cC^{\mu}%
)=T(1[33],(\cC^{\mu}[12])^{\#})=\tilde{h}(e_{3},(\cW_{3}^{\mu
})^{\natural})\neq0$, so $\cC$ is an alternative $M$-torus with  involution.

We now know that $\cC$ can be built from $\mathcal{E}$ by iterating
the Cayley-Dickson process. To be more precise, suppose $M^{\prime}$ is a
subgroup of $\Lm$ such that $\Lm_{-}\subset M^{\prime}\subset M$,
$\mu\in M\setminus M^{\prime}$ and $0\neq x\in\cC^{\mu}$. Then,
$n(\cC^{M^{\prime}},x)\subset\mathcal{E}^{M^{\prime}+\mu}=0$. Using
this fact along with the fact that nonzero homogeneous elements of
$\cC$ are invertible, the standard argument (see for example
\cite[p.~164]{Mc}) shows that $\cC^{\langle M^{\prime},\;\mu\rangle
}=\cC^{M^{\prime}}\oplus x\cC^{M^{\prime}}$ is a graded
subalgebra of $\cC$ that can be identified with $\CD%
(\cC^{M^{\prime}},x^{2})$. If we do this process repeatedly, starting
with $M^{\prime}=\Lm_{-}$, we see that $\cC$ can be obtained from
$\mathcal{E}$ by a finite number of iterations of the Cayley-Dickson process.
Moreover, since $\cC$ is alternative, the number of iterations is
necessarily at most 3 \cite[p.~162]{Mc}.

Since $M/\Lm_{-}$ has exponent 2, it follows that $|M/\Lm_{-}%
|\leq2^{3}=8$. Thus, by \eqref{eq:ML}, we have one of the following cases:
\[
\text{(1) $n=0$ and $k_{0}=1$;\quad(2) $n=0$ and $k_{0}=2$; \quad(3) $n=1$ and
$k_{0}=1$.}%
\]
In all three cases, we choose $0\neq v_{i}\in\cC^{2\lambda_{i}}$ and
set $t_{i}=v_{i}^{2}\in\mathcal{E}^{4\lambda_{i}}$ for $i=1,2$.

In Case 1, we have $M=\langle2\lambda_{1},2\lambda_{2},\eta_{0}\rangle$ and
$\Lm_{-}=\langle4\lambda_{1},4\lambda_{2},\eta_{0}\rangle$. We choose
$0\neq s\in\mathcal{E}^{\eta_{0}}$. Then $\mathcal{E}=F[t_{1}^{\pm1}%
,t_{2}^{\pm1},s^{\pm1}]$ and, by the argument above, we may identify
$\cC=\CD(\mathcal{E},t_{1},t_{2})=\mathcal{Q}$ with
canonical generators $v_{1},v_{2}$ (writing $\sigma=\eta_{0}$ to conform with
the notation in  \S \ref{sec:constructionIIIc}).

In Case 2, we have $M=\langle2\lambda_{1},2\lambda_{2},\eta_{0}\rangle$ and
$\Lm_{-}=\langle4\lambda_{1},4\lambda_{2},2\eta_{0}\rangle$. We choose
$0\neq v\in\cC^{\eta_{0}}$ and set $s=v^{2}\in\mathcal{E}^{2\eta_{0}}$.
Then $\mathcal{E}=F[t_{1}^{\pm1},t_{2}^{\pm1},s^{\pm1}]$ and we may
identify $\cC=\CD(\mathcal{E},t_{1},t_{2},s)=\mathcal{O}$
with canonical generators $v_{1},v_{2},v$ (writing $\lambda=\eta_{0}$).

In Case 3, we have $M=\langle2\lambda_{1},2\lambda_{2},\eta_{0},\eta
_{1}\rangle$ and $\Lm_{-}=\langle4\lambda_{1},4\lambda_{2},\eta_{0}%
,2\eta_{1}\rangle$. We choose $0\neq s\in\mathcal{E}^{\eta_{0}}$ and $0\neq
v\in\cC^{\eta_{1}}$, and we set $r=v^{2}\in\mathcal{E}^{2\eta_{1}}$.
Then $\mathcal{E}=F[t_{1}^{\pm1},t_{2}^{\pm1},s^{\pm1},r^{\pm1}]$ and we may
identify $\cC=\CD(\mathcal{E},t_{1},t_{2},r)=\mathcal{O}^{\prime}$
with canonical generators $v_{1},v_{2},v$ (writing $\sigma=\eta_{0}$, $\lambda=\eta_{1}$).

One sees, using \eqref{eq:ML} and \eqref{eq:G}, that the involution $\ast$ on
$\mathcal{E}$ coincides in each case with the involution defined by
\eqref{eq:def*} in \S \ref{sec:constructionIIIc}.

If we now define $\psi : \Herm \to \Herm$ by $\psi=\eta^{-1}\ph^{-1}\eta$
and $h:\Herm\times \Herm\to \cE$ by
\begin{equation}
\label{eq:hdef}
h(x,y)=T(x,\psi(y))=\tilde{T}(\eta(x),\ph^{-1}\eta(y))=\tilde{h}(\eta
(x),\eta(y)),
\end{equation}
then $h$ is a graded hermitian form. Also, by \eqref{eq:tNpreserved} and
\eqref{eq:etapreserves}, $\psi$ is semi-norm preserving.  Therefore,
by Lemma \ref{lem:psi2},
Case 3 does not occur, and in Case 2 we have $-1\in(F^{\times})^{2}$.
Moreover, the map $\psi$ is determined, as in Lemma \ref{lem:psi2}(ii), by
some choice of $w_{i}\in\cC^{2\lambda_{i}}$ and some choice of $\iota$
with $\iota^{2}=-1$. Let $\cA(h,N)$ be the structurable torus
constructed as in Theorem \ref{thm:constructIIIc} using these choices. Then,
it follows using
\eqref{eq:etapreserves} and \eqref{eq:hdef} that the map
$(a,x)\mapsto(a,\eta(x))$ from $\cA(h,N)=\mathcal{E}\oplus
\Herm$ onto $\cA=\mathcal{E}\oplus\cW$
is an isomorphism of graded algebras with involution. Hence, by the final
statement in Theorem \ref{thm:constructIIIc}, we
have $\cA\simeq_{ig}\cA(\Herm)$.
\end{proof}


\begin{thebibliography}{AABGP}

\bibitem[A1]{A1} B.~Allison, \textit{A class of nonassociative algebras with
involution containing the class of Jordan algebras}, Math. Ann. \textbf{237}
(1978), 133--156. MR0507909 (81h:17003)

\bibitem[A2]{A2} B.~Allison,
\emph{Models of isotropic simple Lie algebras},
Comm. Algebra \textbf{7} (1979), no.~17, 1835--1875. MR0547712 (81d:17005)

\bibitem[AABGP]{AABGP}
B. Allison, S.~Azam, S.~Berman, Y.~Gao and A.~Pianzola,
\emph{Extended affine Lie algebras and their root systems},
Mem. Amer. Math. Soc.  126  (1997), no. 603, x+122 pp. MR1376741
(97i:17015)

\bibitem[AF]{AF} B.~Allison and J.~Faulkner, \textit{A Cayley-Dickson process
for a class of structurable algebras, }Trans. Amer. Math. Soc. \textbf{283}
(1984), 185--210. MR0735416 (85i:17001)

\bibitem[AFY]{AFY} B.~Allison, J. Faulkner, and Y.~Yoshii, \textit{Lie tori of
rank $1$}, Resenhas  6  (2004), no. 2-3, 99--109. MR2215972 (2006m:17021)

\bibitem[AG]{AG} B.~Allison and Y.~Gao, \emph{The root system
and the core of an extended affine Lie algebra}, Selecta
Math.  {\bf 7 } (2001), 149--212.  MR1860013 (2002g:17041)

\bibitem[ABG]{ABG} B.~Allison, Y. Gao and G. Benkart,
\emph{Lie tori of type BC$_r$, $r\ge 3$}, in preparation.

\bibitem[AH]{AH} B.~Allison and W.~Hein,
\textit{Isotopes of some nonassociative algebras with involution},
J. Algebra  69  (1981), no. 1, 120--142. MR0613862 (82k:17013)

\bibitem[AY]{AY}B.~Allison and Y.~Yoshii, \textit{Structurable tori and
extended affine Lie algebras of type } BC$_{1}$, J. Pure and Appl. Algebra
\textbf{184} (2003), 105-138. MR2004970 (2004h:17022)

\bibitem[BZ]{BZ} G.~Benkart and E.~Zelmanov,
\emph{Lie algebras graded by finite root systems and intersection matrix
algebras}, Invent. Math. \textbf{126} (1996), no.~1, 1--45. MR1408554 (97k:17044)

\bibitem[BGK]{BGK} S.~Berman, Y.~Gao and  Y.~Krylyuk, \emph{Quantum tori and
the structure of elliptic quasi-simple Lie algebras},
J. Funct. Anal. {\bf 135} (1996), 339--389.  MR1370607 (97b:17007)

\bibitem[BGKN]{BGKN} S.~Berman, Y.~Gao,  Y.~Krylyuk and E.~Neher, \emph{The
alternative torus and the structure of elliptic quasi-simple Lie
algebras of type $A_2$}, Trans. Amer. Math. Soc. {\bf 347} (1995),
no.~11, 4315--4363. MR1303115 (96b:17009)

\bibitem[D]{D}J.~Dieudonn\'{e}, \textit{La g\'{e}om\'{e}trie des groupes
classiques} (French), Troisi\`{e}me \'{e}dition. Ergebnisse der Mathematik und
ihrer Grenzgebiete, Band 5. Springer-Verlag, Berlin-New York, 1971. MR0310083
(46 \#9186)

\bibitem[F1]{F1}J.~Faulkner, \textit{Generalized quadrangles and cubic forms,}
Comm. Algebra \textbf{29} (2001), 4641--4653. MR1855115 (2002j:17039)

\bibitem[F2]{F2} J.~Faulkner, \textit{Lie tori of type BC$_2$ and structurable
quasitori}, preprint.

\bibitem[J1]{J1} N.~Jacobson, \textit{Triality and Lie algebras of type D$_4$,}
Rend. Circ. Mat. Palermo (2)  \textbf{13}  (1964) 129--153. MR0181705 (31 \#5932)

\bibitem[J2]{J2}N.~Jacobson,\textit{Structure and representations of Jordan algebras},
Amer. Math. Soc. Colloq. Publ. \textbf{39} American Mathematical
Society, Providence, R.I., 1968. MR0251099 (40 \#4330)

\bibitem[K1]{K1} I.~L.~Kantor, \emph{Some generalizations of Jordan algebras},
Trudy Sem. Vektor.~Tenzor.~Anal. {\bf 16} (1972), 407-499 (Russian). MR0321986 (48 \#351)

\bibitem[K2]{K2} I.~L.~Kantor, \emph{Models of exceptional Lie algebras},
Soviet Math. Dokl.~{\bf 14} (1973), 254-258. MR0349779 (50 \#2272)

\bibitem[Mc]{Mc} K.~McCrimmon,
{\textit A taste of Jordan algebras}, Springer, New York, 2004. MR2014924 (2004i:17001)

\bibitem[MY]{MY} J.~Morita and Y.~Yoshii, \textit{Locally extended affine Lie algebras},
J. Algebra 301 (2006), 59--81.  MR2230320 (2007c:17034)

\bibitem[Neeb]{Neeb} K.-H.~Neeb,
\emph{On a normal form of the rational quantum tori and their automorphism groups},
arXiv: math.RA/0511263.

\bibitem[N]{N} E.~Neher, \emph{Extended affine Lie algebras},
C. R. Math. Acad. Sci. Soc. R. Can.  {\bf 26}  (2004), no. 3, 90--96. MR2083842 (2005f:17024)


\bibitem[NY]{NY}E.~Neher and Y.~Yoshii, \textit{Derivations and invariant
forms of Jordan and alternative tori}, Trans. Amer. Math. Soc. \textbf{355}
(2003) 1079--1108. MR1938747 (2003i:16052)

\bibitem[PT]{PT}
S.~E.~Payne and J.A.~Thas, \textit{Finite generalized quadrangles},
Pitman, Boston, MA, 1984. MR0767454 (86a:51029)

\bibitem[S]{S} R.~Schafer, \textit{An introduction to nonassociative algebras,}
Pure and Applied Math., \textbf{22} Academic Press, New York-London. MR0210757
(35 \#1643)

\bibitem[Se]{Se} G.~B.~Seligman, \textit{Constructions of Lie algebras and their
modules}, Lecture Notes in Mathematics, 1300. Springer-Verlag, Berlin, 1988.
vi+190 pp. ISBN: 3-540-18973-4 MR0936841 (89h:17007)

\bibitem[Y1]{Y1}Y.~Yoshii, \textit{Coordinate algebras of extended affine Lie
algebras of type }$A_{1}$, J. Algebra \textbf{234} (2000),
128--168. MR1799481 (2001i:17031)

\bibitem[Y2]{Y2}Y.~Yoshii, \textit{Classification of quantum tori with
involution}, Canad. Math. Bull. \textbf{45} (2002), no. 4, 711--731. MR1941236 (2003i:16060).


\end{thebibliography}
\end{document}